\documentclass[twoside,11pt,reqno]{amsart}
\usepackage{amsmath,amssymb,amscd,mathrsfs}

\makeatletter

\newif\ifcenters@
\centers@false

\@settopoint\oddsidemargin
\@settopoint\marginparwidth
\@settopoint\evensidemargin
\@settopoint\topmargin
\newtheorem{Proposition}{Proposition}[section]
\newtheorem{Lemma}[Proposition]{Lemma}
\newtheorem{Theorem}[Proposition]{Theorem}
\newtheorem{Corollary}[Proposition]{Corollary}

\newtheorem{Remark}[Proposition]{Remark}

\newtheorem{Example}[Proposition]{Example}

\@addtoreset{equation}{section}

\def\op{\operatorname{op}}

\def\row{\operatorname{row}}
\def\Ad{\operatorname{Ad}}
\def\col{\operatorname{col}}
\def\rdet{\operatorname{rdet}}

\def\cdet{\operatorname{cdet}}
\def\id{\operatorname{id}}

\def\rt{{\operatorname{\scriptscriptstyle R}}}
\def\lt{{\operatorname{\scriptscriptstyle L}}}
\def\C{{\mathbb C}}

\def\Z{{\mathbb Z}}

\def\LL{{{\scriptscriptstyle\text{\rm{L}}}}}

\def\F{\mathrm{F}}
\def\L{\mathrm{L}}

\def\pr{{\operatorname{pr}}}
\def\diag{{\operatorname{diag}}}
\def\ad{\operatorname{ad}}
\def\gr{\operatorname{gr}}
\def\Mat{M}
\def\hom{{\operatorname{Hom}}}

\def\End{{\operatorname{End}}}

\def\sgn{{\operatorname{sgn}}}
\def\Mod{\operatorname{\text{-}Mod}}

\def\eps{{\varepsilon}}
\def\phi{{\varphi}}

\def\elly{\mathfrak{l}}

\begin{document}
\title[Shifted Yangians and finite $W$-algebras]{\boldmath 
Shifted Yangians and finite $W$-algebras}
\author{Jonathan Brundan and Alexander Kleshchev}
\address
{Department of Mathematics\\ University of Oregon\\
Eugene\\ OR~97403, USA}
\email{brundan@darkwing.uoregon.edu, klesh@math.uoregon.edu}
\thanks{Research partially supported by the NSF (grant no. DMS-0139019).}
\thanks{
{\em 2000 Subject Classification}: 17B37.}

\begin{abstract}
We give a presentation for the finite $W$-algebra
associated to a nilpotent matrix in the general linear
Lie algebra over $\C$. In the special case that
the nilpotent matrix consists of $n$ Jordan blocks each
of the same size $l$,
the presentation is that of the {\em Yangian of level $l$}
associated to $\mathfrak{gl}_n$,
as was first observed by Ragoucy and Sorba. In the general case,
we are lead to introduce some generalizations of the Yangian
which we call the {\em shifted Yangians}.
\end{abstract}

\maketitle

\section{Introduction}\label{sintro}

Let $\mathfrak{g}$ be a finite dimensional reductive Lie
algebra over $\C$ equipped with a non-degenerate
invariant symmetric bilinear form $(.,.)$. 
Pick a nilpotent element $e \in \mathfrak{g}$, i.e. an element which
acts nilpotently on every finite dimensional $\mathfrak{g}$-module.
A $\Z$-grading $\mathfrak{g} = \bigoplus_{j \in \Z}
\mathfrak{g}_j$ of $\mathfrak{g}$ is called
a {\em good grading} for $e$ if $e \in \mathfrak{g}_2$ and  
the linear map $$
\ad e: \mathfrak{g}_j \rightarrow \mathfrak{g}_{j+2}
$$
is injective for $j \leq -1$, surjective for $j \geq -1$.
This definition originates in \cite{KRW} in the study
of certain $W$-algebras 
defined from affine Lie algebras by quantum Hamiltonian reduction.
A complete classification of all good gradings of
simple Lie algebras up to conjugacy can be found in
\cite{EK}. 

Since
$\ad e:\mathfrak{g}_{-1} \rightarrow \mathfrak{g}_1$
is bijective, the skew-symmetric bilinear 
form $\langle.,.\rangle$ on $\mathfrak{g}_{-1}$
defined by $\langle x,y \rangle := ([x,y],e)$ is non-degenerate.
Pick a Lagrangian subspace $\elly$ of $\mathfrak{g}_{-1}$
with respect to the form $\langle.,.\rangle$ and set
$\mathfrak{m} := \elly \oplus \bigoplus_{j \leq -2} \mathfrak{g}_j$.
This is a nilpotent subalgebra of $\mathfrak{g}$, and 
the map $\chi:\mathfrak{m} \rightarrow \C, x \mapsto (x,e)$ 
defines a representation of $\mathfrak{m}$. Let $I_\chi$ denote the kernel
of the corresponding associative algebra homomorphism 
$U(\mathfrak{m}) \rightarrow \C$,
where $U(\mathfrak{m})$ denotes the universal enveloping algebra of
$\mathfrak{m}$. Let
\begin{equation*}
Q_\chi := 
U(\mathfrak{g}) \otimes_{U(\mathfrak{m})} \C_\chi
\cong
U(\mathfrak{g}) / U(\mathfrak{g}) I_\chi
\end{equation*}
denote the induced $\mathfrak{g}$-module, and 
consider the endomorphism algebra
\begin{equation*}
W(\chi) := \End_{U(\mathfrak{g})}(Q_\chi)^{\op}.
\end{equation*}
Following terminology used in the mathematical physics literature,
we refer to these algebras as {\em finite $W$-algebras}; see e.g. \cite{BT}.
Applying Frobenius reciprocity, it is often more convenient to view
$W(\chi)$ instead as the subspace
of $U(\mathfrak{g}) / U(\mathfrak{g}) I_\chi$ consisting of all
cosets $y + U(\mathfrak{g}) I_\chi$ such that
$[x,y] \in U(\mathfrak{g}) I_\chi$ for all $x \in \mathfrak{m}$.
In this realization, the algebra structure on
$W(\chi)$ is defined by the formula
$(y + U(\mathfrak{g}) I_\chi)(y' + U(\mathfrak{g}) I_\chi)
= yy' + U(\mathfrak{g}) I_\chi$ for $y,y' \in U(\mathfrak{g})$
such that $[x,y], [x,y'] \in U(\mathfrak{g})I_\chi$ for
all $x \in \mathfrak{m}$.

In the special case that our fixed good grading is {\em even},
i.e. $\mathfrak{g}_j =0$ for all odd $j$,
the algebras $W(\chi)$
were already well studied 
by the end of the 1970s 
by Lynch \cite{Ly}, generalizing work of Kostant \cite{K}
treating regular nilpotent elements. 
Of course in the even case, we have simply that
$\mathfrak{m} = \bigoplus_{j \leq -2} \mathfrak{g}_j$.
Letting $\mathfrak{p} := \bigoplus_{j \geq 0} \mathfrak{g}_j$, the
PBW theorem implies that
\begin{equation*}
U(\mathfrak{g}) = U(\mathfrak{p}) \oplus U(\mathfrak{g}) I_\chi.
\end{equation*}
The projection $\pr_\chi:U(\mathfrak{g})
\rightarrow U(\mathfrak{p})$ along this direct sum decomposition
induces an isomorphism 
$U(\mathfrak{g}) / U(\mathfrak{g}) I_\chi 
\stackrel{\sim}{\rightarrow}
U(\mathfrak{p})$ 
which leads to an easier definition of the algebra
$W(\chi)$ in the even case as 
a {\em subalgebra} of $U(\mathfrak{p})$.
To make this precise, 
define a twisted action of $\mathfrak{m}$ on $U(\mathfrak{p})$
by $x \cdot y := \pr_\chi ([x,y])$ for $x \in \mathfrak{m}$
and $y \in U(\mathfrak{p})$. Then $\pr_\chi$ induces an isomorphism
between $W(\chi)$
and subalgebra $U(\mathfrak{p})^{\mathfrak{m}}$ of $U(\mathfrak{p})$
consisting of all twisted $\mathfrak{m}$-invariants. 
This is the original definition used
by Kostant and Lynch.

The most important examples of good gradings arise as follows.
By the Jacobson-Morozov theorem, we can embed $e$
into an $\mathfrak{sl}_2$-triple $(e,h,f)$, so
$[e,f] = h, [h,e] = 2e$ and $[h,f] = -2f$.
Then the representation theory of $\mathfrak{sl}_2$ implies that 
the $\ad h$-eigenspace decomposition of $\mathfrak{g}$
is a good grading for $e$. We refer to a good grading obtained in this
way as a {\em Dynkin grading}. For a Dynkin grading, the module
$Q_\chi$ is a {\em generalized Gelfand-Graev representation}
in the sense of Kawanaka \cite{Kaw} and Moeglin \cite{Moe}.
Its endomorphism algebra $W(\chi)$
has been studied by Premet \cite{P} as an application
of results on Lie algebras in positive characteristic.
Subsequently, a more direct approach has been
given by Gan and Ginzburg \cite{GG} which we follow here. 

Returning to an arbitrary good grading, \cite[Lemma 1.1]{EK} 
shows that we can always embed the given element
$e \in \mathfrak{g}_2$ into 
an $\mathfrak{sl}_2$-triple
$(e,h,f)$ with 
$h \in \mathfrak{g}_0$ and $f \in \mathfrak{g}_{-2}$. 
Letting 
$\mathfrak{c}_{\mathfrak{g}}(f)$
denote the centralizer of $f$ and 
$\mathfrak{m}^\perp
:= \{y \in \mathfrak{g}\:|\:(x,y)=0 \hbox{ for all }x \in \mathfrak{m}\}$, 
the following 
crucial formula is proved as in \cite[$\S$2.3]{GG}, using
\cite[Theorem 1.4]{EK}:
\begin{equation*}
\mathfrak{m}^\perp = [\mathfrak{m}, e] \oplus \mathfrak{c}_{\mathfrak{g}}(f).
\end{equation*}
Remarkably, given this formula, all the arguments 
from \cite{GG} in the context of Dynkin gradings
extend absolutely unchanged to arbitrary good gradings.
Let us just state briefly here the analogue of \cite[Proposition 6.3]{P}
which we regard as the fundamental structure
theorem for $W(\chi)$; see also 
\cite[Theorem 4.1]{GG} and
\cite[Theorem 2.3]{Ly}.
Introduce the {\em Kazhdan filtration}
$\cdots \subseteq 
\F_d U(\mathfrak{g}) \subseteq \F_{d+1} U(\mathfrak{g}) \subseteq 
\cdots$ of $U(\mathfrak{g})$
by declaring that a generator $x \in \mathfrak{g}_j$ is of degree 
$(j+2)$, i.e. $\F_d U(\mathfrak{g})$ is the span of all
monomials $x_1 \cdots x_m$ for $m \geq 0$ and
$x_1 \in \mathfrak{g}_{j_1}, \dots, x_m \in \mathfrak{g}_{j_m}$
with $(j_1+2)+\cdots+(j_m+2) \leq d$. 
Viewing $W(\chi)$ as a subspace of the quotient
$U(\mathfrak{g}) / U(\mathfrak{g}) I_\chi$,
there is an induced Kazhdan filtration on $W(\chi)$;
we denote the associated graded Poisson algebra by $\gr W(\chi)$.
Recall also that the {\em Slodowy slice} through the nilpotent orbit
containing $e$ is the affine subspace $e + \mathfrak{c}_{\mathfrak{g}}(f)$;
see \cite{Slod}. 
It has a natural Poisson structure which may be defined
following \cite[$\S$3.2]{GG} 
as the Hamiltonian reduction of the Kirillov-Kostant 
Poisson structure on $\mathfrak{g}$. 
Now the basic fact is that there is a canonical
isomorphism \begin{equation*}
\nu:\gr W(\chi) 
\stackrel{\sim}{\rightarrow} \C[e+\mathfrak{c}_{\mathfrak{g}}(f)]
\end{equation*} 
of Poisson algebras; see \cite[$\S$4.4]{GG} for its precise definition.
Hence, $W(\chi)$ can be viewed as a quantization of the Slodowy slice
$e + \mathfrak{c}_{\mathfrak{g}}(f)$.
Moreover, up to canonical
isomorphism, the algebra $W(\chi)$ is independent of the
particular choice of the Lagrangian subspace $\elly$ of $\mathfrak{g}_{-1}$;
see \cite[$\S$5.5]{GG}.

Another fundamental result in this subject
is {\em Skryabin's theorem} proved in \cite{Skry} for
Dynkin gradings; see also \cite[Theorem 6.1]{GG}.
Again, Skryabin's proof extends unchanged to any
good grading. To state the result, 
let ${\mathscr W}(\chi)$ 
be the category of all $\mathfrak{g}$-modules on which
$(x-\chi(x))$ act locally nilpotently for all $x\in \mathfrak{m}$
(``generalized Whittaker modules''). 
If $M\in
{\mathscr W}(\chi)$ then the subspace
\begin{equation*}
M^\mathfrak{m}:=\{v\in M\:|\: (x-\chi(x))v=0\ \text{for all $x\in
\mathfrak{m}$}\}
\cong \hom_{U(\mathfrak{g})}(Q_\chi, M)
\end{equation*}
is a $W(\chi)$-module, hence 
$M\mapsto M^\mathfrak{m}$ is a functor $F$ 
from ${\mathscr W}(\chi)$ to the
category $W(\chi)\Mod$ of all left $W(\chi)$-modules.
Also, we have the functor $G:=Q_\chi \otimes_{W(\chi)} ?$ 
from
$W(\chi)\Mod$ to ${\mathscr W}(\chi)$.
{Skryabin's theorem} asserts 
that the functors ${F}$ and
${G}$ are quasi-inverse equivalences
between ${\mathscr W}(\chi)$ and $W(\chi)\Mod$.
Moreover, every $M\in {\mathscr W}(\chi)$ is an injective
$\mathfrak{m}$-module and $Q_\chi$ is a free right $W(\chi)$-module.
Somewhat weaker results in the even case can already be found in
the work of Kostant and Lynch; see 
for example \cite[Theorems 2.4, 4.1]{Ly}.

\vspace{2mm}

In the remainder of the article we study the algebras
$W(\chi)$ in the special case
that $e$ is a nilpotent matrix of Jordan type 
$p_1 \leq \cdots \leq p_n$ inside the Lie
algebra
$\mathfrak{g} = \mathfrak{gl}_N$
over $\C$, taking the bilinear form $(.,.)$ to be the usual
trace form. 
Our main result gives an explicit
set of generators and relations for the algebra $W(\chi)$.
One surprising 
consequence of our presentation is that in fact up to isomorphism
the algebras $W(\chi)$ only depend on the 
conjugacy class of $e$, i.e. the partition $\lambda = (p_n,p_{n-1},\dots,p_1)$
of $N$,
{\em not} on the particular choice of the good grading for $e$.

The classification of good gradings for $e$
is described in \cite[Theorem 4.2]{EK} in terms of certain 
diagrams called {\em pyramids}. A consequence
of this classification is that, in type $A$, 
it is sufficient in order to define {\em all} the algebras
$W(\chi)$ to restrict attention just to {even good gradings}.
More precisely, 
every good grading for $e$ is {\em split} in the sense that 
it is always possible to adjust the grading
to obtain a new, even good grading for $e$ with the property that 
the subalgebra $\mathfrak{m}$ defined from the new grading
coincides with the subalgebra $\mathfrak{m}$ defined from the
original grading for some particular choice of Lagrangian subspace
$\elly$. Hence the algebra $W(\chi)$ defined from the new
grading is {\em equal} to the algebra $W(\chi)$ 
defined from the original grading using this choice of $\elly$. 
(In the language of \cite[$\S$4]{EK}, the pyramid defining
such a new 
grading may be obtained from the pyramid defining
the original grading by shifting one place to the left
all the boxes whose first coordinates
are of different parity to the first coordinates of 
the boxes on the bottom row.)

So we may
assume without loss of generality
that the given good grading is {\em even}.
The even good gradings for $e$ are classified up to conjugacy
in \cite[Proposition 4.3]{EK}; see also \cite[Lemma 7.2]{Ly}.
Again we visualize the classification in terms of some
pyramids as explained in $\S$\ref{swhittaker} below.
We will explain the idea in this introduction just
with one example: the diagram
\vspace{-1mm}
$$
\begin{picture}(50, 45)%
\put(0,0){\line(1,0){60}}
\put(0,16){\line(1,0){60}}
\put(0,30){\line(1,0){45}}
\put(30,45){\line(1,0){15}}
\put(0,0){\line(0,1){30}}
\put(15,0){\line(0,1){30}}
\put(30,0){\line(0,1){45}}
\put(45,0){\line(0,1){45}}
\put(60,0){\line(0,1){15}}
\put(7,23){\makebox(0,0){1}}
\put(7,8){\makebox(0,0){2}}
\put(22,8){\makebox(0,0){4}}
\put(22,23){\makebox(0,0){3}}
\put(37,8){\makebox(0,0){7}}
\put(37,23){\makebox(0,0){6}}
\put(52,8){\makebox(0,0){8}}
\put(37,38){\makebox(0,0){5}}
\put(-16,22){\makebox(0,0){$\pi=$}}
\end{picture}
$$
is a pyramid for $\mathfrak{g} = \mathfrak{gl}_8$ of height $n=3$.
Numbering the bricks $1,\dots,8$
as indicated, the rows $1,\dots,3$ from top to bottom and the columns
$1,\dots,4$ from left to right,
we write $\row(i)$ and $\col(i)$ for the row 
and column numbers of the $i$th brick in the pyramid, respectively.
Denoting the $ij$-matrix unit by $e_{i,j}$, the
nilpotent matrix $e$ associated to $\pi$ is the matrix
$$
e = e_{1,3} + e_{3,6} + e_{2,4}+e_{4,7}+e_{7,8}
$$
defined by reading along the rows of the pyramid,
and the even good grading associated
to $\pi$ is defined by declaring that $e_{i,j}$ is of degree
$(\col(j)-\col(i))$;
actually, according to this definition, $e$ is of degree $1$ not $2$ since 
in the case of an even good grading
we prefer to divide all degrees by $2$.
Note that the Jordan type $(p_1,\dots,p_n)$ of 
$e$ in this example is $(1,3,4)$ (``row lengths'') and 
the parabolic subalgebra $\mathfrak{p} := \bigoplus_{j \geq 0}
\mathfrak{g}_j$ is of standard Levi shape $(2,2,3,1)$ (``column heights''). 
From the pyramid $\pi$ we also read off the
{\em level} $l := p_n$ and 
a certain {\em shift matrix}
$\sigma = (s_{i,j})_{1 \leq i,j \leq n}$
as explained in $\S$\ref{swhittaker}; in our example,
\begin{align*}
l=4,\qquad\qquad
\sigma &= \left(
\begin{array}{lll}
0&0&1\\
2&0&1\\
2&0&0
\end{array}
\right).
\end{align*}
{From now on} we will denote the finite $W$-algebra
$W(\chi)$ instead by $W(\pi)$;
recall it may be defined in the even case as the 
subalgebra $U(\mathfrak{p})^{\mathfrak{m}}$ of 
all twisted $\mathfrak{m}$-invariants in $U(\mathfrak{p})$.

Theorem~\ref{main} below asserts for any pyramid $\pi$ that 
$W(\pi)$ is isomorphic to the
{\em shifted Yangian $Y_{n,l}(\sigma)$ of level $l$}, namely,
the quotient of the {\em shifted Yangian} $Y_n(\sigma)$
by the two-sided ideal
generated by elements $\{D_1^{(r)}\}_{r > p_1}$.
Here, $Y_n(\sigma)$
denotes the algebra defined by generators
$\{D_i^{(r)}\}_{1 \leq i \leq n, r > 0},
\{E_i^{(r)}\}_{1 \leq i < n, r > s_{i,i+1}}$ and
$\{F_i^{(r)}\}_{1 \leq i < n, r > s_{i+1,i}}$
subject to the relations 
(\ref{r2})--(\ref{r13}) recorded below.
In the special case that $p_1 = \cdots = p_n$, i.e. 
all the Jordan
blocks of $e$ are of the same size $l$, the
pyramid $\pi$ is an $n \times l$ rectangle and the shift matrix
$\sigma$ is the zero matrix, the presentation is
a variation on
Drinfeld's presentation \cite{D3} for the {\em Yangian $Y_{n,l}$ 
of level $l$}
considered by Cherednik \cite{Ch0,Ch}.
Hence in this case, $W(\pi)$ is a quotient of the Yangian $Y_n$
associated to the Lie algebra $\mathfrak{gl}_n$,
as was first noticed by Ragoucy and Sorba \cite{RS}. 

The remainder of the article is organized as follows.
In $\S$\ref{syangian}, we define the shifted Yangian
$Y_n(\sigma)$ and prove a PBW theorem for it.
In $\S$\ref{sparabolic}, we introduce some more elaborate
parabolic presentations for $Y_n(\sigma)$ following
\cite{BK}. These are important because they
allow us in $\S$\ref{sbaby} 
to write down an explicit formula for the so-called baby 
comultiplications.
In $\S$\ref{scan}, we introduce the canonical filtration of 
$Y_n(\sigma)$, which eventually turns out to correspond to the
Kazhdan filtration of $W(\pi)$.
In $\S$\ref{struncation}, we prove a PBW theorem
for the finitely generated quotients
$Y_{n,l}(\sigma)$ of $Y_n(\sigma)$.
Then we turn our attention back
to the finite $W$-algebras $W(\pi)$, beginning in
$\S$\ref{swhittaker} by
explaining the classification of even good gradings in terms of
pyramids. In $\S$\ref{sslice}, we recall the setup of
\cite{GG} in our special case in some detail.
The most important section of the article is
$\S$\ref{sinvariants}, where we write down explicit
formulae for elements of $U(\mathfrak{p})$ which eventually
turn out to be precisely the generators $D_i^{(r)}$,
$E_i^{(r)}$ and $F_i^{(r)}$ of $W(\pi)$ that we are after.
The main theorem is then proved by induction 
in $\S$\ref{smain}, the key tool for the induction step
being the baby comultiplications. 
In $\S$\ref{sgrownup} we discuss more general comultiplications. 
Finally, $\S$\ref{sspecial} gives a much simpler
and more direct proof of the main theorem in the
special case $p_1=\cdots=p_n$, using a different 
description of the generators of $W(\pi)$ in this special case
which is closely related to the Capelli determinant.

We finally note that one can 
replace the ground field $\C$ used throughout the article 
with an arbitrary algebraically closed
field of characteristic zero.
In a subsequent article \cite{BK2}, we will use the
presentation for $W(\pi)$ obtained here to study its highest weight
representation theory. 
In particular we will explain the sense in which
finite dimensional (rational) representations of $W(\pi)$\
categorify the polynomial representation of $GL_{\infty}$ 
parametrized by the partition $\lambda$.

\vspace{2mm}
\noindent
{\em Acknowledgements.}
We are grateful to Alexander Premet for illuminating discussions
with the second author
which triggered this work, Pavel Etingof for a useful
discussion of \cite{RS}, and Gerhard Roehrle for showing us the paper
\cite{EK}.
The first author also thanks Jeremy Rickard,
Joseph Chuang and
the Department of Mathematics at the University of Bristol
for their hospitality whilst this paper was written up.

\section{The shifted Yangian}\label{syangian}

Fix $n \geq 1$ 
and a matrix $\sigma = (s_{i,j})_{1 \leq i,j \leq n}$ of
non-negative integers (``shifts'') such that
\begin{equation}\label{shiftcon}
s_{i,j} + s_{j,k} = s_{i,k}
\end{equation}
whenever $|i-j|+|j-k|=|i-k|$.
Note this means that $s_{1,1}=\cdots=s_{n,n} = 0$, and the
matrix $\sigma$ is completely determined by the 
upper diagonal entries
$s_{1,2}, s_{2,3},\dots,s_{n-1,n}$
and the lower diagonal entries
$s_{2,1}, s_{3,2},\dots,s_{n,n-1}$.

The {\em shifted Yangian} associated to the matrix $\sigma$
is the algebra $Y_n(\sigma)$ over $\C$
defined by generators 
\begin{align*}
&\{D_i^{(r)}\}_{1 \leq i \leq n, r > 0},\\
&\{E_i^{(r)}\}_{1 \leq i < n, r > s_{i,i+1}},\\
&\{F_i^{(r)}\}_{1 \leq i < n, r > s_{i+1,i}}
\end{align*}
subject to 
certain relations.
In order to write down these relations,
let 
\begin{equation}\label{didef}
D_i(u) := \sum_{r \geq 0} D_i^{(r)} u^{-r}
\in Y_n(\sigma)[[u^{-1}]]
\end{equation}
where $D_i^{(0)} := 1$, and then
define some  new elements $\widetilde{D}_i^{(r)}$ of $Y_n(\sigma)$
from the equation
\begin{equation}
\widetilde{D}_i(u) =
\sum_{r \geq 0} \widetilde{D}_i^{(r)}u^{-r}
:=-D_i(u)^{-1}.
\end{equation}
For example, $\widetilde{D}_i^{(0)} = -1, 
\widetilde{D}_i^{(1)} = D_i^{(1)},
\widetilde{D}_i^{(2)} = D_i^{(2)}-D_i^{(1)} D_i^{(1)},\dots$.
With this notation, the relations are as follows.
\begin{align}
[D_i^{(r)}, D_j^{(s)}] &=  0,\label{r2}\\
[E_i^{(r)}, F_j^{(s)}] &= \delta_{i,j} 
\sum_{t=0}^{r+s-1} \widetilde D_{i}^{(t)}D_{i+1}^{(r+s-1-t)} ,\label{r3}
\end{align}\begin{align}
[D_i^{(r)}, E_j^{(s)}] &= (\delta_{i,j}-\delta_{i,j+1})
\sum_{t=0}^{r-1} D_i^{(t)} E_j^{(r+s-1-t)},\label{r4}\\
[D_i^{(r)}, F_j^{(s)}] &= (\delta_{i,j+1}-\delta_{i,j})
\sum_{t=0}^{r-1} F_j^{(r+s-1-t)}D_i^{(t)} ,\label{r5}\end{align}\begin{align}
[E_i^{(r)}, E_i^{(s)}] &=
\sum_{t=s_{i,i+1}+1}^{s-1} E_i^{(t)} E_i^{(r+s-1-t)}
-\sum_{t=s_{i,i+1}+1}^{r-1} E_i^{(t)} E_i^{(r+s-1-t)},\label{r6}\\
[F_i^{(r)}, F_i^{(s)}] &=
\sum_{t=s_{i+1,i}+1}^{r-1} F_i^{(r+s-1-t)} F_i^{(t)}-
\sum_{t=s_{i+1,i}+1}^{s-1} 
F_i^{(r+s-1-t)} F_i^{(t)},\label{r7}\end{align}\begin{align}
 [E_i^{(r)}, E_{i+1}^{(s+1)}]- [E_i^{(r+1)}, E_{i+1}^{(s)}] &=
-E_i^{(r)} E_{i+1}^{(s)},\label{r8}\\
[F_i^{(r+1)}, F_{i+1}^{(s)}] - [F_i^{(r)}, F_{i+1}^{(s+1)}] &=
 -F_{i+1}^{(s)} F_i^{(r)},\label{r9}\end{align}\begin{align}
[E_i^{(r)}, E_j^{(s)}] &= 0 \qquad\text{if }|i-j|> 1,\label{r10}\\
[F_i^{(r)}, F_j^{(s)}] &= 0 \qquad\text{if }|i-j|> 1,\label{r11}\\
[E_i^{(r)}, [E_i^{(s)}, E_j^{(t)}]] + 
[E_i^{(s)}, [E_i^{(r)}, E_j^{(t)}]] &= 0 \qquad\text{if }|i-j|=1,\label{r12}\\
[F_i^{(r)}, [F_i^{(s)}, F_j^{(t)}]] + 
[F_i^{(s)}, [F_i^{(r)}, F_j^{(t)}]] &= 0 \qquad\text{if }|i-j|=1,\label{r13}
\end{align}
for all admissible $r,s,t,i,j$.
(For an example of what we mean by ``admissible'' here, the relation (\ref{r8}) should be understood to hold for all
$i=1,\dots,n-2$, $r > s_{i,i+1}$ and $s > s_{i+1,i+2}$.)

\iffalse
\begin{Remark}\rm
The relations (\ref{r6}) are equivalent to
\begin{align}
[E_i^{(r)}, E_i^{(s+1)}] &- [E_i^{(r+1)}, E_i^{(s)}] =
E_i^{(r)} E_i^{(s)} + E_i^{(s)} E_i^{(r)},\label{r6p}\\\intertext{for all 
$i=1,\dots,n-1$ and
$r,s > s_{i,i+1}$, 
and the relations (\ref{r7}) are equivalent to}
[F_i^{(r+1)}, F_i^{(s)}] &- [F_i^{(r)}, F_i^{(s+1)}] =
F_i^{(r)} F_i^{(s)} + F_i^{(s)} F_i^{(r)},\label{r7p}
\end{align}
for all $i=1,\dots,n-1$ and $r,s > s_{i+1,i}$.
\end{Remark}
\fi

If the matrix $\sigma$ is the zero matrix,
we denote $Y_n(\sigma)$ simply by $Y_n$.
In this special case, the above presentation
is a variation on Drinfeld's presentation \cite{D3}
for the usual Yangian $Y(\mathfrak{gl}_n)$ 
associated to the Lie algebra 
$\mathfrak{gl}_n$; see 
\cite[Remark 5.12]{BK}.
We will prove in Corollary~\ref{injcor} below that
the map sending
the generators of $Y_n(\sigma)$ to the elements with the
same name in $Y_n$ is an injective algebra homomorphism. Given this fact,
the algebra $Y_n(\sigma)$ is canonically a {\em subalgebra}
of the usual Yangian $Y_n$. 
By the relations, there is an anti-automorphism 
$\tau:Y_n \rightarrow Y_n$ of order $2$ defined by
\begin{equation}\label{taudef}
\tau(D_i^{(r)}) = D_i^{(r)},\quad
\tau(E_i^{(r)}) = F_i^{(r)},\quad
\tau(F_i^{(r)}) = E_i^{(r)}.
\end{equation}
This obviously interchanges the two subalgebras
$Y_n(\sigma)$ and $Y_n(\sigma^t)$ of $Y_n$, where
$\sigma^t$ denotes the transpose of the matrix $\sigma$.
Hence $\tau$ also induces an anti-isomorphism
$\tau:Y_n(\sigma) \rightarrow Y_n(\sigma^t)$.
\iffalse
Another natural automorphism of $Y_n$ is
$\gamma:Y_n \rightarrow Y_n$ defined on
generators by
\begin{align}\label{gammadef1}
\gamma(D_i^{(r)}) &= (-1)^{(r-1)}\widetilde D_{n+1-i}^{(r)},\\
\gamma(E_i^{(r)}) &= (-1)^{(r-1)}F_{n-i}^{(r)},\\
\gamma(F_i^{(r)}) &= (-1)^{(r-1)}E_{n-i}^{(r)}.
\end{align}
\fi
Suppose instead that
$\dot\sigma = (\dot s_{i,j})_{1 \leq i,j \leq n}$ is another shift matrix
satisfying (\ref{shiftcon}),
such that $\dot s_{i,i+1}+\dot s_{i+1,i}= s_{i,i+1}+s_{i+1,i}$
for all $i=1,\dots,n-1$.
Another check of relations 
shows that there is a unique algebra isomorphism
$\iota:Y_n(\sigma) \rightarrow Y_n({\dot{\sigma}})$
defined by 
\begin{equation}\label{iotadef}
\iota(D_i^{(r)}) = \dot{D}_i^{(r)},\quad
\iota(E_i^{(r)}) = \dot{E}_i^{(r-s_{i,i+1}+\dot s_{i,i+1})},\quad
\iota(F_i^{(r)}) = \dot{F}_i^{(r-s_{i+1,i}+\dot s_{i+1,i})}.
\end{equation}
(Here and later on we denote the generators 
$D_i^{(r)}, E_i^{(r)}$ and $F_i^{(r)}$
of $Y_n(\dot\sigma)$
instead by $\dot{D}_i^{(r)}, \dot{E}_i^{(r)}$
and $\dot{F}_i^{(r)}$ to avoid potential confusion.)

Let us now prove as promised that the canonical map
$Y_n(\sigma) \rightarrow Y_n$ is injective.
Introduce the {\em loop filtration}
$\L_0 Y_n(\sigma) \subseteq \L_1 Y_n(\sigma) \subseteq \cdots$
by declaring that the generators
$D_i^{(r)}, E_i^{(r)}$ and $F_i^{(r)}$ of $Y_n(\sigma)$
are of degree $(r-1)$,
i.e. $\L_d Y_n(\sigma)$ is the span of all monomials in
the generators of total degree $\leq d$.
For $1 \leq i < j \leq n$ and $r > s_{i,j}$, define
elements $E_{i,j}^{(r)} \in Y_n(\sigma)$ recursively by
\begin{equation}
E_{i,i+1}^{(r)} := E_i^{(r)},
\qquad
E_{i,j}^{(r)} := [E_{i,j-1}^{(r - s_{j-1,j})}, E_{j-1}^{(s_{j-1,j}+1)}].
\end{equation}
Similarly, for $1 \leq i < j \leq n$ and $r > s_{j,i}$, define
elements $F_{i,j}^{(r)} \in Y_n(\sigma)$ by
\begin{equation}
F_{i,i+1}^{(r)} := F_i^{(r)},
\qquad
F_{i,j}^{(r)} := [F_{j-1}^{(s_{j,j-1}+1)},F_{i,j-1}^{(r - s_{j,j-1})}].
\end{equation}
It is easy to see that $E_{i,j}^{(r)}$ and $F_{i,j}^{(r)}$
belong to $\L_{r-1} Y_n(\sigma)$.
For $1 \leq i,j \leq n$ and $r \geq s_{i,j}$, define
\begin{equation}\label{newelts}
e_{i,j;r} := \left\{
\begin{array}{ll}
\gr^\LL_{r} D_{i}^{(r+1)} &\hbox{if $i = j$,}\\
\gr^\LL_{r} E_{i,j}^{(r+1)} &\hbox{if $i < j$,}\\
\gr^\LL_{r} F_{j,i}^{(r+1)} &\hbox{if $i > j$,}
\end{array}
\right.
\end{equation}
all elements of the associated graded algebra $\gr^\LL Y_n(\sigma)$.
Let $\mathfrak{gl}_n[t]$ denote the Lie algebra
$\mathfrak{gl}_n \otimes \C[t]$ on basis 
$\{e_{i,j}t^r\}_{1 \leq i,j \leq n, r \geq 0}$, viewed
as a graded Lie algebra so that $e_{i,j}t^r$ is of degree $r$.
In view of the assumption (\ref{shiftcon}), 
the vectors $\{e_{i,j}t^r\}_{1 \leq i,j \leq n, r \geq s_{i,j}}$
span a subalgebra of $\mathfrak{gl}_n[t]$
which we denote by $\mathfrak{gl}_n[t](\sigma)$ (the ``shifted loop algebra'').
The grading on $\mathfrak{gl}_n[t](\sigma)$
induces a grading on the universal enveloping algebra
$U(\mathfrak{gl}_n[t](\sigma))$.

\begin{Theorem}\label{pbw0}
There is an isomorphism
$\pi:U(\mathfrak{gl}_n[t](\sigma)) \rightarrow \gr^\LL
Y_n(\sigma)$ of graded algebras such that
$e_{i,j}t^r \mapsto e_{i,j;r}$ for each $1 \leq i,j \leq n$
and $r \geq s_{i,j}$.
\end{Theorem}

\begin{proof}
Using the relations like in the proof of \cite[Lemma 5.8]{BK}, one shows
for all $1 \leq h,i,j,k \leq n$, $r \geq s_{i,j}$ and
$s \geq s_{h,k}$ that
\begin{equation}\label{assgr}
[e_{i,j;r}, e_{h,k;s}]
= e_{i,k;r+s} \delta_{h,j}
- \delta_{i,k} e_{h,j;r+s},
\end{equation}
equality in $\gr^\LL Y_n(\sigma)$.
Hence there is a well-defined
surjection
$\pi:U(\mathfrak{gl}_n[t](\sigma)) \twoheadrightarrow
\gr^\LL Y_n(\sigma)$
mapping $e_{i,j}t^r \in \mathfrak{gl}_n[t](\sigma)$ to 
$e_{i,j;r} \in \gr^\LL Y_n(\sigma)$.

In the special case $\sigma = 0$,
the PBW theorem for the usual Yangian $Y_n$ 
implies that the ordered monomials in the elements
$\{e_{i,j;r}\}_{1 \leq i,j \leq n, r \geq 0}$ are linearly
independent in $\gr^\LL Y_n$; see the proof of
\cite[Lemma 5.10]{BK}. Hence $\pi$ is an isomorphism
in this case.

In general, the canonical map $Y_n(\sigma) \rightarrow Y_n$ 
is a homomorphism of filtered algebras, so
induces a map $\gr^\LL Y_n(\sigma) \rightarrow \gr^\LL Y_n$
which sends $e_{i,j;r} \in \gr^\LL Y_n(\sigma)$ to
$e_{i,j;r} \in \gr^\LL Y_n$ (despite the fact that
it does {\em not} in general 
send $E_{i,j}^{(r+1)}, F_{i,j}^{(r+1)} \in Y_n(\sigma)$ 
to $E_{i,j}^{(r+1)},F_{i,j}^{(r+1)} \in Y_n$ if $j-i>1$). 
So the previous paragraph
implies that the 
ordered monomials in the elements
$\{e_{i,j;r}\}_{1 \leq i,j \leq n, r \geq s_{i,j}}$ are linearly
independent in $\gr^\LL Y_n(\sigma)$ too.
Hence $\pi$ is an isomorphism in general.
\end{proof}

\begin{Corollary}\label{injcor}
The canonical map $Y_n(\sigma) \rightarrow Y_n$ is injective.
\end{Corollary}

\begin{proof}
We saw in the proof of Theorem~\ref{pbw0} that the 
canonical map $Y_n(\sigma) \rightarrow Y_n$ is filtered 
and the associated graded
map $\gr^\LL Y_n(\sigma) \rightarrow \gr^\LL Y_n$
is injective. 
\end{proof}

The presentation of $Y_n(\sigma)$
is adapted to the natural {\em triangular decomposition}
of this algebra. Let $Y_{(1^n)}$ denote the subalgebra
of $Y_n(\sigma)$ generated by all the 
$D_i^{(r)}$'s,
let $Y_{(1^n)}^+(\sigma)$ denote the subalgebra of $Y_n(\sigma)$
generated by all the $E_{i}^{(r)}$'s
and let $Y_{(1^n)}^-(\sigma)$ denote the subalgebra generated by the
all $F_{i}^{(r)}$'s.

\begin{Theorem}\label{pbw2}
\begin{itemize}
\item[(i)]
The monomials in the elements $\{D_i^{(r)}\}_{1 \leq i \leq n,
r > 0}$ taken in some fixed order form a basis for $Y_{(1^n)}$.
\item[(ii)]
The monomials in the elements $\{E_{i,j}^{(r)}\}_{1 \leq i < j 
\leq n,
r > s_{i,j}}$ taken in some fixed order form a basis for 
$Y^+_{(1^n)}(\sigma)$.
\item[(iii)]
The monomials in the elements $\{F_{i,j}^{(r)}\}_{1 \leq i < j 
\leq n,
r > s_{j,i}}$ taken in some fixed order form a basis for $Y^-_{(1^n)}(\sigma)$.
\item[(iv)]
The monomials in the union of the elements listed in (i)--(iii)
taken in some fixed order form a basis for $Y_n(\sigma)$.
\end{itemize}
\end{Theorem}

\begin{proof}
Part (iv) follows from Theorem~\ref{pbw0} and the PBW theorem
for $U(\mathfrak{gl}_n[t](\sigma))$.
The other parts are proved similarly, going back to (\ref{assgr}).
\end{proof}

\begin{Corollary}\label{triangular}
The natural multiplication map
$Y_{(1^n)}^-(\sigma) \otimes Y_{(1^n)} \otimes Y_{(1^n)}^+(\sigma)
\rightarrow Y_n(\sigma)$ is a vector space isomorphism.
\end{Corollary}

\ifcenters@
\begin{Remark}\rm
Let us describe the center $Z(Y_n(\sigma))$ of $Y_n(\sigma)$.
Recalling the notation (\ref{didef}), let
\begin{equation}\label{centelts}
C_n(u) =
\sum_{r \geq 0} C_n^{(r)} u^{-r}
:= D_1(u)D_2(u-1)\cdots D_n(u-n+1) \in Y_n(\sigma)[[u^{-1}]].
\end{equation}
Then, the elements $C_n^{(1)},C_n^{(2)},\dots$ are algebraically
independent and generate $Z(Y_n(\sigma))$.
Indeed, exploiting the embedding $Y_n(\sigma)\hookrightarrow Y_n$,
it is known by \cite[Theorem 2.13]{MNO} (or 
\cite[Theorem 7.2]{BK}) that the elements
$C_n^{(1)},C_n^{(2)},\dots$ are algebraically independent
and generate $Z(Y_n)$, so they certainly belong to
$Z(Y_n(\sigma))$. The fact that $Z(Y_n(\sigma))$ is no larger
than $Z(Y_n)$ may be proved by 
passing to the associated graded algebra 
$\gr^\LL Y_n(\sigma)$ and using a variation on the trick 
from the proof of \cite[Theorem 2.13]{MNO}.
We will outline a different argument in
Remark~\ref{missing} below.
\end{Remark}
\fi

\section{Parabolic presentations}\label{sparabolic}

We are going to need
various more elaborate presentations of the shifted Yangian which
are analogues of the parabolic presentations of \cite{BK}.
In order to explain the relationship between all these presentations,
we begin with an elementary remark about Gauss factorizations.

Let $T = (T_{i,j})_{1 \leq i,j \leq n}$ be an $n \times n$
matrix with entries in some ring 
such that the submatrices $(T_{i,j})_{1 \leq i,j \leq m}$
are invertible for all $m=1,\dots,n$.
Fix also a tuple $\nu = (\nu_1,\dots,\nu_m)$ of positive integers
summing to $n$, which we think of as the {\em shape} of
the standard Levi subgroup 
$\mathfrak{gl}_\nu
:=
\mathfrak{gl}_{\nu_1}\oplus\cdots\oplus \mathfrak{gl}_{\nu_m}$ 
of $\mathfrak{gl}_n$.
Working with $m \times m$ block matrices so that the 
$ab$-block is of size $\nu_a \times \nu_b$,
the matrix $T$ possesses a unique {Gauss} factorization
$T = {^\nu}F{^\nu}D{^\nu}E$
where ${^\nu}D$ is a block diagonal matrix,
${^\nu}E$ is a block upper unitriangular matrix,
and ${^\nu}F$ is a block lower unitriangular matrix:
$${^\nu}D =
\left(
\begin{array}{lllll}
{^\nu}D_1&0&&\cdots&0\\
0&{^\nu}D_2&&&0\\
\vdots&&\ddots&&\vdots\\
0&&&{^\nu}D_{m-1}&0\\
0&\cdots&&0&{^\nu}D_m
\end{array}
\right),
$$$$
{^\nu}E = 
\left(
\begin{array}{lllll}
I_{\nu_1}&{^\nu}E_1&*&\cdots&*\\
0&I_{\nu_2}&{^\nu}E_2&&\vdots\\
\vdots&&I_{\nu_3}&\ddots&*\\
0&&&\ddots&{^\nu}E_{m-1}\\
0&\cdots&&0&I_{\nu_m}
\end{array}\right),\:\:\:\:
{^\nu}F = 
\left(
\begin{array}{lllll}
I_{\nu_1}&0&&\cdots&0\\
{^\nu}F_1&I_{\nu_2}&&&0\\
*&{^\nu}F_2&I_{\nu_3}&&\vdots\\
\vdots&&\ddots&\ddots&0\\
*&\cdots&*&{^\nu}F_{m-1}&I_{\nu_m}
\end{array}\right).
$$
The diagonal blocks of ${^\nu}D$ define matrices
${^\nu}D_1,\dots,{^\nu}D_m$, the upper diagonal blocks of ${^\nu}E$
define matrices ${^\nu}E_1,\dots,{^\nu}E_{m-1}$, and the lower
diagonal blocks of ${^\nu}F$ define matrices ${^\nu}F_1,\dots,{^\nu}F_{m-1}$.
So ${^\nu}D_a$ is a $\nu_a \times \nu_a$ matrix,
${^\nu}E_a$ is a $\nu_a \times \nu_{a+1}$ matrix, and
${^\nu}F_a$ is a $\nu_{a+1}\times \nu_a$ matrix.
Now consider what happens when we split a block into two:
suppose that
$\nu_b = \alpha+\beta$ for some $1 \leq b \leq m$ and 
$\alpha,\beta \geq 1$,
and let $\mu = (\nu_1,\dots,\nu_{b-1},\alpha,\beta,\nu_{b+1},\dots,\nu_m)$.
The following lemma shows how to compute the matrices
${^\mu}D_1,\dots,{^\mu}D_{m+1},
{^\mu}E_1,\dots,{^\mu}E_m$ and ${^\mu}F_1,\dots,{^\mu}F_m$
just from knowledge of the matrices
${^\nu}D_1,\dots,{^\nu}D_{m},
{^\nu}E_1,\dots,{^\nu}E_{m-1}$ and ${^\nu}F_1,\dots,{^\nu}F_{m-1}$.

\begin{Lemma}\label{quasi}
In the above notation, define an $\alpha \times \alpha$ matrix
$A$, an $\alpha\times \beta$ matrix $B$,
a $\beta \times \alpha$ matrix $C$ and a $\beta \times \beta$ matrix $D$
from the equation
$$
{^\nu}D_b = \left(\begin{array}{ll}I_\alpha&0\\
C&I_\beta\end{array}
\right)\left(\begin{array}{ll}A&0\\
0&D\end{array}
\right)\left(\begin{array}{ll}I_\alpha&B\\
0&I_\beta\end{array}
\right).
$$
Then,
\begin{itemize}
\item[(i)] ${^\mu}D_a = {^\nu}D_a$ for $a < b$, ${^\mu}D_b = A$,
${^\mu}D_{b+1} = D$, and ${^\mu}D_c = {^\nu}D_{c-1}$ for $c > b+1$;
\item[(ii)] ${^\mu}E_a = {^\nu}E_a$ for $a < b-1$,
${^\mu}E_{b-1}$ is the submatrix consisting
of the first $\alpha$ columns
of
${^\nu}E_{b-1}$, ${^\mu}E_{b} = B$, ${^\mu}E_{b+1}$ is the submatrix consisting of the last $\beta$ rows  
of ${^\nu}E_b$, and ${^\mu}E_c = {^\nu}E_{c-1}$ for $c > b+1$;
\item[(iii)] ${^\mu}F_a = {^\nu}F_a$ for $a < b-1$,
${^\mu}F_{b-1}$ is the submatrix consisting of the first
$\alpha$ rows  of
${^\nu}F_{b-1}$, ${^\mu}F_{b} = C$, ${^\mu}F_{b+1}$ is the submatrix 
consisting of the last $\beta$ columns of
${^\nu}F_b$, and ${^\mu}F_c = {^\nu}F_{c-1}$ for $c > b+1$;
\end{itemize}
\end{Lemma}

\begin{proof}
Multiply matrices.
\end{proof}

Now let us briefly recall the parabolic presentations for the
Yangian $Y_n$ from \cite{BK}. 
As in \cite{MNO}, the Yangian can be defined
in terms of the RTT presentation as the algebra over $\C$
defined by generators $\{T_{i,j}^{(r)}\}_{1 \leq i,j \leq n, r >0}$
subject just to the relations
\begin{equation}\label{mr}
[T_{i,j}^{(r)}, T_{h,k}^{(s)}]
= \sum_{t=0}^{\min(r,s)-1}
\left( 
T_{i,k}^{(r+s-1-t)}T_{h,j}^{(t)} -
T_{i,k}^{(t)}T_{h,j}^{(r+s-1-t)}
\right)
\end{equation}
for every $1 \leq h,i,j,k \leq n$ and $r,s > 0$,
where $T_{i,j}^{(0)} := \delta_{i,j}$.
Let $T_{i,j}(u) := \sum_{r \geq 0} T_{i,j}^{(r)} u^{-r}$
and let $T(u)$ denote the $n \times n$ matrix
$(T_{i,j}(u))_{1 \leq i,j \leq n}$.
Given a shape $\nu = (\nu_1,\dots,\nu_m)$, 
consider the Gauss factorization
$T(u) = F(u) D(u) E(u)$ where $D(u)$ is a block diagonal matrix,
$E(u)$ is a block upper unitriangular matrix and $F(u)$
is a block lower unitriangular matrices, all block matrices
being of shape $\nu$ like before. 
We will from now on
omit any extra superscript $\nu$ since it should be clear
from the context which shape $\nu$ we have in mind.
The diagonal blocks of $D(u)$ define matrices
$D_1(u),\dots,D_m(u)$, the upper diagonal blocks of $E(u)$
define matrices $E_1(u),\dots, E_{m-1}(u)$, and the lower diagonal
blocks of $F(u)$ define matrices $F_1(u), \dots, F_{m-1}(u)$.
Thus $D_a(u) = (D_{a;i,j}(u))_{1 \leq i,j \leq \nu_a}$ is a
$\nu_a \times \nu_a$ matrix,
$E_a(u) = (E_{a;i,j}(u))_{1 \leq i \leq \nu_a, 1 \leq j \leq \nu_{a+1}}$
is a $\nu_a \times \nu_{a+1}$ matrix, and
$F_a(u) = (F_{a;i,j}(u))_{1 \leq i \leq \nu_{a+1}, 1 \leq j \leq \nu_a}$
is a $\nu_{a+1} \times \nu_a$ matrix.
Write
$$
D_{a;i,j}(u) = \sum_{r \geq 0} D_{a;i,j}^{(r)} u^{-r},\quad
E_{a;i,j}(u) = \sum_{r > 0} E_{a;i,j}^{(r)} u^{-r},\quad
F_{a;i,j}(u) = \sum_{r > 0} F_{a;i,j}^{(r)} u^{-r},
$$
thus defining elements $D_{a;i,j}^{(r)}, E_{a;i,j}^{(r)}$
and $F_{a;i,j}^{(r)}$ of $Y_n$, all dependent of course on the
fixed choice of $\nu$.
Now \cite[Theorem A]{BK} shows that the elements
\begin{align*}
&\{D_{a;i,j}^{(r)}\}_{1 \leq a \leq m, 1 \leq i,j \leq \nu_a,
r >0},\\
&\{E_{a;i,j}^{(r)}\}_{1 \leq a < m, 1 \leq i \leq \nu_a, 1 \leq j \leq \nu_{a+1},
r >0},\\
&\{F_{a;i,j}^{(r)}\}_{1 \leq a < m, 1 \leq i \leq \nu_{a+1}, 1 \leq j \leq \nu_{a},
r >0}
\end{align*}
generate $Y_n$ subject only to 
the relations (\ref{pr3})--(\ref{pr14}) below (taking
the shift matrix there to be the zero matrix).
For example, the presentation for $Y_n$ from $\S$\ref{syangian} is 
the special case $\nu = (1^n)$, in which case we denote
$D_{i;1,1}^{(r)}, E_{i;1,1}^{(r)}$ and $F_{i;1,1}^{(r)}$ 
simply by $D_i^{(r)}, E_i^{(r)}$ and $F_i^{(r)}$ respectively,
while the RTT presentation from (\ref{mr}) is the special
case $\nu = (n)$, in which case
$D_{1;i,j}^{(r)}=T_{i,j}^{(r)}$.

We are going to adapt these parabolic presentations to the shifted
Yangian $Y_n(\sigma)$. Return to the setup of $\S$\ref{syangian},
assuming that $\sigma = (s_{i,j})_{1 \leq i, j \leq n}$ is a 
fixed shift matrix.
Suppose in addition that the shape $\nu = (\nu_1,\dots,\nu_m)$ is 
{\em admissible} for $\sigma$, meaning that
$s_{i,j} = 0$ for all $\nu_1+\cdots+\nu_{a-1}+1 \leq i,j 
\leq \nu_1+\cdots+\nu_a$
and $a=1,\dots,m$. We will adopt the shorthand
\begin{equation}\label{rels}
s_{a,b}(\nu) := s_{\nu_1+\cdots+\nu_a,\nu_1+\cdots+\nu_b}.
\end{equation}
The shifts $(s_{i,j})_{1 \leq i,j \leq n}$ can be recovered from the
``relative'' shifts $(s_{a,b}(\nu))_{1 \leq a,b \leq m}$ 
given the admissible shape $\nu$.

Define a new algebra ${^\nu}Y_n(\sigma)$ over $\C$, which will shortly be
{identified} with $Y_n(\sigma)$ from $\S$\ref{syangian},
by
generators 
\begin{align*}
&\{D_{a;i,j}^{(r)}\}_{1 \leq a \leq m, 1 \leq i,j \leq \nu_a,
r > 0},\\
&\{E_{a;i,j}^{(r)}\}_{1 \leq a < m, 1 \leq i 
\leq \nu_a, 1 \leq j \leq \nu_{a+1},
r > s_{a,a+1}(\nu)},\\
&\{F_{a;i,j}^{(r)}\}_{1 \leq a < m, 1 \leq i 
\leq \nu_{a+1}, 1 \leq j \leq \nu_{a},
r > s_{a+1,a}(\nu)}
\end{align*}
subject to certain relations. 
To write down these relations,
we must introduce some further notation.
Let $D_{a;i,j}^{(0)} := \delta_{i,j}$, 
$D_{a;i,j}(u) := \sum_{r \geq 0} D_{a;i,j}^{(r)} u^{-r}$
and introduce the matrix 
$D_a(u) := (D_{a;i,j}(u))_{1 \leq i,j \leq \nu_a}$.
Let $\widetilde{D}_a(u) = 
(\widetilde{D}_{a;i,j}(u))_{1 \leq i,j \leq \nu_a}$
denote the matrix $-D_a(u)^{-1}$, and write $\widetilde{D}_{a;i,j}(u) =
\sum_{r \geq 0} \widetilde D_{a;i,j}^{(r)} u^{-r}$,
thus defining elements $\widetilde{D}_{a;i,j}^{(r)}$
of ${^\nu}Y_n(\sigma)$ for each $a=1,\dots,m$, $1 \leq i,j \leq \nu_a$
and $r \geq 0$.
In particular, $\widetilde{D}_{a;i,j}^{(0)} = - \delta_{i,j}$
and $\widetilde{D}_{a;i,j}^{(1)} = D_{a;i,j}^{(1)}$.
Now the relations are:
\begin{align}
[D_{a;i,j}^{(r)}, D_{b;h,k}^{(s)}] &=
\delta_{a,b}
\sum_{t=0}^{\min(r,s)-1}
\left(
D_{a;i,k}^{(r+s-1-t)}D_{a;h,j}^{(t)} 
-D_{a;i,k}^{(t)}D_{a;h,j}^{(r+s-1-t)} \right),\label{pr3}\\
[E_{a;i,j}^{(r)}, F_{b;h,k}^{(s)}]
&=
\delta_{a,b} \sum_{t=0}^{r+s-1}
\widetilde{D}_{a;i,k}^{(r+s-1-t)}D_{a+1;h,j}^{(t)},\label{pr6}
\end{align}\begin{align}
[D_{a;i,j}^{(r)}, E_{b;h,k}^{(s)}] &=
\delta_{a,b} 
\sum_{t=0}^{r-1} 
\sum_{g=1}^{\nu_a}D_{a;i,g}^{(t)} E_{a;g,k}^{(r+s-1-t)}\delta_{h,j}
- \delta_{a,b+1} \sum_{t=0}^{r-1}
D_{b+1;i,k}^{(t)} E_{b;h,j}^{(r+s-1-t)},\label{pr4}\\
[D_{a;i,j}^{(r)}, F_{b;h,k}^{(s)}] &=
\delta_{a,b+1} \sum_{t=0}^{r-1}
 F_{b;i,k}^{(r+s-1-t)}D_{b+1;h,j}^{(t)}-\delta_{a,b} 
\delta_{i,k}\sum_{t=0}^{r-1}\sum_{g=1}^{\nu_a} 
F_{a;h,g}^{(r+s-1-t)}D_{a;g,j}^{(t)},\label{pr5}\end{align}\begin{align}
[E_{a;i,j}^{(r)}, E_{a;h,k}^{(s)}] 
&=
\sum_{t=s_{a,a+1}(\nu)+1}^{s-1} E_{a;i,k}^{(t)} E_{a;h,j}^{(r+s-1-t)}
-\sum_{t=s_{a,a+1}(\nu)+1}^{r-1} E_{a;i,k}^{(t)} E_{a;h,j}^{(r+s-1-t)},\label{pr7}\\
[F_{a;i,j}^{(r)}, F_{a;h,k}^{(s)}] 
&=
\sum_{t=s_{a+1,a}(\nu)+1}^{r-1} F_{a;i,k}^{(r+s-1-t)}F_{a;h,j}^{(t)}-
\sum_{t=s_{a+1,a}(\nu)+1}^{s-1}
F_{a;i,k}^{(r+s-1-t)}F_{a;h,j}^{(t)},\label{pr8}\end{align}\begin{align}
[E_{a;i,j}^{(r)}, E_{a+1;h,k}^{(s+1)}]
-[E_{a;i,j}^{(r+1)}, E_{a+1;h,k}^{(s)}]
&=
-\sum_{g=1}^{\nu_{a+1}}E_{a;i,g}^{(r)} E_{a+1;g,k}^{(s)}\delta_{h,j},\label{pr9}\\
[F_{a;i,j}^{(r+1)}, F_{a+1;h,k}^{(s)}]
-[F_{a;i,j}^{(r)}, F_{a+1;h,k}^{(s+1)}]
&=
-\delta_{i,k}\sum_{g=1}^{\nu_{a+1}}F_{a+1;h,g}^{(s)}F_{a;g,j}^{(r)},\label{pr10}
\end{align}\begin{align}
[E_{a;i,j}^{(r)}, E_{b;h,k}^{(s)}] &= 0\:\:\qquad\text{if $b>a+1$ 
or if $b = a+1$ and $h \neq j$},\label{pr11}\\
[F_{a;i,j}^{(r)}, F_{b;h,k}^{(s)}] &= 0\:\:\qquad\text{if $b > a+1$
or if $b=a+1$ and $i \neq k$},\label{pr12}
\end{align}\begin{align}
[E_{a;i,j}^{(r)}, [E_{a;h,k}^{(s)}, E_{b;f,g}^{(t)}]]
+
[E_{a;i,j}^{(s)}, [E_{a;h,k}^{(r)}, E_{b;f,g}^{(t)}]] &= 0
\qquad\text{if }|a-b|=1,\label{pr13}\\
[F_{a;i,j}^{(r)}, [F_{a;h,k}^{(s)}, F_{b;f,g}^{(t)}]]
+
[F_{a;i,j}^{(s)}, [F_{a;h,k}^{(r)}, F_{b;f,g}^{(t)}]] &= 0
\qquad\text{if }|a-b|=1,\label{pr14}
\end{align}
for all admissible $a,b,f,g,h,i,j,k,r,s,t$.

Observe right away that there is a canonical
homomorphism ${^\nu}Y_n(\sigma) \rightarrow Y_n$
mapping the generators $D_{a;i,j}^{(r)}, 
E_{a;i,j}^{(r)}$ and $F_{a;i,j}^{(r)}$ of ${^\nu}Y_n(\sigma)$
to the elements
of $Y_n$  with the same names (the ones we defined above in 
terms of Gauss factorizations).
We are going to prove that this canonical
homomorphism is injective and that its image is independent
of the particular choice of the admissible shape $\nu$.
In particular this will identify ${^\nu}Y_n(\sigma)$ with the
algebra $Y_n(\sigma)$ from $\S$\ref{syangian}, 
since that is the special case $\nu = (1^n)$ of the present definition.

The proof that the canonical map ${^\nu}Y_n(\sigma) \rightarrow Y_n$
is injective is an extension of the proof given in
$\S$\ref{syangian}. For $1 \leq a < b \leq m$,
$1 \leq i \leq \nu_a, 1 \leq j \leq \nu_b$ and
$r > s_{a,b}(\nu)$, we define elements
$E_{a,b;i,j}^{(r)}$ inductively by
\begin{equation}\label{higheres}
E_{a,a+1;i,j}^{(r)} := E_{a;i,j}^{(r)},
\qquad
E_{a,b;i,j}^{(r)} := [E_{a,b-1;i,k}^{(r-s_{b-1,b}(\nu))},
E_{b-1;k,j}^{(s_{b-1,b}(\nu)+1)}]
\end{equation}
where $1 \leq k \leq \nu_{b-1}$.
By the relations
this definition is independent of the choice of $k$;
see for instance \cite[(6.9)]{BK} for a similar argument.
Similarly for $1 \leq a < b \leq m$,
$1 \leq i \leq \nu_b$, $1 \leq j \leq \nu_a$ and 
$r > s_{b,a}(\nu)$, we define
elements $F_{a,b;i,j}^{(r)}$ by
\begin{equation}\label{higherfs}
F_{a,a+1;i,j}^{(r)} := F_{a;i,j}^{(r)},
\qquad
F_{a,b;i,j}^{(r)} := [F_{b-1;i,k}^{(s_{b,b-1}(\nu)+1)},
F_{a,b-1;k,j}^{(r-s_{b,b-1}(\nu))}]
\end{equation}
where $1 \leq k \leq \nu_{b-1}$.
Also let $Y_\nu$ denote the subalgebra of ${^\nu}Y_n(\sigma)$
generated by all the $D_{a;i,j}^{(r)}$'s, let
$Y_\nu^+(\sigma)$ denote the subalgebra generated by all the
$E_{a;i,j}^{(r)}$'s and let
$Y_\nu^-(\sigma)$ denote the subalgebra generated by all the
$F_{a;i,j}^{(r)}$'s.
The following theorem generalizes Theorem \ref{pbw2}.

\begin{Theorem}\label{pbw3}
\begin{itemize}
\item[(i)] The monomials in the elements 
$\{D_{a;i,j}^{(r)}\}_{a=1,\dots,m, 1 \leq i,j \leq \nu_a,r>0}$ 
taken in some
fixed order form a basis for $Y_\nu$.
\item[(ii)] The monomials in the elements 
$\{E_{a,b;i,j}^{(r)}\}_{1 \leq a < b \leq m, 1 \leq i \leq \nu_a, 1 \leq j \leq \nu_b, r >s_{a,b}(\nu)}$ taken in some fixed
order form a basis for $Y_\nu^+(\sigma)$.
\item[(iii)] The monomials in the elements
$\{F_{a,b;i,j}^{(r)}\}_{1 \leq a < b \leq m, 
1 \leq i \leq \nu_b, 1 \leq j \leq \nu_a, r > s_{b,a}(\nu)}$ taken in some fixed
order form a basis for $Y_\nu^-(\sigma)$.
\item[(iv)] The monomials in the union of
the elements listed in (i)--(iii)
taken in some fixed order form a basis for ${^\nu}Y_n(\sigma)$.
\end{itemize}
\end{Theorem}

\begin{proof}
Introduce the loop filtration 
$\L_0 {^\nu}Y_n(\sigma) \subseteq \L_1 {^\nu}Y_n(\sigma) \subseteq
\cdots$
of ${^\nu}Y_n(\sigma)$ by declaring
that the generators $D_{a;i,j}^{(r)}, E_{a;i,j}^{(r)}$
and $F_{a;i,j}^{(r)}$ are all of degree $(r-1)$. 
Define elements $\{e_{i,j;r}\}_{1 \leq i,j \leq n, r \geq s_{i,j}}$
of the associated graded algebra $\gr^\LL {^\nu}Y_n(\sigma)$
from the equations
\begin{align}
\label{bday1}
\gr^\LL_{r} D_{a;i,j}^{(r+1)} &= 
e_{\nu_1+\cdots+\nu_{a-1}+i,\nu_1+\cdots+\nu_{a-1}+j;r},\\
\label{bday2}
\gr^\LL_{r} E_{a,b;i,j}^{(r+1)} &= 
e_{\nu_1+\cdots+\nu_{a-1}+i, \nu_1+\cdots+\nu_{b-1}+j;r},\\
\label{bday3}
\gr^\LL_{r} F_{a,b;i,j}^{(r+1)} &= 
e_{\nu_1+\cdots+\nu_{b-1}+i,
\nu_1+\cdots+\nu_{a-1}+j;r}.
\end{align}
Following the proof of 
\cite[Lemma 6.7]{BK}, one checks that these elements
satisfy the relations (\ref{assgr}).
Hence, there is a
well-defined surjective homomorphism
$\pi:U(\mathfrak{gl}_n[t](\sigma)) \twoheadrightarrow \gr^\LL {^\nu}Y_n(\sigma)$
mapping $e_{i,j} t^r \in U(\mathfrak{gl}_n[t](\sigma))$
to $e_{i,j;r} \in \gr^\LL {^\nu}Y_n(\sigma)$.

In the special case $\sigma = 0$, when we already know that 
${^\nu}Y_n(\sigma)$ is the usual
Yangian $Y_n$, one checks using Lemma~\ref{quasi} and induction
on the length $m$ of the shape $\nu$ that the element
$e_{i,j;r}$ defined here is equal to $\gr^\LL_r T_{i,j}^{(r+1)}$.
In particular $e_{i,j;r}$ coincides with the element
of $\gr^\LL_r Y_n$ defined by (\ref{newelts}).
Hence like in the proof of Theorem~\ref{pbw0},
the PBW theorem for the usual Yangian implies that the ordered monomials
in the elements
$\{e_{i,j;r}\}_{1 \leq i,j \leq n, r \geq 0}$ are linearly independent
in $\gr^\LL Y_n$, and $\pi$ is an isomorphism.

In general, the canonical map ${^\nu}Y_n(\sigma) \rightarrow Y_n$
is a homomorphism of filtered algebras, so induces a map
$\gr^\LL {^\nu}Y_n(\sigma) \rightarrow \gr^\LL Y_n$ which maps
$e_{i,j;r} \in \gr^\LL {^\nu}Y_n(\sigma)$ to $e_{i,j;r} \in \gr^\LL Y_n$.
So the previous paragraph implies that the ordered monomials
in the elements $\{e_{i,j;r}\}_{1 \leq i,j \leq n, r \geq s_{i,j}}$
are linearly independent in $\gr^\LL {^\nu}Y_n(\sigma)$ too.
Hence $\pi$ is an isomorphism
in general.

The theorem now follows like Theorem~\ref{pbw2}.
\end{proof}

The following two corollaries generalize Corollaries \ref{injcor} and 
\ref{triangular}.

\begin{Corollary}
The canonical map ${^\nu}Y_n(\sigma) \rightarrow Y_n$ is injective.
\end{Corollary}

\begin{Corollary}
Multiplication $Y_\nu^-(\sigma) \times Y_\nu \times Y_\nu^+(\sigma)
\rightarrow {^\nu}Y_n(\sigma)$ is a vector space isomorphism.
\end{Corollary}

So now for each admissible shape $\nu$, we have defined
a subalgebra ${^\nu}Y_n(\sigma)$ of $Y_n$. 
It remains to see that
these subalgebras coincide for different $\nu$.
Suppose $\nu_b = \alpha+\beta$ for some $1 \leq b \leq m$ and
$\alpha,\beta \geq 1$, and let
$\mu = (\nu_1,\dots,\nu_{b-1},\alpha,\beta,\nu_{b+1},\dots,\nu_m)$. 
Then it suffices to show that ${^\nu} Y_n(\sigma) = {^\mu}Y_n(\sigma)$
as subalgebras of $Y_n$.
Using Lemma~\ref{quasi}, one checks that
${^\mu}Y_n(\sigma) \subseteq {^\nu} Y_n(\sigma)$.
Now the equality ${^\mu}Y_n(\sigma) = {^\nu}Y_n(\sigma)$
follows easily because we have already seen in the proof of
Theorem~\ref{pbw3} that 
their associated graded algebras are equal in $\gr^\LL Y_n$.
We have now proved that the relations (\ref{pr3})--(\ref{pr14})
give presentations for the shifted Yangian $Y_n(\sigma)={^\nu}Y_n(\sigma)$
for each admissible shape $\nu$. 

\begin{Remark}\label{parrem}\rm
As a first application of 
these parabolic presentations, one can introduce
analogues of parabolic subalgebras of $Y_n(\sigma)$: for
an admissible shape $\nu$,
define $Y_\nu^\sharp(\sigma) := Y_\nu Y_\nu^+(\sigma)$
and $Y_\nu^\flat(\sigma) := Y_\nu^-(\sigma) Y_\nu$.
By the relations, these are indeed subalgebras of $Y_n(\sigma)$.
Moreover, there are obvious surjective homomorphisms
$Y_\nu^\sharp(\sigma) \twoheadrightarrow Y_\nu$
and 
$Y_\nu^\flat(\sigma) \twoheadrightarrow Y_\nu$ with kernels
generated by all $E_{a;i,j}^{(r)}$ and all
$F_{a;i,j}^{(r)}$ respectively.
\end{Remark}

\begin{Remark}\rm
In this remark, we describe 
the maps $\tau$ and $\iota$ from (\ref{taudef})--(\ref{iotadef}) 
in terms of the parabolic generators. 
The anti-isomorphism $\tau$ satisfies
\begin{equation}
\tau(D_{a;i,j}^{(r)}) =
D_{a;j,i}^{(r)},\quad
\tau(E_{a;i,j}^{(r)}) = 
F_{a;j,i}^{(r)},\quad
\tau(F_{a;i,j}^{(r)}) = 
E_{a;j,i}^{(r)};
\end{equation}
cf. \cite[(6.6)--(6.8)]{BK}.
Also, in the notation of (\ref{iotadef}) and working with a fixed
admissible shape $\nu$ (for both $\sigma$ and $\dot\sigma$), 
the isomorphism $\iota$ satisfies
\begin{multline}
\quad\iota(D_{a;i,j}^{(r)}) =
\dot D_{a;i,j}^{(r)},\qquad
\iota(E_{a;i,j}^{(r)}) =
\dot E_{a;i,j}^{(r-s_{a,a+1}(\nu)+\dot s_{a,a+1}(\nu))},\\
\iota(F_{a;i,j}^{(r)}) =
\dot F_{a;i,j}^{(r-s_{a+1,a}(\nu)+\dot s_{a+1,a}(\nu))}.\quad
\label{iotadefpar}
\end{multline}
To see this, note by the relations 
(\ref{pr3})--(\ref{pr14}) that these formulae
certainly give a well-defined isomorphism 
$\iota:Y_n(\sigma) \rightarrow 
Y_n(\dot\sigma)$. Now use Lemma~\ref{quasi} to check 
inductively that this is the same map as in 
(\ref{iotadef}).
\end{Remark}

\section{Baby comultiplications}\label{sbaby}

Now let us explain the real reason why the parabolic presentations
of $Y_n(\sigma)$ from $\S$\ref{sparabolic} are so important.
Fix a shift matrix $\sigma = (s_{i,j})_{1 \leq i,j \leq n}$.
Recall a shape $\nu = (\nu_1,\dots,\nu_m)$ is
{admissible} if
$s_{i,j} = 0$ for all $\nu_1+\cdots+\nu_{a-1}+1 
\leq i,j \leq \nu_1+\cdots+\nu_a$
and $a=1,\dots,m$.
Throughout the section, $\nu$
will denote the {\em minimal admissible shape} for $\sigma$,
that is, the admissible shape of smallest possible length $m$.
We will always work in terms of the parabolic
presentation defined relative to this fixed shape $\nu$.
For example in the case of the Yangian $Y_n$ itself, this is
the RTT presentation from (\ref{mr}).

Consider first the special case that $\sigma$ is the zero matrix,
i.e. the case when $Y_n(\sigma)$ is the usual Yangian $Y_n$. It
is well known that $Y_n$ is a Hopf algebra, with 
comultiplication $\Delta:Y_n \rightarrow Y_n \otimes Y_n$
which may be defined in terms of the RTT presentation
by the equation
\begin{equation}\label{comult}
\Delta(T_{i,j}^{(r)}) = \sum_{s=0}^r \sum_{k=1}^n
T_{i,k}^{(s)} \otimes T_{k,j}^{(r-s)}.
\end{equation}
There is also the 
{\em evaluation homomorphism}
$\kappa_1:Y_n \rightarrow U(\mathfrak{gl}_n)$
defined by 
\begin{equation}
\kappa_1(T_{i,j}^{(r)}) = 
\left\{
\begin{array}{ll}
e_{i,j}&\hbox{if $r = 1$,}\\
0&\hbox{if $r > 1$.}
\end{array}\right.
\end{equation}
Let $\Delta_{\rt} := (\id\otimes\kappa_1) \circ \Delta$
and $\Delta_{\lt} := (\kappa_1\otimes\id) \circ \Delta$, thus defining
algebra homomorphisms
\begin{align}
\Delta_{\rt}:Y_n \rightarrow Y_n\otimes U(\mathfrak{gl}_n),
\qquad
&T_{i,j}^{(r)} \mapsto T_{i,j}^{(r)} \otimes 1 + \sum_{k=1}^n 
T_{i,k}^{(r-1)} \otimes e_{k,j},
\label{babyr}\\
\Delta_{\lt}:Y_n \rightarrow U(\mathfrak{gl}_n) \otimes Y_n,
\qquad
&T_{i,j}^{(r)} \mapsto 1 \otimes T_{i,j}^{(r)} + \sum_{k=1}^n e_{i,k} \otimes 
T_{k,j}^{(r-1)}.
\label{babyl}
\end{align}
The following theorem 
defines analogous ``baby
comultiplications'' for the 
shifted Yangians in general.

\begin{Theorem}\label{baby}
Let $\nu = (\nu_1,\dots,\nu_m)$ be the minimal admissible shape for $\sigma$
and set $t := \nu_m$.
For $1 \leq i,j \leq t$, let $\tilde e_{i,j} := e_{i,j} + \delta_{i,j}(n-t)
\in U(\mathfrak{gl}_t)$.
\begin{itemize}
\item[(i)] If either 
$t = n$ or $s_{n-t,n-t+1} \neq 0$, define
$\dot\sigma = (\dot s_{i,j})_{1 \leq i,j \leq n}$ 
from 
\begin{equation}\label{sr}
\dot s_{i,j} = \left\{
\begin{array}{ll}
s_{i,j}-1&\hbox{if $i \leq n-t < j$,}\\
s_{i,j}&\hbox{otherwise.}
\end{array}\right.
\end{equation}
Then, there is a unique algebra homomorphism
$\Delta_{\rt}:Y_n(\sigma) \rightarrow 
Y_n(\dot\sigma) 
\otimes U(\mathfrak{gl}_t)$
such that
\begin{align*}
D_{a;i,j}^{(r)} &\mapsto \dot D_{a;i,j}^{(r)} \otimes 1
+ \delta_{a,m} \sum_{k=1}^{t}\dot D_{a;i,k}^{(r-1)} \otimes \tilde e_{k,j},\\
E_{a;i,j}^{(r)} &\mapsto \dot E_{a;i,j}^{(r)} \otimes 1
+ \delta_{a,m-1} \sum_{k=1}^{t}\dot E_{a;i,k}^{(r-1)} \otimes \tilde e_{k,j},\\
F_{a;i,j}^{(r)} & \mapsto \dot F_{a;i,j}^{(r)} \otimes 1.
\end{align*}
\item[(ii)]
If either $t = n$ or $s_{n-t+1,n-t} \neq 0$, define
$\dot\sigma = (\dot s_{i,j})_{1 \leq i,j \leq n}$ 
from 
\begin{equation}\label{sl}
\dot s_{i,j} = \left\{
\begin{array}{ll}
s_{i,j}-1&\hbox{if $j \leq n-t < i$,}\\
s_{i,j}&\hbox{otherwise.}
\end{array}\right.
\end{equation}
Then, there is a unique algebra homomorphism
$\Delta_{\lt}:Y_n(\sigma) \rightarrow U(\mathfrak{gl}_t) \otimes Y_n(\dot \sigma)$
such that
\begin{align*}
D_{a;i,j}^{(r)} &\mapsto 1 \otimes \dot D_{a;i,j}^{(r)}
+ \delta_{a,m} \sum_{k=1}^{t} \tilde e_{i,k} \otimes \dot D_{a;k,j}^{(r-1)},\\
E_{a;i,j}^{(r)} &\mapsto 1 \otimes \dot E_{a;i,j}^{(r)},\\
F_{a;i,j}^{(r)} & \mapsto 1 \otimes \dot F_{a;i,j}^{(r)}
+ \delta_{a,m-1} \sum_{k=1}^{t}\tilde e_{i,k}\otimes \dot F_{a;k,j}^{(r-1)}.
\end{align*}
\end{itemize}
(In both (i) and (ii) we are working in terms of the parabolic generators
relative to the fixed shape $\nu$.)
\end{Theorem}

\begin{proof}
Check the relations (\ref{pr3})--(\ref{pr14})
(to check the relations (\ref{pr13}) and (\ref{pr14}) one needs to use
(\ref{pr7}), (\ref{pr8}) and (\ref{pr9}), (\ref{pr10}) several times).
\end{proof}

The next lemma explains how to compute the
baby comultiplications on the higher root elements
$E_{a,b;i,j}^{(r)}$ and $F_{a,b;i,j}^{(r)}$, still
working in terms of
the minimal admissible shape $\nu$ for $\sigma$.

\begin{Lemma}\label{higherroots}
\begin{itemize}
\item[(i)] Under the hypotheses of Theorem~\ref{baby}(i),
we have for $1 \leq a < b-1 < m$ and all admissible $i,j,r$ that
\begin{align*}
\Delta_{\rt}(E_{a,b;i,j}^{(r)}) &= \left\{
\begin{array}{ll}
\!\dot E_{a,b;i,j}^{(r)} \otimes 1&\hbox{if $b < m$,}\\
\!{[\dot E_{a,m-1;i,h}^{(r-s_{m-1,m}(\nu))},\dot E_{m-1;h,j}^{(s_{m-1,m}(\nu)+1)}]
\otimes 1+\sum_{k=1}^{t}\dot E_{a,m;i,k}^{(r-1)} \otimes \tilde 
e_{k,j}}&\hbox{if $b = m$,}
\end{array}\right.\\
\Delta_{\rt}(F_{a,b;i,j}^{(r)}) &= \dot F_{a,b;i,j}^{(r)} \otimes 1,
\end{align*}
for any $1 \leq h \leq \nu_{m-1}$.
\item[(ii)] Under the hypotheses of Theorem~\ref{baby}(ii),
we have for $1 \leq a < b-1 < m$ and all admissible $i,j,r$ that
\begin{align*}
\!\Delta_{\lt}(E_{a,b;i,j}^{(r)}) &= 1 \otimes \dot E_{a,b;i,j}^{(r)},\\
\!\Delta_{\lt}(F_{a,b;i,j}^{(r)}) &=
\left\{
\begin{array}{ll}
1 \otimes \dot F_{a,b;i,j}^{(r)}&\hbox{if $b < m$,}\\
1 \otimes [\dot F_{m-1;i,h}^{(s_{m,m-1}(\nu)+1)},\dot F_{a,m-1;h,j}^{(r-s_{m,m-1}(\nu))}]+\sum_{k=1}^{t} \tilde e_{i,k} \otimes \dot F_{a,m;k,j}^{(r-1)}&\hbox{if $b = m$,}
\end{array}\right.
\end{align*}
for any $1 \leq h \leq \nu_{m-1}$.
\end{itemize}
\end{Lemma}

\begin{proof}
Let us just explain how to compute
$\Delta_{\rt}(E_{a,m;i,j}^{(r)})$ for $1 \leq a < m-1$, since all the other cases are similar.
By definition, 
$E_{a,m;i,j}^{(r)} = 
[E_{a,m-1;i,h}^{(r-s_{m-1,m}(\nu))},
E_{m-1;h,j}^{(s_{m-1,m}(\nu)+1)}]$
for any $1 \leq h \leq \nu_{m-1}$.
Clearly
$\Delta_{\rt}(E_{a,m-1;i,h}^{(r-s_{m-1,m}(\nu))})
= \dot E_{a,m-1;i,h}^{(r-s_{m-1,m}(\nu))} \otimes 1$. Hence,
\begin{align*}
\Delta_{\rt}(E_{a,m;i,j}^{(r)})&=
\left[\dot E_{a,m-1;i,h}^{(r-s_{m-1,m}(\nu))} \otimes 1,
\dot E_{m-1;h,j}^{(s_{m-1,m}(\nu)+1)} \otimes 1
+ \sum_{k=1}^{t}\dot E_{m-1;h,k}^{(s_{m-1,m}(\nu))} \otimes \tilde e_{k,j}\right]\\
&=
\left[\dot E_{a,m-1;i,h}^{(r-s_{m-1,m}(\nu))},\dot E_{m-1;h,j}^{(s_{m-1,m}(\nu)+1)}\right ] \otimes 1
+ \sum_{k=1}^{t}\dot E_{a,m;i,k}^{(r-1)} \otimes \tilde e_{k,j},
\end{align*}
as claimed.
\end{proof}

\begin{Remark}\label{areinj}\rm
In the notation of Theorem~\ref{baby},
the maps $\Delta_{\rt}$ and $\Delta_{\lt}$ are injective.
One can see this as follows.
Let $\eps:U(\mathfrak{gl}_t) \rightarrow \C$ be the homomorphism
with $\eps(\tilde e_{i,j})=0$ for $1 \leq i,j \leq t$. By definition,
$Y_n(\sigma)$ and $Y_n(\dot \sigma)$ are subalgebras of the
Yangian with $Y_n(\sigma) \subseteq Y_n(\dot\sigma)$.
Now the point is that the compositions
$(\id \bar\otimes \eps) \circ \Delta_{\rt}$
and $(\eps \bar\otimes \id)\circ \Delta_{\lt}$  
coincide (when defined) with
the natural embedding
$Y_n(\sigma) \hookrightarrow Y_n(\dot\sigma)$.
\end{Remark}

\section{The canonical filtration}\label{scan}

In this section, we introduce another important filtration
of the shifted Yangian. It is easiest to start
with the Yangian $Y_n$ itself defined by the RTT presentation
as in $\S$\ref{sparabolic}. The {\em canonical filtration}
$\F_0 Y_n \subseteq \F_1 Y_n \subseteq \cdots$ of $Y_n$
is defined by declaring that the generators $T_{i,j}^{(r)}$
are all of degree $r$, i.e. $\F_d Y_n$ is the span of all monomials
in these generators of total degree $\leq d$. It is obvious from
the relations (\ref{mr}) that the associated graded
algebra $\gr Y_n$ is commutative.

Suppose instead that we are given a shape
$\nu = (\nu_1,\dots,\nu_m)$.
Using the explicit definitions from 
\cite[(6.1)--(6.4)]{BK}, one checks
that the parabolic generators
$D_{a;i,j}^{(r)}, E_{a,b;i,j}^{(r)}$
and $F_{a,b;i,j}^{(r)}$ of $Y_n$ are linear combinations of monomials
in the $T_{i,j}^{(s)}$ of total degree $r$ and 
conversely, if we define
$D_{a;i,j}^{(r)}, E_{a,b;i,j}^{(r)}$
and $F_{a,b;i,j}^{(r)}$ all to be of degree $r$, each $T_{i,j}^{(s)}$
is a linear combination of monomials in these elements
of total degree $s$.
Hence $\F_d Y_n$ can also be described as the span of all monomials
in the elements
$D_{a;i,j}^{(r)}, E_{a,b;i,j}^{(r)}$
and $F_{a,b;i,j}^{(r)}$ of total degree $\leq d$.
For $1 \leq a,b \leq m$,
$1 \leq i \leq \nu_a, 1 \leq j \leq \nu_b$ and $r > 0$, define
\begin{equation}\label{eelts}
e_{a,b;i,j}^{(r)} := \left\{
\begin{array}{ll}
\gr_r D_{a;i,j}^{(r)}&\hbox{if $a=b$,}\\
\gr_r E_{a,b;i,j}^{(r)}&\hbox{if $a < b$,}\\
\gr_r F_{b,a;i,j}^{(r)}&\hbox{if $a > b$,}
\end{array}
\right.
\end{equation}
all elements of $\gr_r Y_{n}$.
Of course, this notation depends implicitly on the fixed shape $\nu$.
The commutativity of $\gr Y_n$ together with
Theorem~\ref{pbw3}(iv) imply:

\begin{Theorem} \label{assgr1}
For any shape $\nu = (\nu_1,\dots,\nu_m)$,
$\gr Y_n$ is the free commutative algebra on generators
$\{e_{a,b;i,j}^{(r)}\}_{1 \leq a, b \leq m,
1 \leq i \leq \nu_a, 1 \leq j \leq \nu_b, r > 0}$.
\end{Theorem}

Now consider the shifted Yangian. So let $\sigma = (s_{i,j})_{1 \leq i,j \leq n}$ be a shift matrix as in $\S$\ref{syangian}.
Viewing $Y_n(\sigma)$ as a subalgebra
of $Y_n$, introduce the {\em canonical filtration} 
$\F_0 Y_n(\sigma)\subseteq \F_1 Y_n(\sigma) \subseteq \cdots$
of $Y_n(\sigma)$ by defining
$\F_d Y_n(\sigma) := Y_n(\sigma) \cap \F_d Y_n$.
Thus, the inclusion of $Y_n(\sigma)$ into $Y_n$ is filtered and the induced
map $\gr Y_n(\sigma) \rightarrow \gr Y_n$ is injective. Hence, 
$\gr Y_n(\sigma)$ is identified with a subalgebra of the commutative
algebra $\gr Y_n$. The following theorem describes this subalgebra
explicitly.

\begin{Theorem} \label{assgr2}
For an admissible shape 
$\nu = (\nu_1,\dots,\nu_m)$,
$\gr Y_n(\sigma)$ is the 
subalgebra of $\gr Y_n$
generated by the elements
$\{e_{a,b;i,j}^{(r)}\}_{1 \leq a, b \leq m,
1 \leq i \leq \nu_a, 1 \leq j \leq \nu_b, r > s_{a,b}(\nu)}$.
\end{Theorem}

\begin{proof}
Note in view of the relations (\ref{pr9})--(\ref{pr10}) that
the element $e_{a,b;i,j}^{(r)}$ of
$\gr Y_n(\sigma)$ is identified with the element of the same
name in $\gr Y_n$.
Given this the theorem follows
directly from Theorems~\ref{assgr1} and \ref{pbw3}(iv).
\end{proof}

\begin{Remark}\label{simply}\rm
Theorem~\ref{assgr2} means that, given an admissible shape $\nu$,
we can define the canonical filtration
on $Y_n(\sigma)$ intrinsically simply by declaring that the
elements $D_{a;i,j}^{(r)}, E_{a,b;i,j}^{(r)}$
and $F_{a,b;i,j}^{(r)}$ of $Y_n(\sigma)$ are all of degree $r$, i.e. $\F_d
Y_n(\sigma)$ is the span of all monomials in these elements of total
degree $\leq d$.
In particular, this definition is independent of the
particular choice of admissible shape $\nu$.
\end{Remark}

\begin{Remark}\label{sumply}
\rm
The comultiplication $\Delta:Y_n \rightarrow Y_n \otimes Y_n$
is a filtered map with respect to the canonical filtration,
as follows immediately from 
(\ref{comult}).
Similarly, 
the baby comultiplications
$\Delta_{\rt}:Y_n(\sigma) \rightarrow Y_n(\dot\sigma)\otimes 
U(\mathfrak{gl}_t)$
and $\Delta_\lt:Y_n(\sigma) \rightarrow U(\mathfrak{gl}_t)\otimes
Y_n(\dot\sigma)$ 
from Theorem~\ref{baby} are filtered maps whenever they are defined,
providing we extend the canonical filtration 
of $Y_n(\dot\sigma)$
to
$Y_n(\dot\sigma)\otimes 
U(\mathfrak{gl}_t)$
and 
$U(\mathfrak{gl}_t)\otimes Y_n(\dot\sigma)$
by declaring that the matrix units $e_{i,j} \in \mathfrak{gl}_t$
are of degree $1$.
This follows from Theorem~\ref{baby} and Lemma~\ref{higherroots}.
The argument explained in Remark~\ref{areinj} shows moreover that the
associated graded maps 
$\gr \Delta_\rt$ and $\gr \Delta_\lt$ are injective.
\end{Remark}

\section{Truncation}\label{struncation}

Continue with a fixed
shift matrix $\sigma = (s_{i,j})_{1 \leq i,j \leq n}$.
Fix also an integer $l \geq  
s_{1,n}+s_{n,1}$ (the ``level''), 
and for $i=1,\dots,n$ let
\begin{equation}\label{pidef}
p_i := l -s_{i,n} - s_{n,i},
\end{equation}
thus defining a tuple $(p_1,\dots,p_n)$ of
integers with $0 \leq p_1 \leq \cdots \leq p_n = l$.
If $\nu = (\nu_1,\dots,\nu_m)$ is an admissible shape, we will
also use the shorthand
\begin{equation}\label{padef}
p_a(\nu) := p_{\nu_1+\cdots+\nu_a}
\end{equation}
for each $a=1,\dots,m$; cf. (\ref{rels}).
(A better way to keep track of all this data will be explained in the next
section.)

The {\em shifted Yangian of level l}, 
denoted $Y_{n,l}(\sigma)$,
is defined to be the quotient of $Y_n(\sigma)$ by the
two-sided ideal generated by the elements
$\{D_1^{(r)}\}_{r > p_1}$.
Alternatively, in terms of the parabolic presentation 
relative to an admissible shape $\nu = (\nu_1,\dots,\nu_m)$,
$Y_{n,l}(\sigma)$ is the quotient of $Y_n(\sigma)$
by the two-sided
ideal generated by the elements
$\{D_{1;i,j}^{(r)}\}_{1 \leq i,j \leq \nu_1, r > p_1}$.
The equivalence of all these definitions is easy to see
since $D_{1;1,1}^{(r)} = D_1^{(r)}$ and 
all other $D_{1;i,j}^{(r)}$ for $r > p_1$ 
and $1 \leq i,j \leq \nu_1$
obviously lie in the two-sided
ideal generated by $\{D_{1;1,1}^{(r)}\}_{r > p_1}$ 
in view of the relation (\ref{pr3}).
For example, in the case that $\sigma$ is the zero matrix, when
we denote $Y_{n,l}(\sigma)$ simply by $Y_{n,l}$, this algebra
is the quotient of $Y_n$ by the two-sided ideal generated
by $\{T_{i,j}^{(r)}\:|\:1 \leq i,j \leq n, r > l\}$,
which is precisely the definition of the 
{\em Yangian of level $l$} from \cite{Ch}.
\iffalse
By convention, we also write 
$Y_{0,0} = Y_{0,0}(\sigma)$ for the trivial algebra $\C$.
\fi

In the hope that it is clear from context whether we are talking about
$Y_n(\sigma)$ or $Y_{n,l}(\sigma)$, we will abuse notation
and use the same symbols
$D_{a;i,j}^{(r)}$,
$\widetilde{D}_{a;i,j}^{(r)}$,
$E_{a,b;i,j}^{(r)}$ and $F_{a,b;i,j}^{(r)}$
both for the generators 
of $Y_n(\sigma)$ and for their canonical images in 
the quotient $Y_{n,l}(\sigma)$.
\ifcenters@ 
Similarly we
use the same notation 
$C_n^{(r)}$ for the images in 
$Y_{n,l}(\sigma)$ of the central elements of $Y_n(\sigma)$
from (\ref{centelts});
clearly these are also central in 
$Y_{n,l}(\sigma)$.
\fi
The anti-isomorphism $\tau$ from (\ref{taudef})
factors through the quotients to induce an anti-isomorphism
\begin{equation}
\tau:Y_{n,l}(\sigma) \rightarrow Y_{n,l}(\sigma^t).
\end{equation}
Similarly, 
given another shift matrix $\dot\sigma$ 
such that $\dot s_{i,i+1}+\dot s_{i+1,i} = s_{i,i+1}+s_{i+1,i}$
for each $i=1,\dots,n-1$,
the isomorphism $\iota$ from (\ref{iotadef}) induces 
an isomorphism 
\begin{equation}\label{iotadef2}
\iota:Y_{n,l}(\sigma) \rightarrow Y_{n,l}(\dot\sigma).
\end{equation}
So the isomorphism type of 
the algebra $Y_{n,l}(\sigma)$ only actually depends on the
tuple $(p_1,\dots,p_n)$.

We will exploit the canonical filtration
$\F_0 Y_{n,l}(\sigma) \subseteq
\F_1 Y_{n,l}(\sigma) \subseteq \cdots$
of $Y_{n,l}(\sigma)$ 
induced by the quotient map
$Y_{n}(\sigma) \twoheadrightarrow Y_{n,l}(\sigma)$
and the canonical filtration of $Y_n(\sigma)$
from $\S$\ref{scan}.
Recalling Remark~\ref{simply},
this may be defined directly
given an admissible shape $\nu = (\nu_1,\dots,\nu_m)$ by
declaring all the elements
$D_{a;i,j}^{(r)}, E_{a,b;i,j}^{(r)}$
and $F_{a,b;i,j}^{(r)}$ of $Y_{n,l}(\sigma)$
to be of degree $r$, and then 
$\F_d Y_{n,l}(\sigma)$ is the span of all monomials
in these elements of total degree $\leq d$.
For $1 \leq a,b \leq m, 1 \leq i \leq \nu_a, 1 \leq j \leq \nu_b$
and $r > s_{a,b}(\nu)$, define elements
$e_{a,b;i,j}^{(r)}$ of the associated graded algebra
$\gr Y_{n,l}(\sigma)$
by the formula (\ref{eelts}).
Since $\gr Y_{n,l}(\sigma)$ is a quotient of 
$\gr Y_n(\sigma)$, Theorem~\ref{assgr2} implies that 
it is commutative and
is generated by all
$\{e_{a,b;i,j}^{(r)}\}_{1 \leq a,b \leq m, 
1 \leq i \leq \nu_a, 1 \leq j \leq \nu_b, r >s_{a,b}(\nu)}$.

\begin{Lemma}\label{yak}
For any admissible shape $\nu = (\nu_1,\dots,\nu_m)$, 
$\gr Y_{n,l}(\sigma)$
is generated just by the elements
$\{e_{a,b;i,j}^{(r)}\}_{1 \leq a,b \leq m, 
1 \leq i \leq \nu_a, 1 \leq j \leq \nu_b, s_{a,b}(\nu) <
r \leq s_{a,b}(\nu) + 
p_{\min(a,b)}(\nu)}.$
\end{Lemma}

\begin{proof}
For $1 \leq  c \leq m$, let $\Omega_c$ denote the set
\begin{align*}
\{D_{a;i,j}^{(r)}\}_{1 \leq a \leq c, 
1 \leq i,j 
\leq \nu_a, 0 < r \leq p_a(\nu)}
&\cup
\{E_{a,b;i,j}^{(r)}\}_{1 \leq a < b \leq c,
1 \leq i \leq \nu_a, 1 \leq j \leq \nu_b, s_{a,b}(\nu) < r \leq
s_{a,b}(\nu)+p_a(\nu)}\\
&\cup\{F_{a,b;i,j}^{(r)}\}_{1 \leq a < b \leq c,
1 \leq j \leq \nu_a, 1 \leq i \leq \nu_b, s_{b,a}(\nu) < r \leq
s_{b,a}(\nu)+p_a(\nu)},\\\intertext{and 
let $\widehat\Omega_c$ denote}
\{D_{a;i,j}^{(r)}\}_{1 \leq a \leq c, 1 \leq i,j 
\leq \nu_a, r > 0}
&\cup
\{E_{a,b;i,j}^{(r)}\}_{1 \leq a < b \leq c,
1 \leq i \leq \nu_a, 1 \leq j \leq \nu_b, r > s_{a,b}(\nu)}\\
&\cup\{F_{a,b;i,j}^{(r)}\}_{1 \leq a < b \leq c,
1 \leq j \leq \nu_a, 1 \leq i \leq \nu_b, r > s_{b,a}(\nu)}.
\end{align*}
To prove the lemma, we show by induction on 
$c=1,\dots,m$ 
that every element
of ${\widehat \Omega}_c$ of degree $r$ can be expressed as a linear
combination of monomials in the elements of 
$\Omega_c$ of total
degree $r$. 
The case $c=1$ is clear, since
$D_{1;i,j}^{(r)} = 0$ for $r > p_1(\nu)$ by definition.
Now take $c > 1$ and assume the result has been proved for all smaller
$c$.

For $1 \leq a < c-1$, we have by definition that
$E_{a,c;i,j}^{(r)} = 
[E_{a,c-1;i,k}^{(r-s_{c-1,c}(\nu))}, 
E_{c-1;k,j}^{(s_{c-1,c}(\nu)+1)}]$ 
for some $1 \leq k \leq \nu_{c-1}$.
By induction, $E_{a,c-1;i,k}^{(r-s_{c-1,c}(\nu))}$ is a linear combination
of monomials in the elements of $\Omega_{c-1}$ of total degree
$(r-s_{c-1,c}(\nu))$. By the relations, the commutator of any such
monomial with $E_{c-1;k,j}^{(s_{c-1,c}(\nu)+1)}$ 
is a linear combination of monomials in the elements of $\Omega_c$
of total degree $r$. Hence, $E_{a,c;i,j}^{(r)}$
is a linear combination of monomials in the elements of $\Omega_c$
of total degree $r$.
Similarly, so is $F_{a,c;i,j}^{(r)}$.

Next, we have by the relations that
\begin{equation}\label{step1}
E_{c-1;i,j}^{(r)} = 
[D_{c-1;i,i}^{(r-s_{c-1,c}(\nu))},
E_{c-1;i,j}^{(s_{c-1,c}(\nu)+1)}]
-\sum_{t=1}^{r-s_{c-1,c}(\nu)-1}
\sum_{h=1}^{\nu_{c-1}}
D_{c-1;i,h}^{(t)} E_{c-1;h,j}^{(r-t)}.
\end{equation}
By induction, $D_{c-1;i,i}^{(r-s_{c-1,c}(\nu))}$
is a linear combination
of monomials in the elements of $\Omega_{c-1}$ of total degree
$(r-s_{c-1,c}(\nu))$. Hence by the relations, the first term on the right
hand side of (\ref{step1})
is a linear combination of monomials in the elements of
$\Omega_{c-1} \cup \{E_{a,c;i,j}^{(r)}\}_{1 \leq a < c, 
1 \leq i \leq \nu_a,
1 \leq j \leq \nu_c, 
s_{a,c}(\nu) < r \leq 
s_{a,c}(\nu)+p_a(\nu)}$
of total degree $r$.
Now using (\ref{step1}) and induction on $r$,
one deduces that 
each $E_{c-1;i,j}^{(r)}$
is a linear combination of monomials in the elements of $\Omega_c$
of total degree $r$.
A similar argument 
\iffalse
using the relation
\begin{equation}\label{step11}
F_{c-1;i,j}^{(r)} = 
[F_{c-1;i,j}^{(s_{c,c-1}(\nu)+1)}, D_{c-1;j,j}^{(r-s_{c,c-1}(\nu))}]
-\sum_{t=1}^{r-s_{c,c-1}(\nu)-1}
\sum_{h=1}^{\nu_{c-1}}
F_{c-1;i,h}^{(r-t)}D_{c-1;h,j}^{(t)}
\end{equation}
\fi
shows that each $F_{c-1;i,j}^{(r)}$
is a linear combination of monomials in the elements of $\Omega_c$
of total degree $r$, too.

Finally we must consider $D_{c;i,j}^{(r)}$. 
Recall that
$p_c(\nu) = p_{c-1}(\nu)+s_{c-1,c}(\nu)+
s_{c,c-1}(\nu)$.
By the relations, we have 
for $1 \leq k \leq \nu_{c-1}$ that
\begin{equation}\label{step2}
D_{c;i,j}^{(r)} = \sum_{t=0}^{r-1} 
\widetilde{D}_{c-1;k,k}^{(r-t)}
D_{c;i,j}^{(t)} - [E_{c-1;k,j}^{(r-s_{c,c-1}(\nu))},
F_{c-1;i,k}^{(s_{c,c-1}(\nu)+1)}].
\end{equation}
By the previous paragraph, $E_{c-1;k,j}^{(r-s_{c,c-1}(\nu))}$ is a linear combination of monomials in the elements from $\Omega_{c-1}
\cup \{E_{a,c;i,j}^{(r)}\}_{1 \leq a < 
c, 1 \leq i \leq \nu_a,
1 \leq j \leq \nu_c, s_{a,c}(\nu) < r \leq s_{a,c}(\nu)+p_a(\nu)}$
of total degree $(r - s_{c,c-1}(\nu))$. Hence by the relations, 
the second term on the right hand side of (\ref{step2}) is a linear
combination of monomials in the elements of $\Omega_c$
of total degree $r$.
Now using (\ref{step2}) and induction on $r$ one deduces
that each 
$D_{c;i,j}^{(r)}$ is a linear combination of monomials in the elements
of $\Omega_c$ of total degree $r$, to complete the induction step.
\end{proof}

Assume additionally for the next paragraph that the level $l$ is
$> 0$.
Let $t \in \{1,\dots,n\}$ be maximal such that
$s_{i,j}=0$ for all $n-t+1 \leq i,j \leq n$.
If either $t = n$ or 
$s_{n-t,n-t+1} \neq 0$,
it is easy to check that the 
baby comultiplication $\Delta_{\rt}$
from Theorem~\ref{baby}(i) factors through the quotients to
induce a map
\begin{equation}\label{deltar}
\Delta_{\rt}:Y_{n,l}(\sigma) 
\rightarrow Y_{n,l-1}(\dot\sigma) \otimes 
U(\mathfrak{gl}_t),
\end{equation}
where $\dot\sigma$ is as in (\ref{sr}).
Similarly, if either $t = n$ or $s_{n-t+1,n-t} \neq 0$,
the baby comultiplication $\Delta_{\lt}$
from Theorem~\ref{baby}(ii) induces a map
\begin{equation}\label{deltal}
\Delta_{\lt}:Y_{n,l}(\sigma) \rightarrow U(\mathfrak{gl}_t) \otimes
Y_{n,l-1}(\dot\sigma),
\end{equation}
where $\dot\sigma$ is as in (\ref{sl}).
Recalling Remark~\ref{sumply}, these maps
$\Delta_{\rt}$ and $\Delta_{\lt}$ are filtered,
so induce homomorphisms
\begin{align}\label{drg}
\gr \Delta_{\rt}:\gr Y_{n,l}(\sigma) &\rightarrow \gr 
(Y_{n,l-1}(\dot\sigma) \otimes 
U(\mathfrak{gl}_t)),\\
\gr \Delta_{\lt}:\gr Y_{n,l}(\sigma) &\rightarrow \gr 
(U(\mathfrak{gl}_t) \otimes Y_{n,l-1}(\dot\sigma))\label{dlg}
\end{align}
of graded algebras.

\begin{Theorem}\label{miuramain}
For any admissible shape $\nu = (\nu_1,\dots,\nu_m)$,
$\gr Y_{n,l}(\sigma)$ is the free commutative algebra
on generators
$\{e_{a,b;i,j}^{(r)}\}_{1 \leq a,b \leq m,
1 \leq i \leq \nu_a, 1 \leq j \leq \nu_b, s_{a,b}(\nu)
< r \leq s_{a,b}(\nu) + p_{\min(a,b)}(\nu)}.$
Moreover, the maps $\gr \Delta_{\rt}$ and $\gr\Delta_{\lt}$ from 
(\ref{drg})--(\ref{dlg})
are injective whenever they are defined,
hence so are the maps (\ref{deltar})--(\ref{deltal})
\end{Theorem}

\begin{proof}
We proceed by induction on the level $l$,
the case $l = 0$ being trivial. Now suppose that $l > 0$
and that the first statement of the theorem has been proved for
all smaller $l$. In view of Lemma~\ref{yak}, it suffices to check
the induction step in the special case that $\nu = (\nu_1,\dots,\nu_m)$
is the minimal admissible shape for $\sigma$.
At least one of the maps $\Delta_{\rt}$ or $\Delta_{\lt}$
is always defined. We will assume without loss
of generality that $\Delta_{\rt}$ is defined; the theorem 
in the other case can be deduced from this one using the
anti-isomorphism $\tau$.
Introduce the following elements
of $\gr (Y_{n,l-1}(\dot \sigma) \otimes U(\mathfrak{gl}_t))$:
\begin{align*}
\dot e_{a,b;i,j}^{(r)} &:= 
\left\{
\begin{array}{ll}
\gr_r \dot D_{a;i,j}^{(r)} \otimes 1&\hbox{if $a = b$,}\\
\gr_r \dot E_{a,b;i,j}^{(r)} \otimes 1&\hbox{if $a < b$,}\\
\gr_r \dot F_{a,b;i,j}^{(r)} \otimes 1&\hbox{if $a > b$.}
\end{array}
\right.
\end{align*}
Also let
$x_{i,j} := \gr_1 1 \otimes e_{i,j}$ for $1 \leq i,j \leq t$.
By Theorem~\ref{baby}(i), Lemma~\ref{higherroots}(i) and Lemma~\ref{yak},
there exist polynomials
$f_{a;i,j}^{(r)}$ in all the variables
$\dot e_{a,b;i,j}^{(r)}$
such that $\gr \Delta_{\rt}$ maps
\begin{align}\label{first}
e_{a,b;i,j}^{(r)} &\mapsto \dot e_{a,b;i,j}^{(r)}\\\intertext{for 
$1 \leq a \leq m, 1 \leq b < m, 1 \leq i\leq \nu_a,
1 \leq j \leq \nu_b$ and $s_{a,b}(\nu) < r \leq s_{a,b}(\nu)+
p_{\min(a,b)}(\nu)$, and}
\label{second}
e_{a,m;i,j}^{(r)} &\mapsto \sum_{k=1}^t \dot e_{a,m;i,k}^{(r-1)}
x_{k,j} + f_{a;i,j}^{(r)}
\end{align}
for $1\leq a \leq m, 1 \leq i \leq \nu_a, 1\leq j \leq \nu_m$
and $s_{a,m}(\nu) < r \leq s_{a,m}(\nu)+p_a(\nu)$,
where $\dot e_{m,m;i,j}^{(0)} := \delta_{i,j}$.
By the induction hypothesis and the PBW theorem for $U(\mathfrak{gl}_t)$, 
the elements
\begin{align*}
&\{x_{i,j}\}_{1 \leq i,j \leq t},\\
&\{\dot e_{a,b;i,j}^{(r)}\}_{1 \leq a \leq m, 1 \leq b < m,
1 \leq i \leq \nu_a, 1 \leq j \leq \nu_b,
s_{a,b}(\nu) < r \leq s_{a,b}(\nu)+p_{\min(a,b)}(\nu)},\\
&\{\dot e_{a,m;i,j}^{(r)}\}_{1 \leq a \leq m,
1 \leq i \leq \nu_a, 1 \leq j \leq \nu_m,
s_{a,m}(\nu)+\delta_{a,m}-1 < r \leq s_{a,m}(\nu)+p_a(\nu)-1}
\end{align*}
and are algebraically independent in 
$\gr (Y_{n,l-1}(\dot \sigma) \otimes U(\mathfrak{gl}_t))$.
Using this and (\ref{first})--(\ref{second}), one verifies 
explicitly that the images of the generators
$$
\{e_{a,b;i,j}^{(r)}\}_{1 \leq a,b \leq m, 
1 \leq i \leq \nu_a, 1 \leq j \leq \nu_b, s_{a,b}(\nu) <
r \leq s_{a,b}(\nu) + p_{\min(a,b)}(\nu)}
$$
of $\gr Y_{n,l}(\sigma)$ from Lemma~\ref{yak}
under the map $\gr \Delta_{\rt}$
are algebraically independent
in $\gr Y_{n,l-1}(\dot \sigma) \otimes U(\mathfrak{gl}_t)$. 
Hence $\gr \Delta_{\rt}$ is
injective and these generators must already be algebraically
independent in $\gr Y_{n,l}(\sigma)$. This completes the proof of
the induction step.
\end{proof}

\begin{Corollary}\label{newpbw}
For any admissible shape $\nu = (\nu_1,\dots,\nu_m)$,
the monomials in the elements
\begin{align*}
&\{D_{a;i,j}^{(r)}\}_{1 \leq a \leq m, 1 \leq i,j 
\leq \nu_a, 0 < r \leq p_a(\nu)},\\
&\{E_{a,b;i,j}^{(r)}\}_{1 \leq a < b \leq m,
1 \leq i \leq \nu_a, 1 \leq j \leq \nu_b, s_{a,b}(\nu) < r \leq
s_{a,b}(\nu)+p_a(\nu)},\\
&\{F_{a,b;i,j}^{(r)}\}_{1 \leq a < b \leq m,
1 \leq j \leq \nu_a, 1 \leq i \leq \nu_b, s_{b,a}(\nu) < r \leq
s_{b,a}(\nu)+p_a(\nu)}
\end{align*}
taken in any fixed order form a basis for $Y_{n,l}(\sigma)$.
\end{Corollary}

\begin{Remark}\label{invsys}\rm
Obviously by their definition there is an inverse system
\begin{equation*}
Y_{n,l}(\sigma) \twoheadleftarrow Y_{n,l+1}(\sigma)
\twoheadleftarrow Y_{n,l+2}(\sigma) \twoheadleftarrow \cdots
\end{equation*}
Moreover, the maps are homomorphisms of filtered algebras
with respect to the
canonical filtrations.
Comparing the basis theorems proved in Corollary~\ref{newpbw}
and Theorem~\ref{pbw3}(iv), it follows that
$Y_n(\sigma)
= \varprojlim Y_{n,l}(\sigma)$
where the inverse limit is taken in the category of filtered algebras.
Hence we can view the shifted Yangian
$Y_n(\sigma)$ as the limiting case $l \rightarrow \infty$
of the shifted Yangian $Y_{n,l}(\sigma)$ of level $l$.
\end{Remark}

\section{Pyramids}\label{swhittaker}

In this section we introduce the combinatorics of pyramids.
(In the more general language of \cite[$\S$4]{EK} 
the things we call pyramids here should be called
even pyramids.)
Suppose to start with that we are given a tuple
$(q_1,q_2,\dots,q_l)$ of positive integers.
We associate a diagram $\pi$ consisting of $q_1$ bricks
stacked in the first (leftmost) column,
$q_2$ bricks stacked in the second column, \dots,
$q_l$ bricks stacked in the $l$th (rightmost) column.
For instance if $l=4$ and
$(q_1,q_2,q_3,q_4) = (4,2,3,3)$, then
$$
\begin{picture}(65, 60)%
\put(0,0){\line(1,0){60}}
\put(0,15){\line(1,0){60}}
\put(0,30){\line(1,0){60}}
\put(0,45){\line(1,0){15}}
\put(0,60){\line(1,0){15}}
\put(30,45){\line(1,0){30}}
\put(0,0){\line(0,1){60}}
\put(15,0){\line(0,1){60}}
\put(30,0){\line(0,1){45}}
\put(45,0){\line(0,1){45}}
\put(60,0){\line(0,1){45}}
\put(7,8){\makebox(0,0){4}}
\put(7,23){\makebox(0,0){3}}
\put(7,38){\makebox(0,0){2}}
\put(7,53){\makebox(0,0){1}}
\put(22,8){\makebox(0,0){6}}
\put(22,23){\makebox(0,0){5}}
\put(37,8){\makebox(0,0){9}}
\put(37,23){\makebox(0,0){8}}
\put(37,38){\makebox(0,0){7}}
\put(52,8){\makebox(0,0){12}}
\put(52,23){\makebox(0,0){11}}
\put(52,38){\makebox(0,0){10}}
\put(-20,27){\makebox(0,0){$\pi =$}}
\end{picture}
$$
Call $\pi$ a {\em pyramid} of {\em level} $l$ 
and {\em height} $\max(q_1,\dots,q_l)$ if each row of the
diagram consists of a single connected horizontal strip.
Equivalently, $\pi$ is a pyramid if the sequence
$(q_1,q_2,\dots,q_l)$ of {\em column heights} is a
unimodal sequence, i.e.
\begin{equation}\label{pyr1}
0 < q_1 \leq \cdots \leq q_k,
\quad
q_{k+1} \geq \cdots \geq q_l > 0
\end{equation} 
for some $0 \leq k \leq l$.
Of course, the above diagram is {\em not} a pyramid, since there
is a gap between the entries $2$ and $7$.

Given a diagram $\pi$ (not necessarily a pyramid), we pick an integer
$n \geq \max(q_1,\dots,q_l)$ and number
the rows of the diagram $1,2,\dots,n$ 
from top to bottom.
Let $p_i$ denote the number of bricks in the $i$th row,
thus defining the tuple $(p_1,\dots,p_n)$ of {\em row lengths} with
\begin{equation}
0 \leq p_1\leq \cdots \leq p_n = l.
\end{equation}
\iffalse
Usually we have in mind that $n$ should
exactly equal $\max(q_1,\dots,q_l)$,
so that {\em either $n=0$ or $p_1 > 0$}, but it is sometimes useful to allow $n$ to be larger.
\fi
Fix also some numbering of the bricks of the diagram by
$1,2,\dots,N$. Usually we have in mind the numbering down columns from
left to right as in the above example but any 
bijective numbering will do.
For $i=1,\dots,N$, let $\row(i)$ and $\col(i)$
denote the {\em row} and {\em column} numbers
of the brick in which $i$ appears, respectively.

Now let $\mathfrak{g}$ denote the Lie algebra $\mathfrak{gl}_N$
over $\C$ and
introduce a $\Z$-grading $\mathfrak{g} = 
\bigoplus_{j \in \Z} \mathfrak{g}_j$ of $\mathfrak{g}$
by declaring that
the $ij$-matrix unit $e_{i,j}$ is of degree $(\col(j)-\col(i))$.
Let
\begin{equation}\label{Ac1}
\mathfrak{p} := \bigoplus_{j \geq 0} \mathfrak{g}_j,\qquad\qquad
\mathfrak{h} := \mathfrak{g}_0,\qquad\qquad
\mathfrak{m} := \bigoplus_{j < 0} \mathfrak{g}_j.
\end{equation}
Thus $\mathfrak{p}$ is a parabolic subalgebra of $\mathfrak{g}$
with Levi factor $\mathfrak{h}
\cong \mathfrak{gl}_{q_1}\oplus\cdots\oplus\mathfrak{gl}_{q_l}$,
and $\mathfrak{m}$
is the nilradical of the opposite parabolic to $\mathfrak{p}$
with respect to $\mathfrak{h}$. 
Let $P, H$ and $M$ be the corresponding closed subgroups of 
the algebraic group $G := GL_N$ over $\C$.
Also let
\begin{equation}\label{edef}
e := \sum_{\substack{1 \leq i,j \leq N \\ \row(i) = \row(j)\\ \col(i)= \col(j)-1}} e_{i,j}.
\end{equation}
By a dimension calculation, one checks that
$e$ belongs to the unique dense orbit of $H$ on $\mathfrak{g}_1$,
regardless of whether the diagram $\pi$ is a pyramid or not.
Moreover, the Jordan block sizes of the matrix $e$ are 
precisely
the lengths of the maximal connected horizontal strips in the
diagram $\pi$, hence $e$ is of Jordan type
$(p_1,\dots,p_n)$ if and only if the diagram $\pi$ is a pyramid.

Recall from \cite[$\S$5.2]{Carter} that an element
$e$ of the nilradical of $\mathfrak{p}$ is called
a {\em Richardson element for $\mathfrak{p}$} if its orbit
under the adjoint action of $P$
is dense in the nilradical of $\mathfrak{p}$;
equivalently,
\begin{equation}\label{rp}
\dim \mathfrak{c}_{\mathfrak{g}}(e) = \dim \mathfrak{h}.
\end{equation}
By another dimension calculation (as observed originally by
Kraft \cite{Kr}) the Jordan type
of a Richardson element for $\mathfrak{p}$ is 
given by the row lengths $(p_1,\dots,p_n)$ of the diagram $\pi$,
again regardless of whether $\pi$ is a pyramid or not.
We say that $\mathfrak{p}$ is a {\em good parabolic}
if $\mathfrak{g}_1$ contains a Richardson element for $\mathfrak{p}$.
Such an element $e$ then clearly belongs both to the
dense orbit of $P$ on the nilradical of $\mathfrak{p}$
and to the dense orbit of $H$ on $\mathfrak{g}_1$.
Hence up to conjugacy, we may assume that $e$ given
by (\ref{edef}). Since its Jordan type must also be
$(p_1,\dots,p_n)$, we obtain the following 
special case of \cite[Lemma 7.2]{Ly}:

\begin{Theorem}\label{pyrclass}
$\mathfrak{p}$ is a good parabolic if and only if the
diagram $\pi$ is a pyramid.
\end{Theorem}

\begin{Remark}\rm
In view of \cite[Theorem 2.1]{EK}, this theorem
implies that there is a bijective map
from pyramids to
conjugacy classes of
even good gradings of $\mathfrak{g} = \mathfrak{gl}_N$,
as defined in the introduction.
The map sends a pyramid
$\pi$ to twice the grading defined here, which is an even good
grading
for the nilpotent matrix $e$ defined by (\ref{edef}).
\end{Remark}

To make the connection with the earlier sections,
we point out that pyramids provide an extremely
convenient way to visualize the data needed to define the
shifted Yangian $Y_{n,l}(\sigma)$ of level $l$.
Indeed, given a shift matrix 
$\sigma = (s_{i,j})_{1 \leq i,j \leq n}$
and a level $l \geq s_{1,n}+s_{n,1}$, 
we define 
$(p_1,\dots,p_n)$ from $p_i := l-s_{i,n}-s_{n,i}$
like in (\ref{pidef}). Now draw a pyramid $\pi$
with $p_i$ bricks on the $i$th row indented
$s_{n,i}$ columns from the left hand edge and
$s_{i,n}$ columns from the right hand edge,
for each $i=1,\dots,n$.
Conversely, given a pyramid $\pi$ of height $\leq n$,
define a shift matrix $\sigma = (s_{i,j})_{1 \leq i,j \leq n}$ 
from the equation
\begin{equation}\label{pyr3}
s_{i,j} = \left\{
\begin{array}{ll}
\#\{c=1,\dots,k\:|\:
i >
n-q_c \geq j\}&\hbox{if $i \geq j$,}\\
\#\{c=k+1,\dots,l\:|\:
i \leq
n-q_c < j\}&\hbox{if $i \leq j$,}
\end{array}\right.
\end{equation}
where 
$(q_1,\dots, q_l)$ are the column heights and
$k$ is chosen as in (\ref{pyr1}).
This definition is independent of the choice of $k$ 
{\em only if} the pyramid $\pi$ is of height {\em exactly} $n$, 
in which case 
$s_{i,j}$ is simply the number of bricks the $i$th row is indented
from the $j$th row at the left edge of the diagram if $i \geq j$,
at the right edge  of the diagram if $i \leq j$.
For example, 
$$
l=7,
\sigma = \left(\begin{array}{llll}
0&1&1&3\\
0&0&0&2\\
1&1&0&2\\
2&2&1&0
\end{array}\right)
\qquad\longleftrightarrow \qquad
n=4,\pi = 
{\begin{picture}(90, 35)%
\put(0,-25){\line(1,0){105}}
\put(0,-10){\line(1,0){105}}
\put(15,5){\line(1,0){60}}
\put(30,20){\line(1,0){45}}
\put(30,35){\line(1,0){30}}
\put(0,-25){\line(0,1){15}}
\put(15,-25){\line(0,1){30}}
\put(30,-25){\line(0,1){60}}
\put(45,-25){\line(0,1){60}}
\put(60,-25){\line(0,1){60}}
\put(75,-25){\line(0,1){45}}
\put(90,-25){\line(0,1){15}}
\put(105,-25){\line(0,1){15}}
\end{picture}}
$$
Now we have a convenient notation to record 
the explicit description of the centralizer 
$\mathfrak{c}_{\mathfrak{g}}(e)$ of the nilpotent element
$e$ associated to a pyramid $\pi$;
see \cite[IV.1.6]{SS}.

\begin{Lemma}\label{centbase}
Let $\pi$ be a pyramid of height $\leq n$ with row
lengths $(p_1,\dots,p_n)$ and associated shift
matrix $\sigma = (s_{i,j})_{1 \leq i,j \leq n}$, and
let $e$ be the nilpotent matrix defined by (\ref{edef}).
For $1 \leq i,j \leq n$ and $r \geq 0$, let
$$
c_{i,j}^{(r)} := \sum_{\substack{1 \leq h,k \leq N \\ 
\row(h) = i, \row(k) = j \\ \col(k) - \col(h) +1 = r}} e_{h,k}.
$$
Then, 
the vectors $\{c_{i,j}^{(r)}\}_{1 \leq i,j \leq n, s_{i,j} < r \leq s_{i,j}+
p_{\min(i,j)}}$ give a basis for
$\mathfrak{c}_{\mathfrak{g}}(e)$.
\end{Lemma}

\section{Finite $W$-algebras}\label{sslice}

Continue with $\mathfrak{g} = \mathfrak{gl}_N$
and $G = GL_N$ acting on $\mathfrak{g}$
by the adjoint action $\Ad$.
Fixing a pyramid $\pi$ with 
bricks numbered $1,\dots,N$, we have the associated $\Z$-grading
on $\mathfrak{g}$ defined as in 
$\S$\ref{swhittaker} so that $e_{i,j}$ is of degree
$(\col(j)-\col(i))$.
Also define $\mathfrak{p},
\mathfrak{h}$, $\mathfrak{m}$ and the nilpotent
matrix $e \in \mathfrak{g}_1$ according to
(\ref{Ac1})--(\ref{edef}).
Define $h \in \mathfrak{g}_0$ and $f \in \mathfrak{g}_{-1}$ 
to be the unique matrices so that
$(e,h,f)$ is an $\mathfrak{sl}_2$-triple in $\mathfrak{g}$.
Let $\mathfrak{m}^\perp$ and $\mathfrak{p}^\perp$
denote the orthogonal complements of $\mathfrak{m}$ and $\mathfrak{p}$ 
with respect to the trace form
$(.,.)$ on $\mathfrak{g}$,
i.e. $\mathfrak{m}^\perp = \bigoplus_{j \leq 0} \mathfrak{g}_j$
and $\mathfrak{p}^\perp = \bigoplus_{j > 0} \mathfrak{g}_j$.

\begin{Lemma}\label{handy}
In the above notation,
\begin{itemize}
\item[(i)] 
$\mathfrak{m}^\perp = [\mathfrak{m},e] \oplus \mathfrak{c}_{\mathfrak{g}}(f)$;
\item[(ii)] 
$\mathfrak{p} = [\mathfrak{p}^\perp,f] \oplus \mathfrak{c}_{\mathfrak{g}}(e)$;
\item[(iii)]
$\mathfrak{c}_{\mathfrak{g}}(f)^\perp = \mathfrak{m}\oplus
[\mathfrak{p}^\perp,f]$.
\end{itemize}
\end{Lemma}

\begin{proof}
Since $e$ is a Richardson element for $\mathfrak{p}$,
we know that $\dim \mathfrak{c}_{\mathfrak{g}}(e)
= \dim \mathfrak{h}$ and
$\mathfrak{c}_{\mathfrak{g}}(e) \subseteq \mathfrak{p}$.
Hence the map $\mathfrak{m} \rightarrow [\mathfrak{m},e],
x \mapsto [x,e]$ is a bijection.
Similarly, $\dim \mathfrak{c}_{\mathfrak{g}}(f)
= \dim \mathfrak{h}$ and
$\mathfrak{c}_{\mathfrak{g}}(f) \subseteq \mathfrak{m}^\perp$.
Hence the map $\mathfrak{p}^\perp \rightarrow [\mathfrak{p}^\perp,f],
y \mapsto [y,f]$ is a bijection.
So
\begin{align*}
\dim [\mathfrak{m},e] + \dim \mathfrak{c}_{\mathfrak{g}}(f)
&= \dim \mathfrak{m} + \dim \mathfrak{h}
= \dim \mathfrak{m}^{\perp},\\
\dim [\mathfrak{p}^\perp,f] + \dim \mathfrak{c}_{\mathfrak{g}}(e)
&= \dim \mathfrak{p}^\perp + \dim \mathfrak{h}
= \dim \mathfrak{p}.
\end{align*}
Also by $\mathfrak{sl}_2$ theory,
$[\mathfrak{m},e] \cap \mathfrak{c}_{\mathfrak{g}}(f)= 
[\mathfrak{p}^\perp,f] \cap \mathfrak{c}_{\mathfrak{g}}(e) = 0$.
Parts (i) and (ii) follow.
For (iii), we have that
$$
\dim \mathfrak{m} + \dim [\mathfrak{p}^\perp,f]
= \dim \mathfrak{g} - \dim \mathfrak{h} = \dim \mathfrak{c}_{\mathfrak{g}}(f)^\perp,
$$
and clearly $\mathfrak{m} \cap [\mathfrak{p}^\perp,f] = 0$,
so we just have to check
that $\mathfrak{m} \subseteq \mathfrak{c}_{\mathfrak{g}}(f)^\perp$
and that
$[\mathfrak{p}^\perp,f] 
\subseteq \mathfrak{c}_{\mathfrak{g}}(f)^\perp$.
The former statement is true since 
$\mathfrak{c}_{\mathfrak{g}}(f) \subseteq \mathfrak{m}^\perp$.
For the latter, note for any $x \in \mathfrak{p}^\perp$
and $y \in \mathfrak{c}_{\mathfrak{g}}(f)$ that
$([x,f],y) = (x,[f,y]) =0$ by the invariance of the trace
form.
\end{proof}

To recall the definition of the algebra $W(\pi)$ from the introduction, 
let $\chi:\mathfrak{m} \rightarrow \C$ denote
the representation mapping $x \mapsto (x,e)$ for each $x \in \mathfrak{m}$.
Let $I_\chi$ denote the kernel of the corresponding
algebra homomorphism $U(\mathfrak{m}) \rightarrow \C$
and $\pr_\chi:U(\mathfrak{g}) \rightarrow U(\mathfrak{p})$
be the projection along the direct sum decomposition
\begin{equation}
U(\mathfrak{g}) = U(\mathfrak{p}) \oplus U(\mathfrak{g}) I_\chi.
\end{equation}
For $x \in \mathfrak{m}$ and $y \in U(\mathfrak{p})$, we set
$x \cdot y := \pr_\chi([x,y]) = \pr_\chi(xy) - \chi(x)y$,
to define the {\em twisted action} 
of $\mathfrak{m}$ on $U(\mathfrak{p})$.
Then, the {\em finite $W$-algebra}
$W(\pi)$ associated to the pyramid $\pi$ is the space $U(\mathfrak{p})^{\mathfrak m}$
of twisted $\mathfrak{m}$-invariants in $U(\mathfrak{p})$.
Equivalently,
\begin{equation}
W(\pi)= \{y \in U(\mathfrak{p})\:|\:
(x - \chi(x))y \in U(\mathfrak{g})I_\chi
\hbox{ for all }x \in \mathfrak{m}\},
\end{equation}
from which it follows easily that $W(\pi)$ is actually
a subalgebra of $U(\mathfrak{p})$. 
For example, in the special case that $\pi$ consists of a single column, 
we obviously have that
$\mathfrak{p} = \mathfrak{g}$, $\mathfrak{m} = 0$
and $e = 0$, hence 
$W(\pi) = U(\mathfrak{gl}_N)$.

We need the
{\em Kazhdan filtration}
$\cdots\subseteq 
\F_d U(\mathfrak{g}) \subseteq \F_{d+1} U(\mathfrak{g}) \subseteq \cdots$
of $U(\mathfrak{g})$ defined
by declaring that the
generator $e_{i,j}$ is of degree
\begin{equation}\label{degdef}
\deg(e_{i,j}) := \col(j)-\col(i)+1
\end{equation}
for each $1 \leq i,j \leq N$, i.e. $\F_d U(\mathfrak{g})$ is
spanned by all monomials 
$e_{i_1,j_1} \cdots e_{i_m,j_m}$ for $m \geq 0$ and
$\deg(e_{i_1,j_1})+\cdots+\deg(e_{i_m,j_m}) \leq d$.
The adjoint action of $\mathfrak{g}$
on $U(\mathfrak{g})$ is filtered in the sense that
$[\mathfrak{g}_j, \F_d U(\mathfrak{g})]
\subseteq \F_{d+j} U(\mathfrak{g})$ for each $j,d \in \Z$.
Hence the associated graded algebra
$\gr U(\mathfrak{g})$ is naturally a graded
$\mathfrak{g}$-module.
Of course by the PBW theorem we can identify
$\gr U(\mathfrak{g})$ 
with the symmetric algebra $S(\mathfrak{g})$ with the usual
$\mathfrak{g}$-action, but viewed as a 
graded algebra via the {\em Kazhdan grading} also defined by (\ref{degdef}).

There are induced Kazhdan filtrations of the subalgebras
$U(\mathfrak p)$  and $W(\pi)$ defined by
setting $\F_d U(\mathfrak{p}) := \pr_\chi (\F_d U(\mathfrak{g}))$
and $\F_d W(\pi) := W(\pi) \cap \F_d U(\mathfrak{p})$.
This time we have that
$\F_d  U(\mathfrak{p})= \F_d W(\pi) = 0$ for all $d < 0$.
The projection $\pr_\chi:U(\mathfrak{g}) \rightarrow U(\mathfrak{p})$ 
is filtered, hence so is the 
twisted action of $\mathfrak{m}$ on $U(\mathfrak{p})$,
and
$\gr \pr_\chi:\gr U(\mathfrak{g}) \rightarrow \gr U(\mathfrak{p})$
is a graded $\mathfrak{m}$-module homomorphism.
If we identify $\gr U(\mathfrak{p})$ with $S(\mathfrak{p})$
graded via (\ref{degdef}),
then we can describe the
map $\gr \pr_\chi$ simply as the algebra homomorphism
$p:S(\mathfrak{g}) \rightarrow S(\mathfrak{p})$
that acts as the identity on elements of $\mathfrak{p}$
and sends $x \in \mathfrak{m}$ to $\chi(x)$.
This explains the top left hand quadrant of the following 
commutative diagram:
\begin{equation}\label{thed}
\begin{CD}
\gr U(\mathfrak{g}) & @= &S(\mathfrak{g}) &@>\sim>\alpha> & \C[\mathfrak g]\\
@V\gr\pr_\chi VV &&@VpVV&& @VrVV\\
\gr U(\mathfrak{p}) &@= &S(\mathfrak{p}) &@>\sim>\beta> & \C[e+\mathfrak{m}^\perp]&\:\:\:\hookleftarrow\:\:\:\: &\C[e+\mathfrak{m}^\perp]^M\\
@A\gr iAA&&@VqVV&&@VsVV\\
\gr W(\pi)&@>>\theta>&S(\mathfrak{c}_{\mathfrak{g}}(e))&@>\sim>\gamma>&\C[e+\mathfrak{c}_{\mathfrak{g}}(f)]
\end{CD}
\begin{picture}(0,0)
\put(-75,-24){\line(1,1){13}}
\put(-75,-24){\makebox(0,0){$\swarrow$}}
\put(-77,-18){\makebox(0,0){$_{\sim}$}}
\end{picture}
\end{equation}
To make sense of the bottom left hand quadrant,
let $i:W(\pi) \rightarrow U(\mathfrak{p})$ be the inclusion,
and let
$q:S(\mathfrak{p}) \rightarrow S(\mathfrak{c}_{\mathfrak{g}}(e))$
be the algebra homomorphism induced
by the projection $\mathfrak{p}\rightarrow\mathfrak{c}_{\mathfrak{g}}(e)$
along the direct sum decomposition from Lemma~\ref{handy}(ii).
Note both $[\mathfrak{p}^\perp,f]$ and $\mathfrak{c}_{\mathfrak{g}}(e)$
are graded subspaces of $\mathfrak{g}$; in the latter case
one way to see this is by Lemma~\ref{centbase} 
since the basis element $c_{i,j}^{(r)}$
there is clearly of degree $r$ with respect to the Kazhdan grading.
Hence the map $q$ is a graded map. Finally let
$\theta:\gr W(\pi) \rightarrow S(\mathfrak{c}_{\mathfrak{g}}(e))$ 
be the composite $q \circ \gr i$.

Now we turn our attention to the right hand half of the diagram.
Imitating \cite[$\S$2.1]{GG}, 
let $\gamma:\C^\times \rightarrow G$ be the group homomorphism
defined by letting
\begin{equation}
\gamma(t):=
\diag(t^{-\col(1)}, t^{-\col(2)}, \cdots, t^{-\col(N)}) \in GL_N
\end{equation}
for each $t \in \C^\times$.
So $\Ad \gamma(t)$ acts on $\mathfrak{g}_j$ by the scalar $t^{j}$
for each $t \in \C^\times$ and $j \in \Z$.
Introduce a linear action $\rho$
of $\C^\times$ on the variety $\mathfrak{g}$
by letting 
\begin{equation}\label{rho1}
\rho(t)(x) := t^{-1} \Ad \gamma(t)(x)
\end{equation}
for each $t \in \C^\times, x \in \mathfrak{g}$.
Thus, $\rho(t) (e_{i,j}) = t^{-\deg(e_{j,i})} e_{i,j}$
for each $1 \leq i,j \leq N$. In particular,
$\rho(t)(e) = e$ for every $t \in \C^\times$ and
\begin{equation}\label{lim1}
\lim_{t \rightarrow \infty} \rho(t)(x) = 0
\end{equation}
for all $x \in \mathfrak{m}^{\perp}$.
We get an induced action $\bar\rho$
of $\C^\times$ on the coordinate
algebra $\C[\mathfrak{g}]$ with
\begin{equation}
(\bar\rho(t) (f))(x) = f(\rho(t^{-1}) (x))\end{equation}
for each $f \in \C[\mathfrak{g}], x \in \mathfrak{g}$.
Using this, we define a grading on $\C[\mathfrak{g}]$
by declaring that $f \in \C[\mathfrak{g}]$ is of degree $j$
if $\bar\rho(t) (f) = t^{j} f$ for each $t \in \C^\times$.
Since $e$ is a $\rho$-fixed point,
it is easy to see that each $\rho(t)$ leaves 
both $e + \mathfrak{m}^\perp$ and
the Slodowy slice $e + \mathfrak{c}_{\mathfrak{g}}(f)$
invariant. So just like for $\C[\mathfrak{g}]$,
we get induced actions $\bar\rho$ of $\C^\times$
on the coordinate algebras $\C[e+\mathfrak{m}^\perp]$
and $\C[e+\mathfrak{c}_{\mathfrak{g}}(f)]$
which we use to introduce a grading on these algebras.
The natural restriction maps
$r:\C[\mathfrak{g}] \twoheadrightarrow
\C[e+\mathfrak{m}^\perp]$ and
$s:\C[e+\mathfrak{m}^\perp] 
\twoheadrightarrow \C[e+\mathfrak{c}_{\mathfrak{g}}(f)]$
are $\bar\rho$-equivariant, hence are graded.

The trace form defines an isomorphism
$\alpha:S(\mathfrak{g}) \rightarrow \C[\mathfrak{g}]$
of graded algebras. It maps the
kernel of the homomorphism $p:S(\mathfrak{g})
\rightarrow S(\mathfrak{p})$ isomorphically onto the annihilator
in $\C[\mathfrak{g}]$ of the closed subvariety $e+\mathfrak{m}^\perp$.
Hence $\alpha$ 
induces a graded algebra
isomorphism $\beta:S(\mathfrak{p}) \rightarrow 
\C[e+\mathfrak{m}^\perp]$.
Similarly by Lemma~\ref{handy}(ii),(iii),
$\beta$ maps the kernel of the homomorphism
$q:S(\mathfrak{p}) \rightarrow S(\mathfrak{c}_{\mathfrak{g}}(e))$
isomorphically onto the annihilator in $\C[e+\mathfrak{m}^\perp]$
of $e+\mathfrak{c}_{\mathfrak{g}}(f)$, hence induces
a graded algebra isomorphism $\gamma:S(\mathfrak{c}_{\mathfrak{g}}(e))
\rightarrow \C[e+\mathfrak{c}_{\mathfrak{g}}(f)]$.

To complete the picture,
we need to introduce an action of the subgroup $M$ of $G$
corresponding to the subalgebra $\mathfrak{m}$ of $\mathfrak{g}$.
The adjoint action of $G$ on $\mathfrak{g}$
induces an action of $G$ on $\C[\mathfrak{g}]$
by algebra automorphisms, such that the derived action 
of $\mathfrak{g}$
corresponds under the isomorphism $\alpha$ to the
usual $\mathfrak{g}$-action on $S(\mathfrak{g})$.
The subgroup $M$ of $G$ leaves 
$e+\mathfrak{m}^\perp$ invariant, hence we get induced
actions of $M$ and $\mathfrak{m}$ on $\C[e+\mathfrak{m}^\perp]$,
such that
the action of $\mathfrak{m}$
agrees under the isomorphism $\beta$ with the
twisted action of $\mathfrak{m}$ on $S(\mathfrak{p})$.
In particular since $S(\mathfrak{p})$
is a graded $\mathfrak{m}$-module, it follows that the
space $\C[e+\mathfrak{m}^\perp]^{\mathfrak{m}}
=\C[e+\mathfrak{m}^\perp]^M$ of $\mathfrak{m}$-invariants/$M$-fixed points
is a {\em graded} subalgebra of $\C[e+\mathfrak{m}^\perp]$.
One can see this directly by introducing an action
$\rho$ of $\C^\times$ on $M \times (e+\mathfrak{m}^\perp)$
defined by
\begin{equation}\label{rho2}
\rho(t) (m,x) := (\gamma(t) m \gamma(t)^{-1}, \rho(t)(x))
\end{equation}
for all $m \in M, x \in e+\mathfrak{m}^\perp$ and $t \in \C^\times$.
The following calculation checks that the 
adjoint action $\varphi:M \times (e+\mathfrak{m}^\perp)
\rightarrow e+\mathfrak{m}^\perp$ is $\rho$-equivariant:
\begin{align*}
\Ad (\gamma(t) m \gamma(t)^{-1}) \rho(t) (x)
&=
t^{-1} \Ad \gamma(t) \Ad m \Ad \gamma(t)^{-1} \Ad \gamma(t) (x)\\
&= 
t^{-1} \Ad \gamma(t) \Ad m (x)
=
\rho(t) (\Ad m(x)).
\end{align*}
This now 
implies that the space of $M$-fixed points in $\C[e+\mathfrak{m}^\perp]$
is invariant under each $\rho(t)$, hence is graded.
Let us also note that
\begin{equation}\label{lim2}
\lim_{t \rightarrow \infty} (\gamma(t) m \gamma(t)^{-1}) = 1
\end{equation}
for all $m \in M$.

It just remains to prove that the composite
of the inclusion
$\C[e+\mathfrak{m}^\perp]^M \hookrightarrow \C[e+\mathfrak{m}^\perp]$
and the projection $s:\C[e+\mathfrak{m}^\perp] \twoheadrightarrow
\C[e+\mathfrak{c}_{\mathfrak{g}}(f)]$ is an isomorphism.
This follows from the following key result,
which is due originally to Kostant \cite[Theorem 1.2]{K},
and is proved in this generality in \cite[Theorem 1.2]{Ly}
using Zariski's Main Theorem.
The alternative argument sketched here is due to Gan and Ginzburg; see 
\cite[Lemma 2.1]{GG}.

\begin{Theorem}\label{adac}
The adjoint action 
$\varphi:M \times (e + \mathfrak{c}_{\mathfrak{g}}(f))
\rightarrow e + \mathfrak{m}^\perp$ is an isomorphism of affine varieties.
\end{Theorem}

\begin{proof}
We just verify the hypothesis needed to 
apply the general result from the proof of 
\cite[Lemma 2.1]{GG}:
{\em An equivariant morphism $\varphi:X_1 \rightarrow X_2$
of smooth affine $\C^\times$-varieties with contracting
$\C^\times$-actions which induces an isomorphism between the tangent
spaces at the $\C^\times$-fixed points must be an isomorphism}.
We have already defined actions $\rho$
of $\C^\times$ on $e+\mathfrak{m}^\perp$ (\ref{rho1})
and on $M \times (e+\mathfrak{c}_{\mathfrak{g}}(f))$ (\ref{rho2})
and checked that $\varphi$ is $\rho$-equivariant.
By (\ref{lim1}) and (\ref{lim2}), 
we have that
$\lim_{t \rightarrow \infty} \rho(t) (m,x) = (1,e)$
for each $(m,x) \in M \times (e+\mathfrak{c}_{\mathfrak{g}}(f))$
and that
$\lim_{t \rightarrow \infty} \rho(t)(x) = e$
for each $x \in e+\mathfrak{m}^\perp$. Hence
the $\C^\times$-actions are both contracting.
So finally we need to check that the differential 
$d\varphi_{(1,e)}$
is an isomorphism between the 
tangent spaces $T_{(1,e)}(M \times(e + \mathfrak{c}_{\mathfrak{g}}(f)))$
and $T_e(e+\mathfrak{m}^\perp)$.
But if we identify the tangent spaces with
$\mathfrak{m}\oplus \mathfrak{c}_{\mathfrak{g}}(f)$
and $\mathfrak{m}^\perp$ respectively, then the differential is the map
$(x,y) \mapsto [x,e]+y$. Hence it is an isomorphism by Lemma~\ref{handy}(i).
\end{proof}

The crucial thing that we can now read off from the
diagram (\ref{thed}) is the following:

\begin{Corollary}\label{phewy}
The map $\theta:\gr W(\pi) \rightarrow S(\mathfrak{c}_{\mathfrak{g}}(e))$
is an injective graded algebra homomorphism.
\end{Corollary}

\begin{proof}
Clearly $\gr i$ maps
$\gr W(\pi)$ injectively into the space of $\mathfrak{m}$-invariants
in $\gr U(\mathfrak{p})$.
Hence, $\beta \circ \gr i$ maps $\gr W(\pi)$
into $\C[e+\mathfrak{m}^\perp]^M$.
Hence $s \circ \beta \circ \gr i$ is injective.
But this is $\gamma \circ \theta$ by the commmutativity
of the diagram, hence $\theta$ is injective too.
\end{proof}

\begin{Remark}\rm
It is known by \cite[Theorem 2.3]{Ly} that the map
$\theta$ is actually an {\em isomorphism},
hence $\gr W(\pi)$ is isomorphic to the coordinate algebra
$\C[e+\mathfrak{c}_{\mathfrak{g}}(f)]$ of the Slodowy slice
as stated in the introduction.
A quicker proof can also be given by 
following the arguments of \cite[$\S$5]{GG}.
However we do not need to use this fact yet, and we will be able
deduce it
later on as a consequence of the main result of the article; see
Corollary~\ref{onto}.
\end{Remark}

\ifcenters@
\begin{Remark}\label{centinj}\rm
The restriction of the map $\pr_\chi:U(\mathfrak{g}) 
\rightarrow U(\mathfrak{p})$ to 
$Z(U(\mathfrak{g}))$ defines an injective
algebra homomorphism $\psi:Z(U(\mathfrak{g})) \hookrightarrow
W(\pi)$ whose image is contained in the center $Z(W(\pi))$.
The fact that $\psi$ is an algebra homomorphism 
with image contained in $Z(W(\pi))$ is easiest to
see using the definition of
$W(\pi)$ as the endomorphism algebra $\End_{U(\mathfrak{g})}(Q_\chi)$
given in the introduction, since in those terms
$\psi$ is just the representation
$Z(U(\mathfrak{g})) \rightarrow \End_{\C}(Q_\chi)$ 
of $Z(U(\mathfrak{g}))$ on the module $Q_\chi$.
The fact that $\psi$ is injective
is proved in \cite[6.2]{P} or \cite[Proposition 2.6]{Ly}
by observing that the (injective) Harish-Chandra homomorphism
factors through the map $\psi$.
In \cite[$\S$6]{BK2} we show moreover using some basic facts 
about the representation theory of $W(\pi)$ (some of 
the proofs of which depend on
knowing the main result of the article below)
that the image of $\psi$ is actually
{\em equal to} $Z(W(\pi))$. 
Hence, $\psi:Z(U(\mathfrak{g})) \rightarrow Z(W(\pi))$
is actually an {\em isomorphism}.
\end{Remark}
\fi

\section{Invariants}\label{sinvariants}

In this section we define some remarkable elements of $U(\mathfrak p)$,
many of which will eventually 
turn out to be $\mathfrak m$-invariant, i.e. to
belong to the subalgebra $W(\pi)$.
Letting $(q_1,\dots,q_{l})$ denote the column heights of our fixed
pyramid $\pi$, pick an integer $n \geq \max(q_1,\dots,q_l)$.
Define $\rho = (\rho_1,\dots,\rho_l)$ by setting
\begin{equation}\label{rhodef}
\rho_r := n-q_{r} - q_{r+1} -\cdots-q_l
\end{equation}
for each $r=1,\dots,l$.
For $1 \leq i,j \leq N$, define
\begin{equation}\label{etildedef}
\tilde e_{i,j} := (-1)^{\col(j)-\col(i)} 
(e_{i,j} + \delta_{i,j} \rho_{\col(i)}),
\end{equation}
so
\begin{equation}\label{erel}
[\tilde e_{i,j}, \tilde e_{h,k}] =
(\tilde{e}_{i,k} - \delta_{i,k} \rho_{\col(i)})\delta_{h,j} 
-
\delta_{i,k} (\tilde e_{h,j} - \delta_{h,j} \rho_{\col(j)}).
\end{equation}
Let us also spell out the effect of the homorphism
$U(\mathfrak{m}) \rightarrow \C$ induced by the
character $\chi$: 
we have that
\begin{equation}\label{chidef}
\tilde e_{i,j} \mapsto \left\{
\begin{array}{ll}
-1&\hbox{if $\row(i)=\row(j)$ and $\col(i) = \col(j)+1$;}\\
0&\hbox{otherwise.}
\end{array}\right.
\end{equation}

For $1 \leq i,j \leq n$ and signs
$\sigma_1,\dots,\sigma_n \in \{\pm\}$,
we let
$T_{i,j;\sigma_1,\dots,\sigma_n}^{(0)} := \delta_{i,j} \sigma_i$
and for $r \geq 1$ define
\begin{equation}\label{thedef}
T_{i,j;\sigma_1,\dots,\sigma_n}^{(r)}
:=
\sum_{s = 1}^r
\sum_{\substack{i_1,\dots,i_s\\j_1,\dots,j_s}}
\sigma_{\row(j_1)} \cdots \sigma_{\row(j_{s-1})}
 \tilde e_{i_1,j_1} \cdots \tilde e_{i_s,j_s}
\end{equation}
where the second sum is over all $1 \leq i_1,\dots,i_s,j_1,\dots,j_s \leq N$
such that
\begin{itemize}
\item[(1)] $\deg(e_{i_1,j_1})+\cdots+\deg(e_{i_s,j_s}) = r$ (recall (\ref{degdef}));
\item[(2)] $\col(i_t) \leq \col(j_t)$ for each $t=1,\dots,s$;
\item[(3)] if $\sigma_{\row(j_t)} = +$ then
$\col(j_t) < \col(i_{t+1})$ for each
$t=1,\dots,s-1$;
\item[(4)]
if $\sigma_{\row(j_t)} = -$ then $\col(j_t) \geq \col(i_{t+1})$
for each
$t=1,\dots,s-1$;
\item[(5)] $\row(i_1)=i$, $\row(j_s) = j$;
\item[(6)]
$\row(j_t)=\row(i_{t+1})$ for each $t=1,\dots,s-1$.
\end{itemize}
Note the assumptions (1) and (2) imply
that $T_{i,j;\sigma_1,\dots,\sigma_n}^{(r)}$
belongs to $\F_r U(\mathfrak p)$.
For $x=0,1,\dots,n$, let $T_{i,j;x}^{(r)}$
denote $T_{i,j;\sigma_1,\dots,\sigma_n}^{(r)}$ 
in the special case that $\sigma_1=\cdots=\sigma_x = -$,
$\sigma_{x+1}=\cdots=\sigma_n = +$.
Define
\begin{equation}\label{earlier}
T_{i,j;x}(u) := \sum_{r \geq 0} T_{i,j;x}^{(r)} u^{-r}
\in U(\mathfrak p) [[u^{-1}]].
\end{equation}
Since this is the most critical definition in the entire paper,
let us give some simple examples.

\begin{Example}\label{eggy}\rm
For any $1 \leq i,j \leq n$ and $x = 0,1,\dots,n$,
\begin{align*}
T_{i,j;x}^{(1)} &= 
\sum_{\substack{1 \leq h,k \leq N \\ \row(h) = i, \row(k) = j
\\ \col(h) = \col(k)}} \tilde e_{h,k},\\
T_{i,j;x}^{(2)} &= 
\sum_{\substack{1 \leq h,k \leq N \\ \row(h) = i, \row(k) = j \\
\col(h) = \col(k)-1}} \tilde e_{h,k}
&&\!\!\!\!\!-\!\!\!\!\! \sum_{\substack{1 \leq h_1,h_2,k_1,k_2 \leq N\\
\row(h_1) = i, \row(k_1)=\row(h_2) \leq x, \row(k_2)=j \\
\col(h_1)=\col(k_1) \geq \col(h_2)=\col(k_2),
}} \!\!\!\!\!\!\!\!\!\!\!\!\!\!\!\!\!\!\!\!\tilde e_{h_1,k_1} \tilde e_{h_2,k_2}\\
&&&\!\!\!\!\!+\!\!\!\!\! \sum_{\substack{1 \leq h_1,h_2,k_1,k_2 \leq N\\
\row(h_1) = i, \row(k_1)=\row(h_2) > x, \row(k_2)=j \\
\col(h_1)=\col(k_1) < \col(h_2)=\col(k_2),
}} \!\!\!\!\!\!\!\!\!\!\!\!\!\!\!\!\!\!\!\!\tilde e_{h_1,k_1} \tilde e_{h_2,k_2}.
\end{align*}
\end{Example}

\begin{Lemma}\label{traninv}
Suppose $0 \leq x < y \leq n$.
\begin{itemize}
\item[(i)] If $x < i \leq y$ and $y < j \leq n$ then
$$
T_{i,j;x}(u) = \sum_{k=x+1}^y T_{i,k;x}(u) T_{k,j;y}(u).
$$
\item[(ii)] If $y < i \leq n$ and $x < j \leq y$ then
$$
T_{i,j;x}(u) = \sum_{k=x+1}^y T_{i,k;y}(u) T_{k,j;x}(u).
$$
\item[(iii)] If $y < i,j \leq n$ then
$$
T_{i,j;x}(u) = T_{i,j;y}(u)
+ \sum_{k,l=x+1}^y T_{i,k;y}(u) T_{k,l;x}(u) T_{l,j;y}(u).
$$
\item[(iv)]
If $x < i,j \leq y$ then
$$
\sum_{k=x+1}^y T_{i,k;x}(u) T_{k,j;y}(u) = -\delta_{i,j}.
$$
\end{itemize}
\end{Lemma}

\begin{proof}
Let $\xi_{i,j} :=\tilde e_{i,j} u^{-\deg (e_{i,j})}$ for short.

(i)
The right hand side of the formula in (i) is a sum of monomials
of the form
\begin{equation}\label{rhs1}
\pm (\xi_{i_1,j_1}\cdots \xi_{i_r,j_r})(\xi_{k_1,l_1}\cdots \xi_{k_s,l_s}),
\end{equation}
for various $r \geq 0, s \geq 1$, where
$\pm \xi_{i_1,j_1}\cdots \xi_{i_r,j_r}$ appears in
$T_{i,k;x}(u)$ and 
$\pm \xi_{k_1,l_1}\cdots \xi_{k_s,l_s}$ appears in
$T_{k,j;y}(u)$ for some $x < k \leq y$.
Let $X$ be the sum of all such monomials for which
$\col(j_r) < \col(k_1)$ if $r > 0$
and for which $\row(l_t) \notin \{x+1,\dots,y\}$ for all
$1 \leq t < s$.
Let $Y$ be the sum of all remaining monomials.
Thus, the right hand side of the formula in (i)
is equal to $X+Y$. Now we proceed to show that
$X = T_{i,j;x}(u)$ and that $Y = 0$.

First consider $X$.
Take a monomial of the form (\ref{rhs1}) appearing in $X$,
so $\col(j_r) < \col(k_1)$ if $r > 0$ 
and $\row(l_t) \notin \{x+1,\dots,y\}$ for all
$1 \leq t < s$.
It is easy to see that 
this monomial also appears in 
the expansion of $T_{i,j;x}(u)$, with the same sign.
Moreover, the monomial (\ref{rhs1}) appears in $X$ exactly once:
otherwise we would be able to obtain this monomial 
in the expansion of $X$ in another way by splitting it {\em either} as
$$
\pm (\xi_{i_1,j_1}\cdots \xi_{i_r,j_r} \xi_{k_1,l_1}\cdots \xi_{k_t,l_t})
(\xi_{k_{t+1},l_{t+1}}\cdots \xi_{k_s,l_s})
$$
for some $1 \leq t < s$ 
so that, writing $h = \row(l_t)$,
$\pm\xi_{i_1,j_1}\cdots \xi_{i_r,j_r} \xi_{k_1,l_1}\cdots \xi_{k_t,l_t}$
appears in $T_{i,h;x}(u)$ and
$\pm\xi_{k_{t+1},l_{t+1}}\cdots \xi_{k_s,l_s}$
appears in $T_{h,j;y}(u)$,
{\em or} as
$$
\pm (\xi_{i_1,j_1} \cdots \xi_{i_{u-1},j_{u-1}})(\xi_{i_u,j_u}
\cdots \xi_{i_r,j_r} \xi_{k_1,l_1}\cdots\xi_{k_s,l_s})
$$
for some $1 \leq u \leq r$
so that, writing $h = \row(i_u)$,
$\xi_{i_1,j_1} \cdots \xi_{i_{u-1},j_{u-1}}$
appears in $T_{i,h;x}(u)$ and
$\pm\xi_{i_u,j_u}
\cdots \xi_{i_r,j_r} \xi_{k_1,l_1}\cdots\xi_{k_s,l_s}$
appears in
$T_{h,j;y}(u)$.
The former case does not happen since 
we have that $h \notin\{x+1,\dots,y\}$
by assumption.
The latter case does not happen since $\col(j_r) < \col(k_1)$
contrary to the definition of the monomials arising in
$T_{h,j;y}(u)$.
To complete the proof that $X = T_{i,j;x}(u)$, 
we need to show
that every monomial
$\pm \xi_{p_1,q_1}\cdots \xi_{p_u,q_u}$ appearing in
$T_{i,j;x}(u)$ also appears in $X$.
Take $r \geq 0$ to be the maximal index such that
$\pm \xi_{p_1,q_1}\cdots \xi_{p_r,q_r}$ appears in $T_{i,k;x}(u)$
for some $x < k \leq y$; such an $r$ exists as for $r = 0$
the monomial $1$ appears in $T_{i,i;x}(u)$.
It remains to observe that $\col(q_r) < \col(p_{r+1})$
and that
$\row(q_t) \notin \{x+1,\dots,y\}$ for all $r < t < u$, for otherwise
$r$ could be made bigger. In particular this means that
$\pm \xi_{p_{r+1},q_{r+1}} \cdots \xi_{p_u,q_u}$ appears in
$T_{k,j;y}(u)$. Hence if we split our monomial as
$\pm (\xi_{p_1,q_1}\cdots \xi_{p_r,q_r})(\xi_{p_{r+1},q_{r+1}} \cdots \xi_{p_u,q_u})$,
we have something of the form (\ref{rhs1}) that appears in $X$.

Now consider $Y$. 
Take a monomial of the form (\ref{rhs1}) appearing in $Y$, i.e.
either $r >0$ and $\col(j_r) \geq \col(k_1)$
or there is some $1 \leq t < s$ such that
$x < \row(l_t) \leq y$.
We show that this monomial appears exactly twice in
the expansion of $Y$, with opposite signs.
There are two cases.

Suppose first that $r >0$ and that $\col(j_r) \geq \col(k_1)$.
Let $t \leq r$ be the maximal index such that
$x < \row(i_t) \leq y$; such a $t$ exists since $\row(i_1)=i$
and $x < i \leq y$.
Let $h := \row(i_t)$.
Then the monomial $\pm \xi_{i_1,j_1}\cdots \xi_{i_{t-1},j_{t-1}}$
appears in $T_{i,h;x}(u)$ and the monomial
$\pm \xi_{i_t,j_t} \cdots \xi_{i_r,j_r}\xi_{k_1,l_1}\cdots\xi_{k_s,l_s}$
appears in $T_{h,j;y}(u)$.
Moreover, using the facts that $\col(j_{t-1}) < \col(i_t)$
if $t > 1$,
$\col(j_r) \geq \col(k_1)$, and the maximality of the choice of $t$,
we see that
$$
\pm(\xi_{i_1,j_1}\cdots \xi_{i_r,j_r})(\xi_{k_1,l_1}\cdots \xi_{k_s,l_s}),\:\:
\pm(\xi_{i_1,j_1}\cdots \xi_{i_{t-1},j_{t-1}})(\xi_{i_t,j_t},\cdots,\xi_{i_r,j_r}\xi_{k_1,l_1}\cdots\xi_{k_s,l_s})
$$
are the only ways to split the monomial (\ref{rhs1}) so the left hand term
occurs in $T_{i,g;x}(u)$ and the right hand term appears in 
$T_{g,j;y}(u)$ for some $x < g \leq y$. It just remains to check that the
two have opposite signs.

Suppose instead that $\col(j_r) < \col(k_1)$ if $r > 0$
and that $x < \row(l_t) \leq y$ for some $1 \leq t < s$.
Choose the minimal such $t$ and 
let $h := \row(l_t)$.
Then the monomial
$\pm \xi_{i_1,j_1}\cdots \xi_{i_r,j_r}
\xi_{k_1,l_1}\cdots\xi_{k_t,l_t}$
appears in $T_{i,h;x}(u)$ and the monomial
$\pm \xi_{k_{t+1},l_{t+1}}\cdots \xi_{k_s,l_s}$
appears in $T_{h,j;y}(u)$.
Moreover using the facts that
$\col(j_r) < \col(k_1)$ if $r > 0$,
$\col(l_t) \geq \col(k_{t+1})$, and the minimality of the choice of $t$,
we have that
$$
\pm(\xi_{i_1,j_1}\cdots \xi_{i_r,j_r})(\xi_{k_1,l_1}\cdots\xi_{k_s,l_s}),\:\:
\pm(\xi_{i_1,j_1}\cdots\xi_{i_r,j_r}\xi_{k_1,l_1}\cdots\xi_{k_t,l_t})(\xi_{k_{t+1},l_{t+1}}\cdots\xi_{k_s,l_s})
$$
are the only ways to split the monomial (\ref{rhs1})
so the first multiple
occurs in $T_{i,g;x}(u)$ and the second multiple appears in 
$T_{g,j;y}(u)$ for some $x < g \leq y$. The
two have opposite signs.

(ii) Similar.

(iii)
Using (ii), we can rewrite (iii) as
\begin{equation}\label{EGP1}
T_{i,j;x}(u)-T_{i,j;y}(u) = \sum_{k=x+1}^y T_{i,k;x}(u) T_{k,j;y}(u).
\end{equation}
The terms on the right hand side of (\ref{EGP1}) look like
\begin{equation}\label{EGP2}
\pm(\xi_{i_1,j_1}\cdots \xi_{i_r,j_r})(\xi_{k_1,l_1}\cdots\xi_{k_s,l_s}),
\end{equation}
for various $r,s \geq 1$,
where $\pm \xi_{i_1,j_1}\dots \xi_{i_r,j_r}$ appears in $T_{i,k;x}(u)$ and
$\pm \xi_{k_1,l_1}\dots  \xi_{k_s,l_s}$ appears in $T_{k,j;y}(u)$ for some
$x < k \leq y$.
Let $X$ be the sum of all such monomials
for which
$\col(j_r) < \col(k_1)$ 
and $\row(l_t) \notin \{x+1,\dots,y\}$ for
all $1 \leq t < s$.
Let $Y$ be the sum of all such monomials
for which $\col(j_r) \geq \col(k_1)$
and $\row(j_t) \notin \{x+1,\dots,y\}$ for all $1 \leq t < r$.
Let $Z$ be the sum of all the remaining monomials,
so the right hand side of (\ref{EGP1}) is equal to $X+Y+Z$.
We will show that $X+Y = T_{i,j;x}(u)-T_{i,j;y}(u)$ and
that $Z = 0$.

So first consider $X+Y$.
By definition,
$T_{i,j;x}(u)$ is a sum of monomials
of the form $\pm \xi_{p_1,q_1}\cdots \xi_{p_u,q_u}$.
Let $A$ denote the sum of all these monomials
with the property 
that $x < \row(q_t) \leq y$ for some $1 \leq t < u$.
Let $B$ be the sum of all the remaining monomials,
so $T_{i,j;x}(u) = A+B$.
Similarly, $T_{i,j;y}(u)$ is a sum of 
monomials $\pm \xi_{p_1,q_1}\cdots \xi_{p_u,q_u}$.
Let $C$ denote the sum of all the ones
with the property that $x < \row(q_t) \leq y$ for some $1 \leq t < u$.
Let $D$ denote the sum of all the rest, so
$T_{i,j;y}(u) = C+D$.
Now one checks that $X = A, Y = -C$ and $B = D$.
Hence $X+Y = (A+B)-(C+D) = T_{i,j;x}(u) - T_{i,j;y}(u)$
as claimed.

It remains to show that $Z = 0$. Recall that $Z$ is the sum of
all monomials of the form (\ref{EGP2})
such that
{\em either} 
$\col(j_r) < \col(k_1)$ and
$x < \row(l_t) \leq y$ for some $1 \leq t < s$,
{\em or}
$\col(j_r) \geq \col(k_1)$
and $x < \row(j_t) \leq y$ for some $1 \leq t < r$. 
In the former case, choose the 
minimal index $t \geq 1$ for which $x < \row(l_t) \leq y$. Then
$$
\pm(\xi_{i_1,j_1}\cdots \xi_{i_r,j_r})(\xi_{k_1,l_1}\cdots
\xi_{k_s,l_s}),\:\:
\pm(\xi_{i_1,j_1}\cdots \xi_{i_r,j_r}\xi_{k_1,l_1}\cdots
\xi_{k_t,l_t})(\xi_{k_{t+1},l_{t+1}}\cdots   \xi_{k_s,l_s})
$$
are the only two ways to split the monomial
(\ref{EGP2}) so that the left hand term appears in
$T_{i,g;x}(u)$ and the right hand term appears in $T_{g,j;y}(u)$
for some $x < g \leq y$.
Moreover, they appear in the expansion of $Z$ with opposite signs.
The latter case is
similar, but one uses the maximal index $t < r$ 
for which $x < \row(j_t)\leq y$.

(iv)
The left hand side of the formula in (iv)
is a sum of monomials of the form
\begin{equation}\label{E221003_1}
\pm (\xi_{i_1,j_1}\cdots \xi_{i_r,j_r})(\xi_{k_1,l_1}\cdots \xi_{k_s,l_s}),
\end{equation}
for various $r,s \geq 0$,
where
$\pm\xi_{i_1,j_1}\cdots \xi_{i_r,j_r}$ appears in $T_{i,k;x}(u)$ and
$\pm\xi_{k_1,l_1}\cdots \xi_{k_s,l_s}$ appears in $T_{k,j;y}(u)$ 
for some $x < k \leq y$. 
Note that $r=0$ is allowed only if $i=k$, and $s=0$ is allowed
only if $k=j$. So if $i=j$ we get a contribution $-1$ from the 
$r=s=0$ term, while if $i \neq j$ then there is no $r=s=0$
term.
Now the formula in (iv) will follow if we can show that
whenever $(r,s) \neq (0,0)$, the monomial 
(\ref{E221003_1}) appears exactly twice 
on the left hand side
with opposite signs. 
We consider
two cases. 

First, suppose that $s > 0$ and that
$\col(j_r) < \col(k_1)$ if $r > 0$. Let
 $t \geq 1$ be the minimal index such that $x < \row(l_t) \leq y$;
such a $t$ exists as $\row(l_s)=j$ and $x < j \leq y$. 
Let $h := \row(l_t)$.
Then
$\pm \xi_{i_1,j_1}\dots \xi_{i_r,j_r}\xi_{k_1,l_1}\dots \xi_{k_t,l_t}$
appears in $T_{i,h;x}(u)$, and
$\pm \xi_{k_{t+1},l_{t+1}}\dots \xi_{k_s,l_s}$
appears in $T_{h,j;y}(u)$.
Moreover, using the facts that 
$\col(l_t)  \geq \col(k_{t+1})$  if $t<s$,
$\col(j_r)<\col(k_1)$ if $r > 0$,
and the minimality of the choice of $t$, we see that
$$
\pm (\xi_{i_1,j_1}\dots \xi_{i_r,j_r})(\xi_{k_1,l_1}\dots \xi_{k_s,l_s}),\:\:
\pm (\xi_{i_1,j_1}\dots \xi_{i_r,j_r}\xi_{k_1,l_1}\dots \xi_{k_t,l_t})
(\xi_{k_{t+1},l_{t+1}}\dots \xi_{k_s,l_s})
$$
are the only ways to split the monomial (\ref{E221003_1}) so that the
left term occurs in $T_{i,g;x}(u)$ and the right term occurs in
$T_{g,j;y}(u)$ for some $x < g \leq y$. 
The two have opposite signs.

For the second case, suppose either that
$s=0$, or that $r,s>0$ and $\col(j_r)\geq \col(k_1)$. Let
 $t \leq r$ be the maximal index such that $x < \row(i_t) \leq y$
and argue in a similar fashion.
\iffalse
such a $t$ exists as $\row(i_1)=i$ and $x < i \leq y$. Let $h := \row(i_t)$.
Then
$\pm \xi_{i_1,j_1}\dots \xi_{i_{t-1},j_{t-1}}$
appears in $T_{i,h;x}(u)$, and
$\pm \xi_{i_t,j_t}\dots \xi_{i_r,j_r}\xi_{k_1,l_1}\dots  \xi_{k_s,l_s}$
appears in $T_{h,j;y}(u)$.
Moreover, using the facts that $\col(j_{t-1})< \col(i_t)$ if $t > 1$,
$\col(j_r)\geq \col(k_1)$ if $s > 0$,
and the maximality of the choice of $t$, we see that
$$
\pm (\xi_{i_1,j_1}\dots \xi_{i_r,j_r})(\xi_{k_1,l_1}\dots \xi_{k_s,l_s}),
\:\:
\pm (\xi_{i_1,j_1}\dots \xi_{i_{t-1},j_{t-1}})
(\xi_{i_t,j_t}\dots \xi_{i_r,j_r}\xi_{k_1,l_1}\dots  \xi_{k_s,l_s})
$$
are the only ways to split the monomial (\ref{E221003_1}) so that the
left term occurs in $T_{i,g;x}(u)$ and the right term occurs in
$T_{g,j;y}(u)$ for some $x < g \leq y$. 
It just remains to check that the two have opposite signs.
\fi\end{proof}

To explain the significance of this lemma, 
let $T(u) := (T_{i,j;0}(u))_{1 \leq i,j \leq n}$,
an $n \times n$ matrix with entries in $U(\mathfrak{p})[[u^{-1}]]$.
Also let $\nu = (\nu_1,\dots,\nu_m)$ be a fixed shape.
Consider the Gauss factorization
$T(u) = F(u)D(u)E(u)$ where $D(u)$ is a block diagonal matrix,
$E(u)$ is a block upper unitriangular matrix, and $F(u)$ is a block
lower unitriangular matrix, all block matrices being of shape $\nu$.
The diagonal blocks of $D(u)$ define
matrices $D_1(u),\dots,D_m(u)$, 
the upper diagonal blocks
of $E(u)$ define matrices $E_1(u),\dots,E_{m-1}(u)$,
and the lower diagonal matrices of $F(u)$ define matrices
$F_1(u),\dots,F_{m-1}(u)$.
Also let $\widetilde{D}_a(u) := -D_a(u)^{-1}$.
Thus $D_a(u) = (D_{a;i,j}(u))_{1 \leq i,j \leq \nu_a}$
and $\widetilde{D}_a(u) = (\widetilde{D}_{a;i,j}(u))_{1 \leq i,j \leq \nu_a}$
are $\nu_a \times \nu_a$-matrices,
$E_a(u) = (E_{a;i,j}(u))_{1 \leq i \leq \nu_a, 1 \leq j \leq \nu_{a+1}}$
is a $\nu_a \times \nu_{a+1}$-matrix, and
$F_a(u) = (F_{a;i,j}(u))_{1 \leq i \leq \nu_{a+1}, 1 \leq j \leq \nu_{a}}$
is a $\nu_{a+1} \times \nu_{a}$-matrix.
Write
\begin{align*}
D_{a;i,j}(u) &= \sum_{r \geq 0} D_{a;i,j}^{(r)} u^{-r},\quad
&\widetilde{D}_{a;i,j}(u) &= \sum_{r \geq 0} 
\widetilde{D}_{a;i,j}^{(r)} u^{-r},\\
E_{a;i,j}(u) &= \sum_{r > 0} E_{a;i,j}^{(r)} u^{-r},\quad
&F_{a;i,j}(u) &= \sum_{r > 0} F_{a;i,j}^{(r)} u^{-r},
\end{align*}
thus defining elements $D_{a;i,j}^{(r)}, E_{a;i,j}^{(r)}$
and $F_{a;i,j}^{(r)}$ of $U(\mathfrak{p})$, 
all dependent of course on the fixed choice of $\nu$.
All this
parallels the definition
of the elements of $Y_n$ with the same names 
in the paragraph following Lemma~\ref{quasi}.

\begin{Theorem}\label{transtheorem}
With $\nu = (\nu_1,\dots,\nu_m)$ fixed as above and all admissible
$a,i,j$, we have that
\begin{align*}
D_{a;i,j}(u) &= T_{\nu_1+\cdots+\nu_{a-1}+i,\nu_1+\cdots+\nu_{a-1}+j;
\nu_1+\cdots+\nu_{a-1}}(u),\\
\widetilde{D}_{a;i,j}(u) &=
T_{\nu_1+\cdots+\nu_{a-1}+i,\nu_1+\cdots+\nu_{a-1}+j;
\nu_1+\cdots+\nu_a}(u),\\
E_{a;i,j}(u) &= T_{\nu_1+\cdots+\nu_{a-1}+i,\nu_1+\cdots+\nu_{a}+j;
\nu_1+\cdots+\nu_{a}}(u),\\
F_{a;i,j}(u) &= T_{\nu_1+\cdots+\nu_{a}+i,\nu_1+\cdots+\nu_{a-1}+j;
\nu_1+\cdots+\nu_{a}}(u).
\end{align*}
\end{Theorem}

\begin{proof}
Note it suffices to prove the formulae for $D, E$ and $F$, since the
one for $\widetilde{D}$ follows from the one for $D$ 
by Lemma~\ref{traninv}(iv) taking $x =\nu_1+\cdots+\nu_{a-1}$
and $y = \nu_1+\cdots+\nu_a$.
Now we proceed by induction on $m$, 
the base case $m=1$ being trivial.
For the induction step, suppose the theorem has been proved for
the shape $\nu = (\nu_1,\dots,\nu_m)$.
So, in terms of matrices, we have that
\begin{align*}
{^\nu}D_{a}(u) &= \left(T_{\nu_1+\cdots+\nu_{a-1}+i,\nu_1+\cdots+\nu_{a-1}+j;
\nu_1+\cdots+\nu_{a-1}}(u)\right)_{1 \leq i,j \leq \nu_a},\\
{^\nu}E_{a}(u) &= \left(T_{\nu_1+\cdots+\nu_{a-1}+i,\nu_1+\cdots+\nu_{a}+j;
\nu_1+\cdots+\nu_{a}}(u)\right)_{1 \leq i \leq \nu_a, 1 \leq j \leq \nu_{a+1}},\\
{^\nu}F_{a}(u) &= \left(T_{\nu_1+\cdots+\nu_{a}+i,\nu_1+\cdots+\nu_{a-1}+j;
\nu_1+\cdots+\nu_{a}}(u)\right)_{1 \leq i \leq \nu_{a+1}, 1 \leq j \leq \nu_{a}},
\end{align*}
where we have added a superscript $\nu$ for clarity.
Write 
$\nu_b = \alpha+\beta$ for some $1 \leq b \leq m$ and 
$\alpha,\beta \geq 1$, and let $\mu = (\nu_1,\dots,\nu_{b-1},\alpha,\beta,
\nu_{b+1},\dots,\nu_m)$.
Define matrices $A,B,C$ and $D$ by
\begin{align*}
A &= \left( T_{\nu_1+\cdots+\nu_{b-1}+i,\nu_1+\cdots+\nu_{b-1}+j;\nu_1+\cdots+\nu_{b-1}}(u)\right)_{1 \leq i,j \leq \alpha},\\
B &= \left( T_{\nu_1+\cdots+\nu_{b-1}+i,\nu_1+\cdots+\nu_{b-1}+\alpha+j;\nu_1+\cdots+\nu_{b-1}+\alpha}(u)\right)_{1 \leq i \leq \alpha, 1 \leq j \leq \beta},\\
C &= \left( T_{\nu_1+\cdots+\nu_{b-1}+\alpha+i,\nu_1+\cdots+\nu_{b-1}+j;\nu_1+\cdots+\nu_{b-1}+\alpha}(u)\right)_{1 \leq i \leq \beta, 1 \leq 
j \leq \alpha},\\
D &= \left( T_{\nu_1+\cdots+\nu_{b-1}+\alpha+i,\nu_1+\cdots+\nu_{b-1}+\alpha+j;\nu_1+\cdots+\nu_{b-1}+\alpha}(u)\right)_{1 \leq i,j \leq \beta}.
\end{align*}
Then Lemma~\ref{traninv}(i)--(iii) with $x = \nu_1+\cdots+\nu_{b-1}$
and $y = \nu_1+\cdots+\nu_{b-1}+\alpha$ tell us that
$$
{^\nu}D_b(u) = 
\left(
\begin{array}{cc}
A&AB\\
CA&D+CAB
\end{array}\right)
=
\left(\begin{array}{ll}I_\alpha&0\\ C&I_\beta\end{array}\right)
\left(\begin{array}{ll}A&0\\0&D\end{array}\right)
\left(\begin{array}{ll}I_\alpha&B\\0&I_\beta\end{array}\right).
$$
Lemma~\ref{quasi} explains how to 
read off the matrices
${^\mu}D_a(u), {^\mu}E_a(u)$ and ${^\mu}F_a(u)$
from this factorization
to get that 
\begin{align*}
{^\mu}D_{a}(u) &= \left(T_{\mu_1+\cdots+\mu_{a-1}+i,\mu_1+\cdots+\mu_{a-1}+j,
\mu_1+\cdots+\mu_{a-1}}(u)\right)_{1 \leq i,j \leq \mu_a},\\
{^\mu}E_{a}(u) &= \left(T_{\mu_1+\cdots+\mu_{a-1}+i,\mu_1+\cdots+\mu_{a}+j,
\mu_1+\cdots+\mu_{a}}(u)\right)_{1 \leq i \leq \mu_a, 1 \leq j \leq \mu_{a+1}},\\
{^\mu}F_{a}(u) &= \left(T_{\mu_1+\cdots+\mu_{a}+i,\mu_1+\cdots+\mu_{a-1}+j,
\mu_1+\cdots+\mu_{a}}(u)\right)_{1 \leq i \leq \mu_{a+1}, 1 \leq j \leq \mu_{a}}.
\end{align*}
This completes the induction step.
\end{proof}

In the extreme case that $\nu = (1^n)$, we
write simply $D_i^{(r)},
\widetilde{D}_i^{(r)},E_i^{(r)}$ and
$F_i^{(r)}$ for the elements $D_{i;1,1}^{(r)}$,
$\widetilde{D}_{i;1,1}^{(r)}$,
$E_{i;1,1}^{(r)}$ and $F_{i;1,1}^{(r)}$
of $U(\mathfrak{p})$, respectively.

\begin{Corollary}\label{transcor}
$D_i^{(r)} = T_{i,i;i-1}^{(r)}$,
$\widetilde{D}_i^{(r)} = T_{i,i;i}^{(r)}$,
$E_i^{(r)} = T_{i,i+1;i}^{(r)}$ and
$F_i^{(r)} = T_{i+1,i;i}^{(r)}$.
\end{Corollary}

\section{Main theorem}\label{smain}

Let $\pi$ be a pyramid of level $l$ with 
column heights $(q_1,\dots,q_l)$. 
Pick an integer $n \geq \max(q_1,\dots,q_l)$ and read off from the pyramid $\pi$ a shift matrix
$\sigma = (s_{i,j})_{1 \leq i,j \leq n}$ according to (\ref{pyr3}).
Define the
finite $W$-algebra $W(\pi) = U(\mathfrak{p})^{\mathfrak{m}}$
viewed as a filtered algebra via the Kazhdan filtration as
in $\S$\ref{sslice}.
Define the shifted Yangian $Y_{n,l}(\sigma)$ of level $l$
viewed as a filtered algebra via the canonical filtration
as in $\S\S$\ref{syangian}, \ref{scan} and \ref{struncation}.
Suppose also that $\nu = (\nu_1,\dots,\nu_m)$ is an admissible
shape for $\sigma$, and recall the notation $s_{a,b}(\nu)$ and $p_a(\nu)$
from (\ref{rels}) and (\ref{padef}).
We have the elements $D_{a;i,j}^{(r)}$,
$\widetilde{D}_{a;i,j}^{(r)}$,
$E_{a;i,j}^{(r)}$ and $F_{a;i,j}^{(r)}$
of $U(\mathfrak{p})$  defined as in $\S$\ref{sinvariants} 
relative to this fixed shape;
see Theorem~\ref{transtheorem} for the explicit formulae.
We also have the parabolic generators
$D_{a;i,j}^{(r)}$,
$\widetilde{D}_{a;i,j}^{(r)}$,
$E_{a;i,j}^{(r)}$ and $F_{a;i,j}^{(r)}$ 
of $Y_{n,l}(\sigma)$ 
as in $\S$\ref{sparabolic}.
The main result of the article is as follows.

\begin{Theorem}\label{main}
There is a unique isomorphism
$Y_{n,l}(\sigma) \stackrel{\sim}{\rightarrow} W(\pi)$ 
of filtered algebras
such that
for any admissible shape $\nu = (\nu_1,\dots,\nu_m)$
the generators
\begin{align*}
&\{D_{a;i,j}^{(r)}\}_{1 \leq a \leq m,
1 \leq i,j \leq \nu_a, r > 0},\\
&\{E_{a;i,j}^{(r)}\}_{1 \leq a < m, 1 \leq i \leq \nu_a,
1 \leq j \leq \nu_{a+1}, r > s_{a,b}(\nu)},\\
&\{F_{a;i,j}^{(r)}\}_{1 \leq a < m, 1 \leq i \leq \nu_{a+1},
1 \leq j \leq \nu_{a}, r > s_{b,a}(\nu)}
\end{align*}
 of $Y_{n,l}(\sigma)$
map to the elements of $U(\mathfrak{p})$ with the same names.
In particular, these elements of $U(\mathfrak{p})$
are $\mathfrak{m}$-invariants and they generate $W(\pi)$.
\end{Theorem}

The proof will take up most of the rest of the section.
First however let us record a couple of corollaries of the theorem.

\begin{Corollary}\label{onto}
The map $\theta:\gr W(\pi) \rightarrow S(\mathfrak{c}_{\mathfrak{g}}(e))$
from (\ref{thed})
is an isomorphism.
\end{Corollary}

\begin{proof}
Since we already know by 
Corollary~\ref{phewy} that the map
$\theta$ is an injective, graded homomorphism, it suffices
given Theorem~\ref{main} to show that 
$\gr Y_{n,l}(\sigma)$ and $S(\mathfrak{c}_{\mathfrak{g}}(e))$
have the same dimension in each degree. 
This follows from Theorem~\ref{miuramain}
and Lemma~\ref{centbase} (the element
$c_{i,j}^{(r)}$ there is of degree $r$ with respect to the
Kazhdan grading).
\end{proof}

For the next corollary, let $\dot\pi$ be another pyramid with the same
row lengths as $\pi$, i.e. the nilpotent matrix $\dot e$
defined from $\dot\pi$ is conjugate to the nilpotent matrix $e$
defined from $\pi$.
Let $\dot\sigma = (\dot s_{i,j})_{1 \leq i,j \leq n}$ be a shift matrix
corresponding to the pyramid $\dot\pi$.
By Theorem~\ref{main}, the formulae in $\S$\ref{sinvariants} define
generators
$\dot D_{a;i,j}^{(r)}$,
$\dot E_{a;i,j}^{(r)}$ and $\dot F_{a;i,j}^{(r)}$
of $W(\dot\pi)$ for each admissible shape $\nu$.

\begin{Corollary}\label{oneiso}
There is an algebra isomorphism $\iota:W(\pi)\rightarrow W(\dot\pi)$
defined on parabolic generators with respect to an admissible shape
$\nu$ by (\ref{iotadefpar}).
\end{Corollary}

\begin{proof}
This follows from Theorem~\ref{main} and (\ref{iotadef2}).
\end{proof}

So now we must prove Theorem~\ref{main}.
The first reduction is to observe that it suffices
to prove it in the special
case that $\nu$ is the minimal admissible shape for $\sigma$.
It then follows for all other admissible shapes by Lemma~\ref{quasi}
and induction on the length of the shape.
So we assume from now on that $\nu=(\nu_1,\dots,\nu_m)$ 
is the minimal admissible shape for $\sigma$, 
and let $t := \nu_m = \min(q_1,q_l)$.
We proceed to prove Theorem~\ref{main}
by induction on the level $l$. 

Consider first the base case $l=1$, i.e. 
$\pi$ consists of a single column of height $t$.
Since the nilpotent element $e$ from (\ref{edef})
is zero in this case, we have by definition that
$W(\pi) = U(\mathfrak{gl}_t)$. Moreover,
according to Theorem~\ref{transtheorem}, we have that
$D_{m;i,j}^{(1)} = \tilde e_{i,j}
= e_{i,j} + \delta_{i,j}(n-t)$ for $1 \leq i,j \leq t$,
and all other elements $D_{a;i,j}^{(r)}, E_{a;i,j}^{(r)}$ and
$F_{a;i,j}^{(r)}$ of $W(\pi)$ are equal to zero. 
On the other hand, it follows from the relations and Corollary~\ref{newpbw} 
that there is an 
isomorphism $Y_{n,1}(\sigma) \stackrel{\sim}{\rightarrow} 
U(\mathfrak{gl}_t)$ such that
$D_{m,i,j}^{(1)} \mapsto e_{i,j}$, with all the other parabolic generators
$D_{a;i,j}^{(r)},
E_{a;i,j}^{(r)}$ and $F_{a;i,j}^{(r)}$ 
of $Y_{n,1}(\sigma)$ mapping to zero.
Composing with the automorphism $e_{i,j} \mapsto \tilde e_{i,j}$
of $U(\mathfrak{gl}_t)$ gives the 
required isomorphism $Y_{n,1}(\sigma) \stackrel{\sim}{\rightarrow}
W(\pi)$. 

Assume from now on that $l > 1$ 
and that the theorem has been proved for all
smaller levels. 
The verification of the induction step
splits into two cases corresponding to the two possible 
baby comultiplications $\Delta_{\rt}$ and $\Delta_{\lt}$ that were defined
in (\ref{deltar})--(\ref{deltal}):

\noindent
{\em Case $\Delta_{\rt}$:} 
either $t=n$ or $s_{n-t,n-t+1} \neq 0$;

\noindent
{\em Case $\Delta_{\lt}$:} either $t=n$ or $s_{n-t+1,n-t} \neq 0$.

\noindent
The argument in the two cases is quite similar.  We will  explain it in 
detail in
case $\Delta_{\rt}$, then briefly indicate the changes in case $\Delta_{\lt}$.

So assume that 
either $t = n$ or
$s_{n-t,n-t+1} = s_{m-1,m}(\nu) \neq 0$; in particular, $q_1 \geq q_l$.
For notational convenience, assume that the numbering of the
bricks of the pyramid $\pi$ is the standard numbering
down columns from left to right.
Let $\dot\pi$ be the pyramid
obtained from $\pi$ by removing the rightmost column,
i.e. the bricks numbered $(N-t+1),(N-t+2),\dots,N$, from the pyramid $\pi$.
Let $\dot\sigma = 
(\dot s_{i,j})_{1 \leq i,j \leq n}$ be 
the shift matrix corresponding to the pyramid $\dot\pi$ defined by (\ref{sr}).
Define $\dot{\mathfrak{p}}, \dot{\mathfrak{m}}$ and $\dot e$
in $\dot{\mathfrak{g}} = \mathfrak{gl}_{N-t}$
according to (\ref{Ac1})--(\ref{edef})
and let $\dot\chi:\dot{\mathfrak{m}}\rightarrow \C$
be the character $x \mapsto (x,\dot e)$.
Let $\dot D_{a;i,j}^{(r)}, \dot{{\widetilde D}}_{a;i,j}^{(r)}$,
$\dot E_{a;i,j}^{(r)}$ and $\dot F_{a;i,j}^{(r)}$ denote the
elements of $U(\dot{\mathfrak{p}})$ as defined in $\S$\ref{sinvariants}
relative to the shape $\nu$. 
By the induction hypothesis, Theorem~\ref{main} holds for $\dot\pi$,
so we know already that the following elements of
$U(\dot{\mathfrak{p}})$ are
invariant under the twisted action of $\dot{\mathfrak{m}}$, 
i.e. they belong to
finite $W$-algebra
$W(\dot\pi) = U(\dot{\mathfrak{p}})^{\dot{\mathfrak{m}}}$:
\begin{itemize}
\item[(i)] $\dot D_{a;i,j}^{(r)}$ and $\dot{{\widetilde D}}_{a;i,j}^{(r)}$
for  $1 \leq a \leq m$, $1 \leq i,j \leq \nu_a$ and $r > 0$;
\item[(ii)] $\dot E_{a;i,j}^{(r)}$ for 
$1 \leq a < m$, $1 \leq i \leq \nu_a$, $1 \leq j \leq \nu_{a+1}$
and $r > s_{a,a+1}(\nu) - \delta_{a,m-1}$;
\item[(iii)] $\dot F_{a;i,j}^{(r)}$ 
 for 
$1 \leq a< m$, $1 \leq i \leq \nu_{a+1}$, $1 \leq j \leq \nu_{a}$
and $r > s_{a+1,a}(\nu)$.
\end{itemize}
Recalling (\ref{etildedef}), 
we must work now with the non-standard embedding
of $U(\dot{\mathfrak{g}})$ into $U(\mathfrak{g})$ under which the
generators $\tilde e_{i,j}$ of $U(\dot{\mathfrak{g}})$ defined from 
the pyramid $\dot\pi$ map to the generators
$\tilde e_{i,j}$ of $U(\mathfrak{g})$ defined from the pyramid $\pi$,
for $1 \leq i,j \leq N-t$.
This also embeds $U(\dot{\mathfrak{p}})$ into $U(\mathfrak{p})$ and
$\dot{\mathfrak{m}}$ into $\mathfrak{m}$.
Moreover, recalling (\ref{chidef}),
the character $\dot\chi$ of
$\dot{\mathfrak{m}}$ is the restriction of the character
$\chi$ of $\mathfrak{m}$, hence the twisted action of
$\dot{\mathfrak{m}}$ on $U(\dot{\mathfrak{p}})$ is the restriction of the
twisted action of $\mathfrak{m}$ on $U(\mathfrak{p})$.
The following crucial lemma gives us an inductive 
description of the elements 
$D_{a;i,j}^{(r)}$,
$E_{a;i,j}^{(r)}$ and $F_{a;i,j}^{(r)}$ of
$U(\mathfrak{p})$; 
this should be compared with Theorem~\ref{baby}(i). 

\begin{Lemma}\label{superbaby}
The following equations hold for $r > 0$, all admissible $a,i,j$ and 
any fixed $1 \leq h \leq t$:
\begin{align}\label{sb1}
D_{a;i,j}^{(r)} & = \dot D_{a;i,j}^{(r)} + \delta_{a,m}
\left(\sum_{k=1}^t
\dot D_{a;i,k}^{(r-1)} \tilde e_{N-t+k,N-t+j}
+
[\dot D_{a;i,h}^{(r-1)}, \tilde e_{N-2t+h,N-t+j}]\right),\\
E_{a;i,j}^{(r)} & = \dot E_{a;i,j}^{(r)} + \delta_{a,m-1}\left(
\sum_{k=1}^t \dot E_{a;i,k}^{(r-1)} \tilde e_{N-t+k,N-t+j}
+ [\dot E_{a;i,h}^{(r-1)}, \tilde e_{N-2t+h,N-t+j}]\right),\label{sb2}\\
F_{a;i,j}^{(r)} & = \dot F_{a;i,j}^{(r)},
\end{align}
where for (\ref{sb2}) we are assuming that $r > 1$ if $a=m-1$.
\end{Lemma}

\begin{proof}
This follows using Theorem~\ref{transtheorem} and the explicit form of the
elements $T_{i,j;x}^{(r)}$ from (\ref{thedef}).
\end{proof}

In the next few lemmas
we will use these inductive descriptions to show that the elements
$D_{a;i,j}^{(r)}, E_{a;i,j}^{(r)}$ and $F_{a;i,j}^{(r)}$ 
of $U(\mathfrak{p})$
are $\mathfrak{m}$-invariants for the appropriate $r$.

\begin{Lemma}\label{hor1}
The following elements of $U(\mathfrak{p})$ are
invariant under the twisted action of $\mathfrak{m}$:
\begin{itemize}
\item[(i)] $D_{a;i,j}^{(r)}$ and ${{\widetilde D}}_{a;i,j}^{(r)}$
for  $1 \leq a \leq m-1$, $1 \leq i,j \leq \nu_a$ and $r > 0$;
\item[(ii)] $E_{a;i,j}^{(r)}$ for 
$1 \leq a < m-1$, $1 \leq i \leq \nu_a$, $1 \leq j \leq \nu_{a+1}$
and $r>s_{a,a+1}(\nu)$;
\item[(iii)] $F_{a;i,j}^{(r)}$ 
 for 
$a=1 \leq a < m$, $1 \leq i \leq \nu_{a+1}$, $1 \leq j \leq \nu_{a}$
and $r > s_{a+1,a}(\nu)$.
\end{itemize}
\end{Lemma}

\begin{proof}
By Lemma~\ref{superbaby} and the definitions of
${\widetilde D}_{a;i,j}^{(r)}$ and $\dot{{\widetilde D}}_{a;i,j}^{(r)}$, 
all these elements of $U(\mathfrak{p})$
coincide with the corresponding elements of $U(\dot{\mathfrak{p}})$.
Hence by the induction hypothesis we already know they are invariant under the
twisted action of $\dot{\mathfrak{m}}$.
It remains to show that the elements are invariant under the twisted action
of all $\tilde e_{f,g}$ with $1 \leq g \leq N-t < f \leq N$.
By Theorem \ref{transtheorem} and the explicit form of (\ref{thedef}),
all the elements we are considering are linear combinations
of monomials of the form
$\tilde e_{i_1,j_1} \cdots \tilde e_{i_r,j_r} \in U(\dot{\mathfrak{p}})$
with $1 \leq i_s \leq N-t$ and $1 \leq j_s \leq N-2t$ for
all $s=1,\dots,r$.
Using just the fact that $\chi(e_{f,g}) = 0$
for all $1 \leq g \leq N-2t$ and $N-t < f \leq N$, it is easy
to see that all such monomials are invariant under the twisted
action of all $\tilde e_{f,g}$ with $1 \leq g \leq N-t<f\leq N$. 
\end{proof}

\begin{Lemma}\label{hor2}
The following elements of $U(\mathfrak{p})$ 
are invariant under the twisted action of
$\dot{\mathfrak{m}}$:
\begin{itemize}
\item[(i)] $D_{m;i,j}^{(r)}$ for $1 \leq i,j \leq \nu_m$
and $r > 0$;
\item[(ii)] (assuming $t < n$) 
$E_{m-1;i,j}^{(r)}$ for $1 \leq i \leq \nu_{m-1}$,
$1 \leq j \leq \nu_m$ and $r > s_{m-1,m}(\nu)$.
\end{itemize}
\end{Lemma}

\begin{proof}
(i) 
Take $x \in \dot{\mathfrak{m}}$.
According to Lemma~\ref{superbaby}, we have that
$$
D_{m;i,j}^{(r)} = \dot D_{m;i,j}^{(r)}
+ \sum_{k=1}^t \dot D_{m;i,k}^{(r-1)} \tilde e_{N-t+k,N-t+j}
+ [\dot D_{m;i,h}^{(r-1)}, \tilde e_{N-2t+h,N-t+j}].
$$
Noting that
$[x,\tilde e_{N-t+k,N-t+j}] =
[x,\tilde e_{N-2t+h,N-t+j}] = 0$ and using the induction
hypothesis, one deduces easily from this equation that
$\pr_\chi([x,  D_{m;i,j}^{(r)}]) = 0$ as required.

(ii) Similar.
\end{proof}

\begin{Lemma}\label{expinvs}
The following elements of $U(\mathfrak{p})$ 
are invariant under the twisted action
of $\tilde e_{f,g}$ 
for all $1 \leq g \leq N-t < f \leq N$:
\begin{itemize}
\item[(i)] $D_{m;i,j}^{(1)}$
for $1 \leq i,j \leq \nu_m$;
\item[(ii)] (assuming $t < n$ and $s_{m-1,m}(\nu)=1$)
$D_{m;i,j}^{(2)}$
for $1 \leq i,j \leq \nu_m$;
\item[(iii)] (assuming $t < n$ and $s_{m-1,m}(\nu)=1$)
$E_{m-1;i,j}^{(2)}$ for $1 \leq i \leq \nu_{m-1}$,
$1 \leq j \leq \nu_m$.
\end{itemize}
\end{Lemma}

\begin{proof}
Part (i) is easily checked directly
using (\ref{erel}), (\ref{chidef}) and the explicit
formula for $D_{m;i,j}^{(1)}$ given by Theorem~\ref{transtheorem}
and Example~\ref{eggy}.
The proofs of (ii) and (iii) are similar, though the
calculations are not so easy.
\end{proof}

\begin{Lemma}\label{hor3}
Suppose that 
$t < n$ and $s_{m-1,m}(\nu) = 1$.
Then, the following equations hold in $U(\mathfrak{p})$
for $r > 1$, all admissible $i,j$ and any fixed
$1 \leq g \leq \nu_{m-1}$:
\begin{align}\label{me1}
D_{m;i,j}^{(r+1)} &= 
[F_{m-1;i,g}^{(2)}, E_{m-1;g,j}^{(r)}]
+ \sum_{s=0}^r {\widetilde D}_{m-1;g,g}^{(r+1-s)} D_{m;i,j}^{(s)},\\
E_{m-1;i,j}^{(r+1)} &= [D_{m-1;i,g}^{(2)}, E_{m-1;g,j}^{(r)}]
- \sum_{f=1}^{\nu_{m-1}} D_{m-1;i,f}^{(1)} E_{m-1;f,j}^{(r)}.
\label{me2}
\end{align}
\end{Lemma}

\begin{proof}
We prove (\ref{me2}), the first equation being a similar trick.
By the induction hypothesis and the relation (\ref{pr4}), we know  
for all $r > 0$ that
$$
[\dot D_{m-1;i,g}^{(2)}, \dot E_{m-1;g,j}^{(r)}]
= \dot E_{m-1;i,j}^{(r+1)} + \sum_{f=1}^{\nu_{m-1}}
\dot D_{m-1;i,f}^{(1)} \dot E_{m-1;f,j}^{(r)}.
$$
By Lemma~\ref{superbaby}, we have for $r > 1$ that
$$
E_{m-1;g,j}^{(r)}
=
\dot E_{m-1;g,j}^{(r)}
+ \sum_{k=1}^t \dot E_{m-1;g,k}^{(r-1)} \tilde e_{N-t+k,N-t+j}
+ [\dot E_{m-1;g,h}^{(r-1)}, \tilde e_{N-2t+h,N-t+j}].
$$
Obviously,
$[\dot D_{m-1;i,g}^{(2)}, \tilde e_{N-t+k,N-t+j}]
=0$.
Moreover, by Theorem~\ref{transtheorem} and the form of
(\ref{thedef}), no monomial
in the expansion of $\dot D_{m-1;i,g}^{(2)}$ involves any
matrix unit of the form $\tilde e_{?, N-2t+h}$, hence 
$[\dot D_{m-1;i,g}^{(2)}, \tilde e_{N-2t+h,N-t+j}]
= 0$ too. 
Now we can commute with $\dot D_{m-1;i,g}^{(2)} = D_{m-1;i,g}^{(2)}$
to deduce that
\begin{align*}
[D_{m-1;i,g}^{(2)}, E_{m-1;g,j}^{(r)}] &=
\dot E_{m-1;i,j}^{(r+1)} + \sum_{f=1}^{\nu_{m-1}} \dot D_{m-1;i,f}^{(1)} \dot E_{m-1;f,j}^{(r)}\\
 &+ \sum_{k=1}^t \left\{
\dot E_{m-1;i,k}^{(r)} + 
\sum_{f=1}^{\nu_{m-1}} \dot D_{m-1;i,f}^{(1)}
\dot E_{m-1;f,k}^{(r-1)}
\right\}\tilde e_{N-t+k,N-t+j} \\
&+ 
\left[\dot E_{m-1;i,h}^{(r)}
+ \sum_{f=1}^{\nu_{m-1}} \dot D_{m-1;i,f}^{(1)} \dot E_{m-1;f,h}^{(r-1)},
\tilde e_{N-2t+h,N-t+j}\right].
\end{align*}
Finally rewrite the right hand side using Lemma \ref{superbaby} again to see
that it equals $E_{m-1;i,j}^{(r+1)} + \sum_{f=1}^{\nu_{m-1}}
D_{m-1;i,f}^{(1)}E_{m-1;f,j}^{(r)}$.
\end{proof}

\begin{Lemma}\label{hor4}
Suppose that $t = n$ and $l > 1$ or that $t < n$ and $s_{m-1,m}(\nu) > 1$.
Then the following elements of $U(\mathfrak{p})$ are
invariant under the twisted action of $\tilde e_{N-t+f,N-2t+g}$
for all $1 \leq f,g \leq t$.
\begin{itemize}
\item[(i)] $D_{m;i,j}^{(r)}$ for $1 \leq i,j \leq\nu_m$ and $r > 1$;
\item[(ii)] (assuming $t < n$)
$E_{m-1;i,j}^{(r)}$ for $1 \leq i \leq \nu_{m-1}$,
$1 \leq j \leq \nu_m$ and $r > s_{m-1,m}(\nu)$.
\end{itemize}
\end{Lemma}

\begin{proof}
(i)
Let $\ddot{\pi}$ denote the pyramid obtained by removing
the rightmost column from the pyramid $\dot\pi$, i.e. the bricks
numbered $(N-2t+1),(N-2t+2),\dots,(N-t)$.
Define $\ddot{\mathfrak{p}}, \ddot{\mathfrak{m}}$ and $\ddot e$
in $\ddot{\mathfrak{g}} = \mathfrak{gl}_{N-2t}$
according to (\ref{Ac1})--(\ref{edef}), and
embed $U(\ddot{\mathfrak{g}})$ 
into $U(\dot{\mathfrak{g}})$ in exactly the same (non-standard) way as we 
embedded $U(\dot{\mathfrak{g}})$ into $U(\mathfrak{g})$.
We will make use of the elements
$\ddot{D}_{a;i,j}^{(r)}$
of $W(\ddot{\pi})$.
By Lemma~\ref{superbaby} applied to $\dot\pi$,
the following equation holds for $r > 0$,
$1 \leq i,j \leq \nu_m$ and
any fixed $1 \leq h \leq t$:
\begin{align*}
\dot D_{m;i,j}^{(r)} & = \ddot{D}_{m;i,j}^{(r)} + 
\sum_{c=1}^t
\ddot{D}_{m;i,c}^{(r-1)} \tilde e_{N-2t+c,N-2t+j}
+
[\ddot{D}_{m;i,h}^{(r-1)}, \tilde e_{N-3t+h,N-2t+j}].
\end{align*}
Substituting this into (\ref{sb1}) and simplifying
a little using (\ref{erel}) we deduce for $r > 1$ that
$D_{m;i,j}^{(r)}= A+B+C+D+E+F+G+H$
where
\begin{align*}
A &= \ddot{D}_{m;i,j}^{(r)},
&E &=\sum_{k,c=1}^t \ddot{D}_{m;i,c}^{(r-2)}
\tilde e_{N-2t+c,N-2t+k}\tilde e_{N-t+k,N-t+j},\\
B &= \sum_{c=1}^t \ddot{D}_{m;i,c}^{(r-1)}\tilde e_{N-2t+c,N-2t+j},
&F &=\sum_{k=1}^t
[\ddot{D}_{m;i,h}^{(r-2)}, \tilde e_{N-3t+h,N-2t+k}]\tilde e_{N-t+k,N-t+j},\\
C &=[\ddot{D}_{m;i,h}^{(r-1)}, \tilde e_{N-3t+h,N-2t+j}],
&G &=\sum_{c=1}^t\ddot{D}_{m;i,c}^{(r-2)} \tilde e_{N-2t+c,N-t+j},\\
D &=\sum_{k=1}^t
\ddot{D}_{m;i,k}^{(r-1)}\tilde e_{N-t+k,N-t+j},
&H &=[\ddot{D}_{m;i,h}^{(r-2)}, \tilde e_{N-3t+h,N-t+j}].
\end{align*}
Now commute each of these elements in turn with 
$x:=  \tilde e_{N-t+f,N-2t+g}$ then apply $\pr_\chi$, 
using (\ref{rhodef}), (\ref{erel}), (\ref{chidef}) and the observation that
$x$ commutes with all elements of $U(\ddot{\mathfrak{p}})$, 
to deduce that
\begin{align*}
\pr_\chi([x, A]) &=
0,
\qquad\qquad\qquad\qquad\qquad\qquad\:
\pr_\chi([x, C]) =0,\\
\pr_\chi([x, B]) &=
- \delta_{f,j} \ddot{D}_{m;i,g}^{(r-1)}
\qquad\qquad\qquad\qquad
\pr_\chi([x, D]) =
\delta_{f,j} \ddot{D}_{m;i,g}^{(r-1)},\\
\pr_\chi([x, E]) &=
t \delta_{f,j} \ddot{D}_{m;i,g}^{(r-2)}
- \ddot{D}_{m;i,g}^{(r-2)} \tilde e_{N-t+f,N-t+j}
+ \delta_{f,j} \sum_{c=1}^t \ddot{D}_{m;i,c}^{(r-2)}
\tilde e_{N-2t+c,N-2t+g},\\
\pr_\chi([x, F]) &=
\delta_{f,j} [ \ddot{D}_{m;i,h}^{(r-2)}, \tilde e_{N-3t+h,N-2t+g}],\\
\pr_\chi([x, G]) &=
\ddot{D}_{m;i,g}^{(r-2)} \tilde e_{N-t+f,N-t+j}
-t \delta_{f,j} \ddot{D}_{m;i,g}^{(r-2)}
- \delta_{f,j} \sum_{c=1}^t \ddot{D}_{m;i,c}^{(r-2)}
\tilde e_{N-2t+c,N-2t+g}
,\\
\pr_\chi([x, H]) &=-\delta_{f,j} [ \ddot{D}_{m;i,h}^{(r-2)}, \tilde e_{N-3t+h,N-2t+g}].
\end{align*}
These sum to zero, hence $\pr_\chi([x, {D}_{m;i,j}^{(r)}]) = 0$.

(ii) Similar.
\end{proof}

\begin{Lemma}\label{sofar}
The following elements 
of $U(\mathfrak{p})$ are invariant under the twisted action of
$\mathfrak{m}$:
\begin{itemize}
\item[(i)] $D_{a;i,j}^{(r)}$ 
for  $1 \leq a \leq m$, $1 \leq i,j \leq \nu_a$ and $r > 0$;
\item[(ii)] $E_{a;i,j}^{(r)}$ for 
$1 \leq a < m$, $1 \leq i \leq \nu_a$, $1 \leq j \leq \nu_{a+1}$
and $r > s_{a,a+1}(\nu)$;
\item[(iii)] $F_{a;i,j}^{(r)}$ 
 for 
$1 \leq a< m$, $1 \leq i \leq \nu_{a+1}$, $1 \leq j \leq \nu_{a}$
and $r > s_{a+1,a}(\nu)$.
\end{itemize}
\end{Lemma}

\begin{proof}
This is just a matter of assembling the pieces.
Lemma~\ref{hor1} covers all the elements except
$D_{m;i,j}^{(r)}$
and (assuming $t < n$) $E_{m-1;i,j}^{(r)}$.
Since $\mathfrak{m}$ is generated by
$\dot{\mathfrak{m}}$ and the elements
$\tilde e_{f,g}$ for all $1 \leq f,g \leq N$ with
$\col(f) = l, \col(g) = l-1$,
Lemma~\ref{hor2} reduces the problem to showing that
the elements $D_{m;i,j}^{(r)}$ for $r> 0$ and (assuming $t < n$)
$E_{m-1;i,j}^{(r)}$ for $r > s_{m-1,m}(\nu)$ 
are invariant under all such $\tilde e_{f,g}$.
Suppose first that $t=n$. 
Then the required invariance is checked in Lemma~\ref{expinvs}(i)
and Lemma~\ref{hor4}(i).
Now assume that $t < n$.
If $s_{m-1,m}(\nu) = 1$, then
the invariance of $D_{m;i,j}^{(1)},
D_{m;i,j}^{(2)}$ and $E_{m-1;i,j}^{(2)}$ is checked in
Lemma~\ref{expinvs}.
The invariance of all higher
$D_{m;i,j}^{(r)}$ and 
$E_{m-1;i,j}^{(r)}$ then follows
by Lemma~\ref{hor1}, Lemma~\ref{hor3} and induction on $r$.
Finally if $s_{m-1,m}(\nu) > 1$ then the invariance
of $D_{m;i,j}^{(1)}$ is checked in Lemma~\ref{expinvs}(i), and 
the remaining elements are covered by Lemma~\ref{hor4}.
\end{proof}

Now we complete the proof of the induction step
in case $\Delta_{\rt}$. 
By the induction hypothesis, we can identify
the shifted Yangian $Y_{n,l-1}(\dot\sigma)$
with $W(\dot\pi)\subseteq U(\dot{\mathfrak{p}})$, so that the generators
$\dot D_{a;i,j}^{(r)}, \dot E_{a;i,j}^{(r)}$
and $\dot F_{a;i,j}^{(r)}$ of 
$Y_{n,l-1}(\dot\sigma)$
coincide with the elements of $W(\dot\pi)$ 
with the same name.
Theorem~\ref{miuramain} then
shows that there is an injective algebra homomorphism
$\Delta_{\rt}:Y_{n,l}(\sigma) \rightarrow 
U(\dot{\mathfrak{p}})
\otimes U(\mathfrak{gl}_t)$.
Observe moreover comparing Theorem~\ref{miuramain}
with Lemma~\ref{centbase} (like in the proof of Corollary~\ref{onto}) that
for each $d \geq 0$,
\begin{equation}\label{dim1}
 \dim \Delta_{\rt} (\F_d Y_{n,l}(\sigma))
=\dim \F_d Y_{n,l}(\sigma)
= \dim \F_d S(\mathfrak{c}_{\mathfrak{g}}(e)),
\end{equation}
where $\F_d S(\mathfrak{c}_{\mathfrak{g}}(e))$
denotes the sum of all the graded pieces of 
$S(\mathfrak{c}_{\mathfrak{g}}(e))$ of degree $\leq d$
in the Kazhdan grading.
Define elements $E_{a,b;i,j}^{(r)}$ and $F_{a,b;i,j}^{(r)}$
of $\F_r U(\mathfrak{p})$
recursively by the formulae (\ref{higheres})--(\ref{higherfs}), making 
arbitrary choices for the integers $k$ appearing there.
Let $X_d$ denote the subspace of $U(\mathfrak{p})$ spanned by
all monomials in
\begin{align*}
&\{D_{a;i,j}^{(r)}\}_{1 \leq a \leq m,
1 \leq i,j \leq \nu_a, 0 < r \leq s_a(\nu)},\\
&\{E_{a,b;i,j}^{(r)}\}_{1 \leq a < b\leq m,
1 \leq i \leq \nu_a, 1 \leq j \leq \nu_b,
s_{a,b}(\nu) < r \leq s_{a,b}(\nu)+p_a(\nu)},\\
&\{F_{a,b;i,j}^{(r)}\}_{1 \leq a < b\leq m,
1 \leq i \leq \nu_b, 1 \leq j \leq \nu_a,
s_{b,a}(\nu) < r \leq s_{b,a}(\nu)+p_a(\nu)}
\end{align*}
taken in some fixed order and of total degree $\leq d$.
By Lemma~\ref{sofar}, $X_d$ is actually 
a subspace of $\F_d W(\pi)$.
Define an algebra homomorphism
$\varphi_{\rt}:U(\mathfrak{p}) \rightarrow U(\dot{\mathfrak{p}}) \otimes
U(\mathfrak{gl}_t)$
by 
\begin{equation}\label{phirt}
\varphi_{\rt}(\tilde e_{i,j}) = \left\{
\begin{array}{ll}
\tilde e_{i,j} \otimes 1&\hbox{if $\col(i) \leq \col(j) 
\leq l-1$,}\\
0&\hbox{if $\col(i) \leq l-1, l=\col(j)$,}\\
1 \otimes \tilde e_{i-N+t,j-N+t}&\hbox{if $l=\col(i)=\col(j)$,}
\end{array}\right.
\end{equation}
where 
$\tilde e_{i-N+t,j-N+t} \in U(\mathfrak{gl}_t)$
denotes $e_{i-N+t,j-N+t} + \delta_{i,j}(n-t)$ like in Theorem~\ref{baby}.
By Lemma~\ref{superbaby}, we have that
\begin{align}\label{barbie1}
\varphi_{\rt}(D_{a;i,j}^{(r)}) &= \dot D_{a;i,j}^{(r)} \otimes 1
+ \delta_{a,m} \sum_{k=1}^{t}\dot D_{a;i,k}^{(r-1)} \otimes \tilde e_{k,j},\\
\varphi_{\rt}(E_{a;i,j}^{(r)}) &= \dot E_{a;i,j}^{(r)} \otimes 1
+ \delta_{a,m-1} \sum_{k=1}^{t}\dot E_{a;i,k}^{(r-1)} \otimes \tilde e_{k,j},\\
\varphi_{\rt}(F_{a;i,j}^{(r)}) &= \dot F_{a;i,j}^{(r)} \otimes 1.\label{barbie3}
\end{align}
Comparing this with Theorem~\ref{baby}(i) and recalling the
PBW basis for $\F_d Y_{n,l}(\sigma)$ from
Theorem~\ref{miuramain} and 
Corollary~\ref{newpbw},
we see that 
$\varphi_{\rt}(X_d)
=
\Delta_{\rt}(\F_d Y_{n,l}(\sigma))$.
Combining this with (\ref{dim1}) and
Corollary~\ref{phewy},
we deduce that
$$
\dim \F_d S(\mathfrak{c}_{\mathfrak{g}}(e))
= \dim \varphi_{\rt}(X_d) \leq \dim X_d 
\leq \dim \F_d W(\pi) \leq \dim \F_d S(\mathfrak{c}_{\mathfrak{g}}(e)).
$$
Hence equality holds everywhere, so we get that
 $X_d = \F_d W(\pi)$ for each $d \geq 0$, and the
map $\varphi_\rt:W(\pi) \rightarrow U(\dot{\mathfrak{p}}) \otimes U(\mathfrak{gl}_t)$ is injective
with the same image as $\Delta_\rt:Y_{n,l}(\sigma)
\rightarrow U(\dot{\mathfrak{p}}) \otimes U(\mathfrak{gl}_t)$.
In particular this shows that the elements listed in Theorem~\ref{main}
generate $W(\pi)$,
and the map $\varphi_{\rt}^{-1} \circ \Delta_{\rt}:Y_{n,l}(\sigma)
\rightarrow W(\pi)$ is exactly the filtered algebra isomorphism
described in the statement of Theorem~\ref{main}.
This completes the proof of the induction step
in the case $\Delta_{\rt}$. 

Finally we must sketch the proof of the induction step in the case
$\Delta_{\lt}$. So assume that either $t=n$ or $s_{n-t+1,n-t} = s_{m,m-1}(\nu) \neq 0$;
in particular, $q_1 \leq q_l$.
This time it is 
notationally convenient to number of the bricks of the pyramid $\pi$ 
down columns from {\em right to left}, not 
from left to right as before. 
For instance, the entries down the first 
(leftmost) column of the pyramid are $(N-t+1),(N-t+2),\dots,N$ in this new 
numbering. 
Let $\dot\pi$ be the 
pyramid
obtained from $\pi$ by removing this column from the pyramid $\pi$.
Let $\dot\sigma = (\dot s_{i,j})_{1 \leq i,j \leq n}$ be the
shift matrix corresponding to the pyramid $\dot\pi$ defined by (\ref{sl}).
Define $\dot{\mathfrak{p}}, \dot{\mathfrak{m}}$ and $\dot e$
in $\dot{\mathfrak{g}} = \mathfrak{gl}_{N-t}$
according to (\ref{Ac1})--(\ref{edef}).
This time, we can work simply with the natural embedding
of $U(\dot{\mathfrak{g}})$ into $U(\mathfrak{g})$ induced by
the usual embedding of $\dot{\mathfrak{g}}$ into $\mathfrak{g}$, since this
already maps
the elements $\tilde e_{i,j}$ of $U(\dot{\mathfrak{g}})$ 
to the elements $\tilde e_{i,j}$ of $U(\mathfrak{g})$
for $1 \leq i,j \leq N-t$.
The algebra $W(\dot\pi) = U(\dot{\mathfrak{p}})^{\dot{\mathfrak{m}}}$
is a subalgebra of $U(\dot{\mathfrak{p}}) \subset U(\mathfrak{p})$, 
$\dot{\mathfrak{m}}$ is a subalgebra of
$\mathfrak{m}$, and the twisted action of
$\dot{\mathfrak{m}}$ on $U(\dot{\mathfrak{p}})$ is the restriction of the
twisted action of $\mathfrak{m}$ on $U(\mathfrak{p})$.
Let $\dot D_{a;i,j}^{(r)}, \dot{{\widetilde D}}_{a;i,j}^{(r)}$,
$\dot E_{a;i,j}^{(r)}$ and $\dot F_{a;i,j}^{(r)}$ denote the
elements of $U(\dot{\mathfrak{p}})$ as defined in $\S$\ref{sinvariants}
relative to the shape $\nu$.
The all-important analogue of Lemma~\ref{superbaby} is as follows.

\begin{Lemma}\label{superbaby2}
The following equations hold for $r > 0$, all admissible $a,i,j$ and 
any fixed $1 \leq l \leq t$:
\begin{align}\label{sb21}
D_{a;i,j}^{(r)} & = \dot D_{a;i,j}^{(r)} + \delta_{a,m}
\left(\sum_{k=1}^t
\tilde e_{N-t+i,N-t+k} \dot D_{a;k,j}^{(r-1)}
+[\tilde e_{N-t+i,N-2t+h}, \dot D_{a;h,j}^{(r-1)}]\right),\\
E_{a;i,j}^{(r)} & = \dot E_{a;i,j}^{(r)},\\
F_{a;i,j}^{(r)} & = \dot F_{a;i,j}^{(r)} + \delta_{a,m-1}\left(
\sum_{k=1}^t \tilde e_{N-t+i,N-t+k} \dot F_{a;k,j}^{(r-1)}
+ [\tilde e_{N-t+i,N-2t+h},\dot F_{a;h,j}^{(r-1)}]\right),\label{sb22}
\end{align}
where for (\ref{sb22}) we are assuming that $r > 1$ if $a=m-1$.
\end{Lemma}

Combining this with the induction hypothesis and imitating 
the arguments in Lemmas~\ref{hor1}--\ref{sofar} one can now check:

\begin{Lemma}\label{sofar2}
The following elements 
of $U(\mathfrak{p})$ are invariant under the twisted action of
$\mathfrak{m}$:
\begin{itemize}
\item[(i)] $D_{a;i,j}^{(r)}$ 
for  $1 \leq a \leq m$, $1 \leq i,j \leq \nu_a$ and $r > 0$;
\item[(ii)] $E_{a;i,j}^{(r)}$ for 
$1 \leq a < m$, $1 \leq i \leq \nu_a$, $1 \leq j \leq \nu_{a+1}$
and $r > s_{a,a+1}(\nu)$;
\item[(iii)] $F_{a;i,j}^{(r)}$ 
 for 
$1 \leq a< m$, $1 \leq i \leq \nu_{a+1}$, $1 \leq j \leq \nu_{a}$
and $r > s_{a+1,a}(\nu)$.
\end{itemize}
\end{Lemma}

Finally,
define an algebra homomorphism
$\varphi_{\lt}:U(\mathfrak{p}) \rightarrow U(\mathfrak{gl}_t) 
\otimes U(\dot{\mathfrak{p}})$
by 
\begin{equation}\label{philt}
\varphi_{\lt}(\tilde e_{i,j}) = \left\{
\begin{array}{ll}
\tilde e_{i-N+t,j-N+t}\otimes 1&\hbox{if $\col(i)=\col(j)=1$,}\\
0&\hbox{if $\col(i) =1, 2 \leq \col(j)$,}\\
1 \otimes \tilde e_{i,j}&\hbox{if $2 \leq\col(i) \leq \col(j)$,}
\end{array}\right.
\end{equation}
where 
$\tilde e_{i-N+t,j-N+t} \in U(\mathfrak{gl}_t)$
denotes $e_{i-N+t,j-N+t} + \delta_{i,j}(n-t)$ like in Theorem~\ref{baby}.
By Lemma~\ref{superbaby2}, we have that
\begin{align}\label{cindy1}
\varphi_{\lt}(D_{a;i,j}^{(r)}) &= 1 \otimes
\dot D_{a;i,j}^{(r)}
+ \delta_{a,m} \sum_{k=1}^{t} \tilde e_{i,k} \otimes \dot D_{a;k,j}^{(r-1)},\\
\varphi_{\lt}(E_{a;i,j}^{(r)}) &= 1 \otimes \dot E_{a;i,j}^{(r)}
+ \delta_{a,m-1} \sum_{k=1}^{t}\tilde e_{i,k} \otimes \dot E_{a;k,j}^{(r-1)},\\
\varphi_{\lt}(F_{a;i,j}^{(r)}) &= 1 \otimes \dot F_{a;i,j}^{(r)}.
\label{cindy3}
\end{align}
If we identify $Y_{n,l-1}(\dot\sigma)$
with $W(\dot\pi) \subseteq U(\dot{\mathfrak{p}})$ using the induction 
hypothesis,
these are the same as the images of the corresponding
elements of $Y_{n,l}(\sigma)$ under the baby comultiplication
$\Delta_\lt$ from Theorem~\ref{baby}(ii).
So now the proof of the 
induction step in the case $\Delta_{\lt}$ can be completed
like before.
Theorem~\ref{main} is proved.

\section{Grown-up comultiplication}\label{sgrownup}

Fix a pyramid $\pi$ of height $\leq n$ with column heights
$(q_1,\dots,q_l)$.
Throughout the section, we will work with the numbering
of the bricks of $\pi$ down columns from left to right.
Define $\mathfrak{p}, \mathfrak{h}$, $\mathfrak{m}$ and $e$
from (\ref{Ac1}), and
let $W(\pi) := U(\mathfrak{p})^{\mathfrak{m}}$
be the corresponding finite $W$-algebra.
Suppose we are given $l',l'' \geq 0$
with $l'+l'' = l$.
Let $\pi'$ and $\pi''$
denote the pyramids consisting just of
the leftmost $l'$ and the rightmost $l''$
columns of $\pi$, respectively.
We write $\pi = \pi' \otimes \pi''$ whenever
a pyramid is split in this way; for example,
$$
\begin{picture}(90, 60)%
\put(0,0){\line(1,0){75}}
\put(0,15){\line(1,0){75}}
\put(15,30){\line(1,0){45}}
\put(30,45){\line(1,0){30}}
\put(30,60){\line(1,0){15}}
\put(0,0){\line(0,1){15}}
\put(15,0){\line(0,1){30}}
\put(30,0){\line(0,1){60}}
\put(45,0){\line(0,1){60}}
\put(60,0){\line(0,1){45}}
\put(75,0){\line(0,1){15}}
\put(7,8){\makebox(0,0){1}}
\put(22,8){\makebox(0,0){3}}
\put(37,8){\makebox(0,0){7}}
\put(52,8){\makebox(0,0){10}}
\put(67,8){\makebox(0,0){11}}
\put(22,23){\makebox(0,0){2}}
\put(37,23){\makebox(0,0){6}}
\put(52,23){\makebox(0,0){9}}
\put(37,38){\makebox(0,0){5}}
\put(52,38){\makebox(0,0){8}}
\put(37,53){\makebox(0,0){4}}
\end{picture}
\quad
\begin{picture}(45, 60)%
\put(0,0){\line(1,0){45}}
\put(0,15){\line(1,0){45}}
\put(15,30){\line(1,0){30}}
\put(30,45){\line(1,0){15}}
\put(30,60){\line(1,0){15}}
\put(0,0){\line(0,1){15}}
\put(15,0){\line(0,1){30}}
\put(30,0){\line(0,1){60}}
\put(45,0){\line(0,1){60}}
\put(7,8){\makebox(0,0){1}}
\put(22,8){\makebox(0,0){3}}
\put(37,8){\makebox(0,0){7}}
\put(22,23){\makebox(0,0){2}}
\put(37,23){\makebox(0,0){6}}
\put(37,38){\makebox(0,0){5}}
\put(37,53){\makebox(0,0){4}}
\put(53,7){\makebox(0,0){$\otimes$}}
\put(-12,6){\makebox(0,0){$=$}}
\end{picture}
\begin{picture}(70, 60)%
\put(15,0){\line(1,0){30}}
\put(15,15){\line(1,0){30}}
\put(15,30){\line(1,0){15}}
\put(15,45){\line(1,0){15}}
\put(30,0){\line(0,1){45}}
\put(15,0){\line(0,1){45}}
\put(45,0){\line(0,1){15}}
\put(22,8){\makebox(0,0){3}}
\put(37,8){\makebox(0,0){4}}
\put(22,23){\makebox(0,0){2}}
\put(22,38){\makebox(0,0){1}}
\end{picture}
\begin{picture}(30, 60)%
\put(0,0){\line(1,0){15}}
\put(0,15){\line(1,0){15}}
\put(0,0){\line(0,1){15}}
\put(15,0){\line(0,1){15}}
\put(7,8){\makebox(0,0){1}}
\put(23,7){\makebox(0,0){$\otimes$}}
\put(53,7){\makebox(0,0){$\otimes$}}
\put(83,7){\makebox(0,0){$\otimes$}}
\put(113,7){\makebox(0,0){$\otimes$}}
\put(-12,6){\makebox(0,0){$=$}}
\end{picture}
\begin{picture}(30, 60)%
\put(0,0){\line(1,0){15}}
\put(0,15){\line(1,0){15}}
\put(0,30){\line(1,0){15}}
\put(0,0){\line(0,1){30}}
\put(15,0){\line(0,1){30}}
\put(7,8){\makebox(0,0){2}}
\put(7,23){\makebox(0,0){1}}
\end{picture}
\begin{picture}(30, 60)%
\put(0,0){\line(1,0){15}}
\put(0,15){\line(1,0){15}}
\put(0,30){\line(1,0){15}}
\put(0,45){\line(1,0){15}}
\put(0,60){\line(1,0){15}}
\put(0,0){\line(0,1){60}}
\put(15,0){\line(0,1){60}}
\put(7,8){\makebox(0,0){4}}
\put(7,23){\makebox(0,0){3}}
\put(7,38){\makebox(0,0){2}}
\put(7,53){\makebox(0,0){1}}
\end{picture}
\begin{picture}(30, 60)%
\put(0,0){\line(1,0){15}}
\put(0,15){\line(1,0){15}}
\put(0,30){\line(1,0){15}}
\put(0,45){\line(1,0){15}}
\put(0,0){\line(0,1){45}}
\put(15,0){\line(0,1){45}}
\put(7,8){\makebox(0,0){3}}
\put(7,23){\makebox(0,0){2}}
\put(7,38){\makebox(0,0){1}}
\end{picture}
\begin{picture}(30, 60)%
\put(0,0){\line(1,0){15}}
\put(0,15){\line(1,0){15}}
\put(0,0){\line(0,1){15}}
\put(15,0){\line(0,1){15}}
\put(7,8){\makebox(0,0){1}}
\end{picture}
$$
Let 
${\mathfrak{p}'}, {\mathfrak{m}}', e'$
and 
${\mathfrak{p}}'', {\mathfrak{m}}'', e''$
be defined 
from the pyramids $\pi'$ and $\pi''$ 
inside the Lie algebras
${\mathfrak{g}}' = \mathfrak{gl}_{N'}$
and
${\mathfrak{g}}'' = \mathfrak{gl}_{N''}$, respectively.
So $N = N'+N''$.
Let 
$W(\pi') = U(\mathfrak{p}')^{{\mathfrak{m}'}}$
and $W(\pi'') = U({\mathfrak{p}''})^{{\mathfrak{m}''}}$.
Recall also the elements $\tilde e_{i,j}$ of 
$U(\mathfrak{p}), U(\mathfrak{p}')$ and $U(\mathfrak{p}'')$
defined by (\ref{etildedef})
working with the pyramids $\pi$, $\pi'$ and $\pi''$, respectively.
Define a homomorphism
$\Delta_{l',l''}:U(\mathfrak{p}) \rightarrow U(\mathfrak{p}')
\otimes U(\mathfrak{p}'')$
by declaring that
\begin{equation}\label{dede}
\Delta_{l',l''}(\tilde e_{i,j}) = \left\{
\begin{array}{ll}
\tilde e_{i,j}\otimes 1&\hbox{if $\col(j) \leq l'$,}\\
0&\hbox{if $\col(i) \leq l', l'+1 \leq \col(j)$,}\\
1 \otimes \tilde e_{i-N',j-N'}&\hbox{if $l'+1 \leq \col(i)$,}\\
\end{array}
\right.
\end{equation}
for all $1 \leq i,j \leq N$ with $\col(i) \leq \col(j)$.
This map is obviously filtered with respect to the
Kazhdan filtrations.
The following lemma describes the effect of 
$\Delta_{l',l''}$ in terms of the elements $T_{i,j;0}^{(r)}$ 
of $U(\mathfrak{p})$, $U(\mathfrak{p}')$ and
$U(\mathfrak{p}'')$ defined as in $\S$\ref{sinvariants}, again
working with the pyramids $\pi, \pi'$ and $\pi''$ respectively.
This should be compared with (\ref{comult}).

\begin{Lemma}\label{howithappens}
For $1 \leq i,j \leq n$ and $r > 0$,
$\displaystyle \Delta_{l',l''}(T_{i,j;0}^{(r)})
=
\sum_{s=0}^r \sum_{k=1}^n T_{i,k;0}^{(s)} \otimes T_{k,j;0}^{(r-s)}.$
\end{Lemma}
\begin{proof}
Clear from (\ref{thedef}).
\end{proof}

One should visualize the map $\Delta_{l',l''}$ as follows.
The standard embedding of
$\mathfrak{g}' \oplus \mathfrak{g}''$
into $\mathfrak{g}$ also embeds
$\mathfrak{p}' \oplus \mathfrak{p}''$ into 
$\mathfrak{p}$ and $\mathfrak{m}' \oplus \mathfrak{m}''$
into $\mathfrak{m}$. Then the map $\Delta_{l',l''}$
is just induced by the obvious projection
$\mathfrak{p} \twoheadrightarrow \mathfrak{p}' \oplus \mathfrak{p}''$
followed by a constant shift.
The character $\chi' \oplus \chi''$ of $\mathfrak{m}' \oplus \mathfrak{m}''$
defined by taking the trace form with $e'+e''$
is the restriction of the character $\chi$ of $\mathfrak{m}$.
So the map $\Delta_{l',l''}$ sends
twisted $\mathfrak{m}$-invariants in $U(\mathfrak{p})$
to twisted $(\mathfrak{m}' \oplus \mathfrak{m}'')$-invariants
in $U(\mathfrak{p}') \otimes U(\mathfrak{p}'')$.
This shows that the restriction of $\Delta_{l',l''}$ defines an
algebra homomorphism
\begin{equation}\label{comulty}
\Delta_{l',l''}:W(\pi) \rightarrow W(\pi') \otimes W(\pi''),
\end{equation}
which is again a filtered map with respect to the
Kazhdan filtrations.
This defines a comultiplication $\Delta_{l',l''}$ 
between the finite $W$-algebras
which is 
coassociative in the following sense:

\begin{Lemma}\label{coass}
If $\pi = \pi' \otimes \pi'' \otimes \pi'''$ is a pyramid,
where $\pi', \pi''$ and $\pi'''$ are of levels
$l',l''$ and $l'''$ respectively, then the
following diagram commutes:
\begin{equation*}
\begin{CD}
W(\pi' \otimes \pi'' \otimes \pi'') & @>\Delta_{l'+l'',l'''}>> & W(\pi' \otimes \pi'') \otimes W(\pi''')\\
@V\Delta_{l',l''+l'''} VV&&@VV\Delta_{l',l''} \otimes 1V\\
W(\pi') \otimes W(\pi'' \otimes \pi''') &@>1 \otimes \Delta_{l'',l'''}>> &W(\pi') \otimes W(\pi'')
\otimes W(\pi''')
\end{CD}
\end{equation*}
\end{Lemma}

\begin{proof}
Obvious from (\ref{dede}).
\end{proof}

Read off a shift matrix $\sigma = (s_{i,j})_{1 \leq i,j \leq n}$
from the pyramid $\pi$ according to (\ref{pyr3}).
From now on we will identify the algebra
$W(\pi)$ with the shifted Yangian $Y_{n,l}(\sigma)$ of level $l$
using the isomorphism from Theorem~\ref{main}.
Let $t := \min(q_1,q_l)$.
Suppose first that 
the baby comultiplication $\Delta_{\rt}:
Y_{n,l}(\sigma)\rightarrow Y_{n,l-1}(\dot\sigma) \otimes U(\mathfrak{gl}_t)$ 
from (\ref{deltar}) is defined, i.e.
either $t=n$ or $s_{n-t,n-t+1} \neq 0$.
The pyramid corresponding to the level $(l-1)$
and the shift matrix $\dot\sigma$
here (which is defined by (\ref{sr})) is simply the pyramid $\pi'$ 
obtained from $\pi$ by removing the rightmost
column. So, identifying $Y_{n,l-1}(\dot\sigma)$ with $W(\pi')$
according to Theorem~\ref{main} again,
the map $\Delta_{\rt}$ is therefore identified with a map
$\Delta_{\rt}:W(\pi) \rightarrow W(\pi') \otimes U(\mathfrak{gl}_t)$,
just like the map $\Delta_{l-1,1}$.
Suppose instead that the baby comultiplication
$\Delta_{\lt}:Y_{n,l}(\sigma) \rightarrow U(\mathfrak{gl}_t)
\otimes Y_{n,l-1}(\dot\sigma)$ from (\ref{deltal})
is defined, i.e.
either $t=n$ or $s_{n-t+1,n-t} \neq 0$.
This time, 
the pyramid corresponding to the level $(l-1)$ and the
shift matrix $\dot\sigma$ (defined now by (\ref{sl}))
is the pyramid $\pi''$ 
obtained from $\pi$ by removing the leftmost
column.
So $\Delta_{\lt}$ is identified with a map
$W(\pi) \rightarrow U(\mathfrak{gl}_t) \otimes W(\pi'')$,
just like the map $\Delta_{1,l-1}$.

\begin{Lemma}\label{speccase}
Whenever the baby comultiplications 
$\Delta_{\rt}$ and $\Delta_{\lt}$ are defined, they are equal to
$\Delta_{l-1,1}$ and $\Delta_{1,l-1}$,
respectively.
\end{Lemma}

\begin{proof}
This is obvious in the case $l=1$, so assume that $l > 1$.
Assume the map $\Delta_{\rt}$ is defined.
Comparing (\ref{dede}) with (\ref{phirt}), it is clear
that the map $\Delta_{l-1,1}$ coincides with the map
$\varphi_{\rt}$ defined in the proof of Theorem~\ref{main}.
Comparing (\ref{barbie1})--(\ref{barbie3}) with
Theorem~\ref{baby}(i), the map $\varphi_{\rt}$
coincides with $\Delta_\rt$ under the identifications
of $W(\pi)$ with $Y_{n,l}(\sigma)$ and $W(\pi')$ 
with $Y_{n,l-1}(\dot\sigma)$.
Hence, $\Delta_{l-1,1} = \Delta_{\rt}$.
The proof that $\Delta_{1,l-1} = \Delta_{\lt}$ is similar,
using (\ref{philt}), (\ref{cindy1})--(\ref{cindy3})
and Theorem~\ref{baby}(ii).
\end{proof}

If we iterate the comultiplication
(in any order by coassociativity) a total of $(l-1)$ times to
split the pyramid $\pi$ into its individual columns, we obtain a 
homomorphism
\begin{equation}\label{miura}
\mu:W(\pi) \rightarrow U(\mathfrak{gl}_{q_1})
\otimes \cdots \otimes U(\mathfrak{gl}_{q_{l}}).
\end{equation}
Let us give a direct description of this map.
The elements
$\{e_{i,j}^{[r]}\}_{1 \leq r \leq l, 1 \leq i,j \leq q_r}$
defined from $e_{i,j}^{[r]} := e_{q_1+\cdots+q_{r-1}+i,q_1+\cdots+q_{r-1}+j}$
form a basis for the Levi subalgebra 
$\mathfrak{h}$ of $\mathfrak{p}$. 
Identify $U(\mathfrak{h})$
with
$U(\mathfrak{gl}_{q_1})
\otimes \cdots \otimes U(\mathfrak{gl}_{q_{l}})$
so that $e_{i,j}^{[r]}$ is identified with
$1^{\otimes (r-1)} \otimes e_{i,j} \otimes 1^{\otimes (l-r)}$.
The map (\ref{miura}) is then identified with a homomorphism
\begin{equation}\label{miura2}
\mu:W(\pi) \rightarrow U(\mathfrak{h})
\end{equation}
which is a filtered map with respect to the
Kazhdan filtration of $W(\pi)$ and the
standard filtration of
$\F_0 U(\mathfrak{h}) \subseteq \F_1 U(\mathfrak{h})
\subseteq \cdots$ of
$U(\mathfrak{h})$; we will write $\gr U(\mathfrak{h})$
for the associated graded algebra here.
Let
\begin{equation}\label{etadef}
\eta:U(\mathfrak{h}) \rightarrow U(\mathfrak{h}),
\quad
e_{i,j}^{[r]} \mapsto e_{i,j}^{[r]} + \delta_{i,j} (q_{r+1}+\cdots
+ q_l).
\end{equation}
Let $\xi:U(\mathfrak{p}) \twoheadrightarrow U(\mathfrak{h})$
be the algebra homomorphism induced by the natural projection
$\mathfrak{p} \twoheadrightarrow \mathfrak{h}$.
Then, it is easy to see from (\ref{dede}) and (\ref{rhodef})--(\ref{etildedef})
that $\mu$ is precisely
the restriction of the map
$\eta\circ\xi$ to $W(\pi)$.
By analogy with the language used in \cite{BT},
we call $\mu$ the {\em Miura transform}; in
\cite{Ly} it is called the {\em generalized Harish-Chandra
homomorphism}.
The following result is due to Lynch \cite[Corollary 2.3.2]{Ly}.

\begin{Theorem}
The map $\gr \mu:\gr W(\pi) \rightarrow
\gr U(\mathfrak{h})$ is injective, hence so is the Miura transform
$\mu:W(\pi) \rightarrow U(\mathfrak{h})$ itself.
\end{Theorem}

\begin{proof}
By Lemma~\ref{speccase} and the definition (\ref{miura}), 
we can factor $\mu$ as a
composition of $(l-1)$ maps of the form $\Delta_\rt$ or $\Delta_\lt$.
Now the theorem follows from
the injectivity of $\gr \Delta_\rt$ and $\gr \Delta_\lt$
proved in Theorem~\ref{miuramain}.
\end{proof}

\begin{Corollary}\label{injy}
The map
$\gr\Delta_{l',l''}:
\gr W(\pi) \rightarrow \gr W(\pi') \otimes W(\pi'')$ 
is injective for any $l'+l'' = l$.
Hence so is the comultiplication
$\Delta_{l',l''}:W(\pi) \rightarrow W(\pi') \otimes W(\pi'')$.
\end{Corollary}

\begin{proof}
We can factor $\mu$ as a composition of $\Delta_{l',l''}$
followed by $(l-2)$ more maps.
\end{proof}

For the next lemma, we recall from Corollary~\ref{oneiso}
that if $\pi$ and $\dot\pi$ are two pyramids with the same
row lengths, then there is a canonical isomorphism
$\iota:W(\pi) \rightarrow W(\dot\pi)$.
The following lemma shows that the comultiplication
is compatible with these isomorphisms.

\begin{Lemma}\label{cube}
Suppose that $\pi = \pi' \otimes \pi''$
and $\dot\pi = \dot\pi' \otimes \dot\pi''$ 
are pyramids such that $\dot\pi'$ and $\dot\pi''$ 
have the same
row lengths as $\pi'$ and $\pi''$, respectively.
Then the following diagram commutes:
$$
\begin{CD}
W(\pi) & @>\Delta_{l',l''}>> &W(\pi') \otimes W(\pi'')\\
@V\iota VV&&@VV\iota\otimes\iota V\\
W(\dot\pi) & @>\Delta_{l',l''}>> &W(\dot\pi') \otimes W(\dot\pi''),
\end{CD}
$$
where $\pi'$ and $\pi''$ are of levels $l'$ and $l''$,
respectively.
\end{Lemma}

\begin{proof}
Proceed by induction on $l = l'+l''$.
If either $l'=0$ or $l'' = 0$,  the statement of the lemma is
vacuous, so the base case $l=1$ is trivial.
Now suppose that $l',l'' > 0$, so $l > 1$, 
and that the lemma has been proved for all
smaller levels. Read off a shift matrix $\sigma$
from the pyramid $\pi$ and identify $W(\pi)$ with
$Y_{n,l}(\sigma)$ as usual.
At least one of the baby comultiplications $\Delta_{\rt}$
and $\Delta_{\lt}$ from (\ref{deltar})--(\ref{deltal})
is always defined;
we explain the proof of the induction step just 
in the case that $\Delta_{\rt}$ is defined, the argument being
entirely similar in the other case.
Also read off a shift matrix $\dot\sigma$ from the
pyramid $\dot\pi$ and identify $W(\dot\pi)$
with $Y_{n,l}(\dot\sigma)$; we can ensure in 
doing this that the baby comultiplication
$\Delta_{\rt}$ is defined for $Y_{n,l}(\dot\sigma)$
too.

Suppose first that $l'' = 1$.
Then by Lemma~\ref{speccase}, $\Delta_{l',l''}$ is equal in either case
to the map
$\Delta_{\rt}$, and the commutativity of the diagram
is easy to check explicitly on
generators of $W(\pi)$, using (\ref{iotadefpar}) and
Theorem~\ref{baby}(i).

Now suppose that $l'' > 1$.
Let $\rho, \dot\rho, \rho'$ and $\dot\rho'$ 
denote the pyramids obtained from 
$\pi,\dot\pi, \pi''$ and $\dot\pi''$, respectively,
by removing the rightmost column (which is of height
$q_l$ in all cases).
Consider the following cube
$$
\begin{picture}(290,130)
\put(10,10){\makebox(0,0){$W(\pi')\otimes W(\pi'')$}}
\put(180,10){\makebox(0,0){$W(\pi')\otimes W(\rho')\otimes U(\mathfrak{gl}_{q_l})$}}
\put(15,78){\makebox(0,0){$W(\pi)$}}
\put(180,78){\makebox(0,0){$W(\rho)\otimes U(\mathfrak{gl}_{q_l})$}}

\put(85,47){\makebox(0,0){$W(\dot\pi')\otimes W(\dot\pi'')$}}
\put(245,47){\makebox(0,0){$W(\dot\pi')\otimes W(\dot\rho')\otimes U(\mathfrak{gl}_{q_l})$}}
\put(72,115){\makebox(0,0){$W(\dot\pi)$}}
\put(240,115){\makebox(0,0){$W(\dot\rho)\otimes U(\mathfrak{gl}_{q_l})$}}
\put(43,33){\makebox(0,0){$\nearrow$}}
\put(43,33){\line(-1,-1){13}}
\put(43,103){\makebox(0,0){$\nearrow$}}
\put(43,103){\line(-1,-1){13}}
\put(203,33){\makebox(0,0){$\nearrow$}}
\put(203,33){\line(-1,-1){13}}
\put(203,103){\makebox(0,0){$\nearrow$}}
\put(203,103){\line(-1,-1){13}}
\put(10,27){\makebox(0,0){$\downarrow$}}
\put(9.9,27){\line(0,1){38}}
\put(170,27){\makebox(0,0){$\downarrow$}}
\put(169.9,27){\line(0,1){38}}
\put(110,79){\makebox(0,0){$\rightarrow$}}
\put(110,79.56){\line(-1,0){49}}
\put(110,9){\makebox(0,0){$\rightarrow$}}
\put(110,9.56){\line(-1,0){49}}
\put(70,64){\makebox(0,0){$\downarrow$}}
\put(69.9,64){\line(0,1){13}}
\put(69.9,82){\line(0,1){20}}
\put(225,64){\makebox(0,0){$\downarrow$}}
\put(224.9,64){\line(0,1){38}}
\put(178,114){\makebox(0,0){$\rightarrow$}}
\put(178,114.56){\line(-1,0){49}}
\put(178,46){\makebox(0,0){$\rightarrow$}}
\put(167,46.56){\line(-1,0){38}}
\end{picture}
$$
where the maps on the front and back faces are defined
from the comultiplications, and the remaining maps are 
isomorphisms built from $\iota$.
The front and back faces commute by Lemma~\ref{coass}.
The top and bottom faces commute by the special case
considered in the preceeding paragraph. The right hand
face commutes by the induction hypothesis.
Since the comultiplication map $W(\dot\pi')
\otimes W(\dot\pi'') \rightarrow W(\dot\pi') \otimes
W(\dot\rho') \otimes U(\mathfrak{gl}_{q_l})$ is injective
by Corollary~\ref{injy}, it therefore
follows that the left face commutes too. This completes the proof
of the induction step.
\end{proof}

In the remainder of the section, we are going to
lift the comultiplication 
to the shifted Yangian $Y_n(\sigma)$ itself.
First, we must explain one other basic operation
on pyramids, that of {\em column removal}.
Let
$1 \leq i_1 < \cdots < i_{\dot l} \leq l$ be some subset of the columns
of the pyramid $\pi$,
and let 
$\dot\pi$ be the pyramid with column heights
$(q_{i_1}, \dots, q_{i_{\dot l}})$.
We can read off shift matrices
$\sigma = (s_{i,j})_{1 \leq i,j \leq n}$ 
and $\dot\sigma = (\dot s_{i,j})_{1 \leq i,j \leq n}$
from the pyramids $\pi$ and $\pi'$ respectively
in such a way as to ensure $\dot s_{i,j} \leq s_{i,j}$
for all $1 \leq i,j \leq n$.
Hence embedding
$Y_n(\sigma)$ and $Y_n(\dot\sigma)$ into
$Y_n$ in the canonical way, we have that
$Y_n(\sigma) \subseteq Y_n(\dot\sigma)$.
This inclusion $Y_n(\sigma) \hookrightarrow Y_n(\dot\sigma)$
factors through the quotients to induce a map
\begin{equation}
\zeta:Y_{n,l}(\sigma) \rightarrow Y_{n,\dot l}(\dot\sigma).
\end{equation}
Equivalently, identifying
$Y_{n,l}(\sigma)$ with $W(\pi)$ and
$Y_{n,\dot l}(\dot\sigma)$ with $W(\dot\pi)$
by Theorem~\ref{main}, this defines a homomorphism
\begin{equation}
\zeta:W(\pi) \rightarrow W(\dot\pi)
\end{equation}
sending the generators $D_i^{(r)}, E_i^{(r)}$ and
$F_i^{(r)}$ of $W(\pi)$ 
to the elements
$\dot D_i^{(r)}, \dot E_i^{(r)}$
and $\dot F_i^{(r)}$ of $W(\dot\pi)$, respectively.
In order to understand the relationship between
this ``column removal homomorphism'' $\zeta$ and the comultiplication,
we need another description of $\zeta$.
Let $\mu:W(\pi) \rightarrow U(\mathfrak{h})$
and $\dot\mu:W(\dot\pi) \rightarrow U(\dot{\mathfrak{h}})$
be the Miura transforms defined by (\ref{miura2}).
There is an obvious projection
$\hat\zeta: \mathfrak{h} \twoheadrightarrow \dot{\mathfrak{h}}$
defined by
\begin{equation*}
\hat\zeta(e_{i,j}^{[r]}) =
\left\{
\begin{array}{ll}
e_{i,j}^{[s]}&\hbox{if $r = i_s$ for some $s=1,\dots,\dot l$,}\\
0&\hbox{otherwise.}
\end{array}
\right.
\end{equation*}

\begin{Lemma}\label{aftermu}
With the above notation, the following diagram commutes
$$
\begin{CD}
W(\pi) &@>\mu >> &U(\mathfrak{h}) \\
@V\zeta VV&&@VV\hat\zeta V\\
W(\dot\pi)&@>\dot\mu >> & U(\dot{\mathfrak{h}})
\end{CD}
$$
\end{Lemma}

\begin{proof}
One checks explicitly from Corollary~\ref{transcor} and the definition
(\ref{thedef}) that
$\hat\zeta$ maps the elements
$\mu(D_i^{(r)})$,
$\mu(E_i^{(r)})$ and $\mu(F_i^{(r)})$
to $\dot\mu(\dot D_i^{(r)})$,
$\dot\mu(\dot E_i^{(r)})$ and
$\dot\mu(\dot F_i^{(r)})$.
\end{proof}

\begin{Corollary}\label{stable}
Suppose that $\pi = \pi' \otimes \pi''$ 
and $\dot \pi =\dot\pi' \otimes \dot\pi''$,
where $\dot\pi'$ and $\dot\pi''$ are obtained by removing
columns from the pyramids $\pi'$ and $\pi''$, respectively.
Then, the following diagram commutes
$$
\begin{CD}
W(\pi) &@>\Delta_{l',l''}>> &W(\pi') \otimes W(\pi'')\\
@V\zeta VV&&@VV\zeta \otimes \zeta V\\
W(\dot\pi) &@>\Delta_{\dot l', \dot l''}>> &W(\dot\pi') \otimes W(\dot\pi'')
\end{CD}
$$
where $\pi',\pi''$ are of widths $l',l''$
and $\dot\pi',\dot\pi''$ are of widths
$\dot l', \dot l''$.
\end{Corollary}

\begin{proof}
In view of Lemma~\ref{aftermu} and the injectivity of the Miura
transforms, this follows from the commutativity of the following
diagram:
$$
\begin{CD}
U(\mathfrak{h}) &@>\sim >> &U(\mathfrak{h}') \otimes U(\mathfrak{h}'')\\
@V\hat\zeta VV&&@VV\hat\zeta \otimes \hat\zeta V\\
U(\dot{\mathfrak{h}}) &@>\sim >> &U(\dot{\mathfrak{h}}') \otimes
U(\dot{\mathfrak{h}}'')
\end{CD}
$$
where the horizontal maps are induced by the 
obvious isomorphisms
$\mathfrak{h} \cong \mathfrak{h}' \oplus
\mathfrak{h}''$ and
$\dot{\mathfrak{h}} \cong \dot{\mathfrak{h}}' \oplus
\dot{\mathfrak{h}}''$.
\end{proof}

Now we can prove the main theorem of the section:

\begin{Theorem}\label{grownup}
Let $\sigma$ be a shift matrix and write $\sigma = \sigma'
+ \sigma''$ where $\sigma'$ is strictly lower triangular and
$\sigma''$ is strictly upper triangular. Embedding
$Y_n(\sigma), Y_n(\sigma')$ and $Y_n(\sigma'')$ into $Y_n$
in the standard way, the restriction of the comultiplication
$\Delta:Y_n\rightarrow Y_n \otimes Y_n$ gives a 
homomorphism
$$
\Delta:Y_n(\sigma)\rightarrow Y_n(\sigma') \otimes Y_n(\sigma'').
$$
Moreover, for $l' \geq s_{n,1}, l'' \geq s_{1,n}$ and $l = l' + l''$,
this map $\Delta$ factors through
the quotients to define a homomorphism
$Y_{n,l}(\sigma)
\rightarrow Y_{n,l'}(\sigma') \otimes Y_{n,l''}(\sigma'')$
which, on identifying
$Y_{n,l}(\sigma)$ with $W(\pi)$, 
$Y_{n,l'}(\sigma')$ with $W(\pi')$ and
$Y_{n,l''}(\sigma'')$ with $W(\pi'')$ as usual,
is precisely the comultiplication $\Delta_{l',l''}$
from (\ref{comulty}).
\end{Theorem}

\begin{proof}
Start from the map $\Delta_{l',l''}:Y_{n,l}(\sigma)
\rightarrow Y_{n,l'}(\sigma') \otimes Y_{n,l''}(\sigma'')$ from
(\ref{comulty}).
Recall from Remark~\ref{invsys} how $Y_n(\sigma), Y_n(\sigma')$ and
$Y_n(\sigma'')$ are identified with the inverse limits 
$\varprojlim Y_{n,l}(\sigma), \varprojlim Y_{n,l'}(\sigma')$
and 
$\varprojlim Y_{n,l''}(\sigma'')$,
respectively. Corollary~\ref{stable} is exactly what is needed to ensure that
the maps $\Delta_{l',l''}$ are stable as $l',l'' \rightarrow \infty$.
Hence there is an induced homomorphism
$\varprojlim \Delta_{l',l''}:Y_n(\sigma) \rightarrow Y_n(\sigma') \otimes
Y_n(\sigma'')$ lifting the maps $\Delta_{l',l''}$
for all $l' \geq s_{n,1}, l'' \geq s_{1,n}$.
Now we just need to show that this map
$\varprojlim \Delta_{l',l''}$ agrees with the restriction of the comultiplication $\Delta$ on $Y_n$.

Let $X$ be some generator $D_i^{(r)}, E_i^{(r)}$ or $F_i^{(r)}$
of $Y_n(\sigma) \subseteq Y_n$.
Consider how one computes $\Delta(X)$ in practice: first one expresses
$X$ as a linear combination of monomials in the $T_{i,j}^{(r)}$ as explained
in $\S$\ref{sparabolic} (see also the explicit formulae 
\cite[(5.2)--(5.4)]{BK}); then one computes $\Delta(X)$ in terms of
$T_{h,k}^{(s)}$'s
using the definition (\ref{comult}); finally
one rewrites all $T_{h,k}^{(s)}$ appearing in the resulting expression
back in terms of the generators $D_i^{(r)}, E_i^{(r)}$ and $F_i^{(r)}$.
But in view of Lemma~\ref{howithappens}, exactly the same procedure
may be used to compute the effect of $\Delta_{l',l''}$ on the image
of $X$ in $Y_{n,l}(\sigma)$ for each $l$. The point is that 
the formula expressing $T_{i,j;0}^{(r)}$
in terms of the elements $D_i^{(r)}, E_i^{(r)}$ and $F_i^{(r)}$
of $U(\mathfrak{p})$
explained in $\S$\ref{sinvariants} is identical to the formula
doing the same thing in $Y_n$ explained in $\S$\ref{sparabolic}.
It follows that $\Delta = \varprojlim \Delta_{l',l''}$.
\end{proof}

\begin{Remark}\label{allco}\rm
Using Theorem~\ref{grownup} and Lemma~\ref{cube}, one can 
express {\em all} of the comultiplications
$\Delta_{l',l''}:W(\pi)\rightarrow W(\pi') \otimes W(\pi'')$
in terms of the comultiplication $\Delta$ of $Y_n$, as follows.
Let $\dot\pi'$ denote the right-justified pyramid with the same
row lengths as $\pi'$, and let $\dot\pi''$ denote the left-justified
pyramid with the same row lengths as $\pi''$. 
Let $\dot \pi := \dot\pi' \otimes \dot\pi''$.
By Lemma~\ref{cube},
the comultiplication $\Delta_{l',l''}:W(\pi) \rightarrow W(\pi') \otimes W(\pi'')$ can be recovered from $\Delta_{l',l''}:W(\dot\pi)
\rightarrow W(\dot\pi') \otimes W(\dot\pi'')$ using the isomorphisms
$\iota$. Finally 
$\Delta_{l',l''}:W(\dot\pi)
\rightarrow W(\dot\pi') \otimes W(\dot\pi'')$ 
is one of the comultiplications described by Theorem~\ref{grownup}.
\end{Remark}

To end the main part of the article, we record one corollary of
Theorem~\ref{grownup}. This gives an application of shifted Yangians
to prove a result just involving the Yangian $Y_n$ itself.
For the statement, let $Y_{(1^n)}^{\sharp}$ resp.
$Y_{(1^n)}^{\flat}$ denote the Borel subalgebra of the Yangian $Y_n$
generated by the $D_i^{(r)}$'s together with the $E_i^{(r)}$'s resp. the
$F_i^{(r)}$'s.

\begin{Corollary}\label{provo}
The comultiplication $\Delta:Y_n\rightarrow Y_n\otimes Y_n$
has the following properties:
\begin{itemize}
\item[(i)]
$\Delta(D_i^{(r)}) \in Y_{(1^n)}^\sharp \otimes Y_{(1^n)}^\flat$;
\item[(ii)]
$\Delta(E_i^{(r)}) \in Y_{(1^n)}^\sharp \otimes Y_{(1^n)}^\flat + \C(1 \otimes E_i^{(r)})$.
\item[(iii)]
$\Delta(F_i^{(r)}) \in Y_{(1^n)}^\sharp \otimes Y_{(1^n)}^\flat + \C(F_i^{(r)} \otimes 1)$.
\end{itemize}
\end{Corollary}

\begin{proof}
(i) Let $\sigma = (s_{i,j})_{1 \leq i,j \leq n}$ be a shift matrix with 
$s_{j,j+1}, s_{j+1,j} \geq r$ for all $j=1,\dots,n-1$.
The element $D_i^{(r)}$ belongs to $Y_n(\sigma)$.
Hence by the first part of Theorem~\ref{grownup}, 
we know that $\Delta(D_i^{(r)})$ lies in the subalgebra
$Y_n(\sigma') \otimes Y_n(\sigma'')$ of $Y_n \otimes Y_n$. 
Moreover, $\Delta$ is a filtered map 
with respect to
the canonical filtation, hence in fact we have that
$\Delta(D_i^{(r)})  \in \sum_{s=0}^r \F_s Y_n(\sigma') \otimes \F_{r-s}
Y_n(\sigma'')$.
But by the choice of $\sigma$, $\F_s Y_n(\sigma') \subseteq Y_{(1^n)}^\sharp$
and $\F_{r-s} Y_n(\sigma'') \subseteq Y_{(1^n)}^\flat$ for all $s=0,\dots,r$,
so we are done.

(ii) Similar,
except one needs to take $s_{i,i+1}=r-1$ in the shift matrix $\sigma$.

(iii) Similar.
\end{proof}

\section{A special case}\label{sspecial}

In this section, we give a much more direct proof
of the main results of the article
in the special case that the nilpotent matrix $e$
consists of $n$ Jordan blocks all of the same size $l$, i.e. the
pyramid $\pi$ is an $n \times l$ rectangle and $N  = nl$.
The main theorem in this case was first noticed by Ragoucy and Sorba
\cite{RS}; see also \cite{BR} which is closer to the present exposition.
We will identify $\mathfrak{g} = \mathfrak{gl}_{N}$
with the tensor product 
$\mathfrak{gl}_l \otimes \mathfrak{gl}_n$
so that the matrix unit $e_{(r-1)n+i, (s-1)n+j} \in \mathfrak{gl}_{N}$
is identified with $e_{r,s} \otimes e_{i,j}
\in \mathfrak{gl}_l \otimes \mathfrak{gl}_n$
for $1 \leq r,s \leq l, 1 \leq i,j \leq n$.
Numbering the bricks of the pyramid $\pi$ down columns from left to right
as usual,
the $\mathfrak{sl}_2$-triple $(e,h,f)$ 
introduced at the beginning of $\S$\ref{sslice}
is given explicitly in this case by
$$
e = \sum_{r=1}^{l-1} e_{r,r+1} \otimes I_n,\quad
h = \sum_{r=1}^l (l+1-2r) e_{r,r} \otimes I_n,\quad
f = \sum_{r=1}^{l-1} r(l-r) e_{r+1,r} \otimes I_n.
$$
\iffalse
Also, 
$$
\mathfrak p = 
\left(
\begin{array}{c|c|c|c}
*&*&*&*\\
\hline
0&*&*&*\\
\hline
0&0&*&*\\
\hline
\vdots&\vdots&\vdots&\ddots
\end{array}
\right),
\qquad
\mathfrak h = 
\left(
\begin{array}{c|c|c|c}
*&0&0&0\\
\hline
0&*&0&0\\
\hline
0&0&*&0\\
\hline
\vdots&\vdots&\vdots&\ddots
\end{array}
\right),
\qquad
\mathfrak m = 
\left(
\begin{array}{c|c|c|c}
0&0&0&0\\
\hline
*&0&0&0\\
\hline
*&*&0&0\\
\hline
\vdots&\vdots&\vdots&\ddots
\end{array}
\right),
$$
viewing matrices in $\mathfrak g$ as
block matrices consisting of $l \times l$ blocks each of size $n \times n$.\fi
Let $\Mat_n$ denote the algebra of $n \times n$ matrices over $\C$
and let $T(\mathfrak{gl}_l)$ denote the tensor algebra on the vector
space $\mathfrak{gl}_l$.
Define an algebra homomorphism
\begin{equation}\label{t1}
T: T(\mathfrak{gl}_l) \rightarrow \Mat_n \otimes U(\mathfrak{gl}_{N})
\end{equation}
sending a generator $x \in \mathfrak{gl}_l$ to 
$T(x) = \sum_{i,j=1}^n 
e_{i,j} \otimes (x \otimes e_{i,j}) \in \Mat_n \otimes U(\mathfrak{gl}_{N})$.
Also define maps 
\begin{equation}\label{ijv}
T_{i,j}: T(\mathfrak{gl}_l) 
\rightarrow U(\mathfrak{gl}_{N})
\end{equation}
for each $1 \leq i,j \leq n$ from the equation
$T(x) = \sum_{i,j=1}^n e_{i,j} \otimes T_{i,j}(x)$,
for any $x \in T(\mathfrak{gl}_l)$.
Thus, thinking of $T(x)$ as an $n\times n$ matrix with entries in
$U(\mathfrak{gl}_{N})$, $T_{i,j}(x)$ is the $ij$-entry of $T(x)$.
We have by definition that $T(xy) =T(x) T(y)$
for any $x,y \in T(\mathfrak{gl}_l)$, which implies that
\begin{equation}
T_{i,j} (x y) = \sum_{k=1}^n T_{i,k}(x) T_{k,j}(y).
\label{split}
\end{equation}
In particular, $T_{i,j}(1) = \delta_{i,j}$ and
\begin{equation}\label{ijv3}
T_{i,j}(x_1 \otimes \cdots \otimes x_r)= \sum_{\substack{1 \leq 
i_0,i_1,\dots,i_r \leq n \\ i_0 = i, i_r = j}}\!\!\!\!\!
(x_1 \otimes e_{i_0,i_1})(x_2 \otimes e_{i_1, i_2}) 
\cdots (x_r \otimes e_{i_{r-1},i_r})
\end{equation}
for arbitrary $x_1,\dots,x_r \in \mathfrak{gl}_l$.
If $u$ is an indeterminate, we will also write 
$T_{i,j}$
for the map obtained from $T_{i,j}$ 
by extending scalars from $\C$ to $\C[u]$.

\begin{Lemma}\label{Easy}
For $1 \leq h,i,j,k \leq n$
and $x,y_1,\dots,y_r \in \mathfrak{gl}_l$,
\begin{align*}
[T_{i,j}(x), T_{h,k}(y_1\otimes\cdots\otimes y_r)]
=
&\sum_{s=1}^r 
T_{h,j}(y_1\otimes\cdots\otimes y_{s-1}) T_{i,k}(x y_s \otimes
y_{s+1}
\otimes \cdots \otimes y_r)
\\
-&\sum_{s=1}^r
T_{h,j}(y_1\otimes\cdots\otimes y_{s-1}\otimes y_{s} x) T_{i,k}(
y_{s+1}
\otimes \cdots \otimes y_r),
\end{align*}
where $xy_s$ and $y_sx$ denote the usual products of matrices  
in $\Mat_l$.
\end{Lemma}

\begin{proof}
Use (\ref{split}) and induction on $r$.
\end{proof}

Define $\Omega(u)$
to be the following $l \times l$ matrix 
with entries in $T(\mathfrak{gl}_l)[u]$:
\begin{equation*}
\Omega(u) := \left(
\begin{array}{cccccc}
e_{1,1}+u+\rho_1 & e_{1,2} & e_{1,3}&\cdots & e_{1,l}\\
1&e_{2,2}+u+\rho_2&&&\vdots\\
0&&\ddots&&e_{l-2,l}\\
\vdots&&1&e_{l-1,l-1} +u+\rho_{l-1}&e_{l-1,l}\\
0&\cdots&0&1&e_{l,l}  +u+\rho_l
\end{array}
\right)
\end{equation*}
where $\rho_r = -(l-r)n$ as in (\ref{rhodef}).
\ifcenters@
This matrix should be compared with
(\ref{znr}).
Recalling the definitions (\ref{detdefr}) and (\ref{ijv}), define 
\else
Define the row-wise determinant of an $l \times l$ matrix $A = (a_{i,j})$ 
with entries in a non-commutative ring $R$ by the formula
\begin{equation}\label{detdefr}
\rdet A  = \sum_{\pi \in S_l} \sgn(\pi) a_{1,\pi
1} \cdots a_{l,\pi l}.
\end{equation}
Then $\rdet \Omega(u)$ makes sense as an element of 
$T(\mathfrak{gl}_l)[u]$.
Now we use the maps (\ref{ijv}) to
define 
\fi
\begin{equation}
T_{i,j}(u) = \sum_{r = 0}^l T_{i,j}^{(r)} u^{l-r}
:= T_{i,j}(\rdet \Omega(u))
\end{equation}
for each $i,j=1,\dots,n$.
It is easy to see for $r = 0,1,\dots,l$ that
$T_{i,j}^{(r)}$ is precisely the element
$T_{i,j;0}^{(r)}$ of $\F_r U(\mathfrak{p})$ defined in (\ref{earlier}).
We also set $T_{i,j}^{(r)} := 0$ for $r > l$.
In the next lemma, we 
check directly that these elements 
belong to the algebra $W(\pi)= U(\mathfrak{p})^{\mathfrak{m}}$
introduced in $\S$\ref{sslice}.

\begin{Lemma}\label{inv}
For each $1 \leq i,j \leq n$ and $r > 0$,
$T_{i,j}^{(r)}$ is invariant under the twisted action
of $\mathfrak{m}$.
\end{Lemma}

\begin{proof}
Since $\mathfrak m$ is generated by elements of the form
$e_{r+1,r} \otimes e_{i,j}$, it suffices to show that 
$\pr_\chi ([T_{i,j}(e_{r+1,r}),T_{h,k}(\rdet \Omega(u))]) = 0$
for every $1 \leq h,i,j,k \leq n$ and 
$r = 1,\dots,l-1$.
Let us write $\Omega_{r, s}(u)$ for the submatrix of
$\Omega(u)$ consisting only of rows and columns 
numbered $r,\dots,s$.
We compute using Lemma~\ref{Easy} to get that
\begin{multline*}
[T_{i,j}(e_{r+1,r}), T_{h,k}(\rdet \Omega(u))]
=T_{h,j}(\rdet \Omega_{1, r-1}(u))\times\\
T_{i,k}
\left(
\rdet
\left(
\begin{array}{ccccc}
e_{r+1,r}&e_{r+1,r+1}&\cdots&e_{r+1,l}\\
1&e_{r+1,r+1}+u+\rho_{r+1} &\cdots & e_{r+1,l}\\\vdots&&\ddots&\vdots\\
0&\cdots&1&e_{l,l}+u+\rho_l
\end{array}
\right)\right)\\
-
T_{h,j}
\left(
\rdet 
 \left(
\begin{array}{ccccc}
e_{1,1}+u+\rho_1 &\cdots & e_{1,r}&e_{1,r}\\
1&\ddots&&\vdots\\
\vdots&&e_{r,r} +u+\rho_r&e_{r,r}\\
0&\cdots&1&e_{r+1,r}
\end{array}
\right)\right)\\
\times T_{i,k}(\rdet \Omega_{r+2, l}(u)).
\end{multline*}
In order to apply $\pr_\chi$ to the right hand side, we observe 
that for any $1 \leq m \leq n$,
$$
\pr_\chi \left(T_{i,m}\left(e_{r+1,r} \left(e_{r+1,r+1}+u +\rho_{r+1}\right)\right)\right)
=
T_{i,m}\left(e_{r+1,r+1} +u +\rho_r\right).
$$
Hence,
we get that
\begin{multline*}
\pr_\chi ([T_{i,j}(e_{r+1,r}), T_{h,k}(\rdet \Omega(u))])
=
T_{h,j}(\rdet \Omega_{1,r-1}(u))\times\\
T_{i,k}
\left(
\rdet
\left(
\begin{array}{ccccc}
1&e_{r+1,r+1}&\cdots&e_{r+1,l}\\
1&e_{r+1,r+1}+u +\rho_r &\cdots & e_{r+1,l}\\\vdots&&\ddots&\vdots\\
0&\cdots&1&e_{l,l}+u+\rho_l
\end{array}
\right)\right)\\
-
T_{h,j}
\left(
\rdet 
 \left(
\begin{array}{ccccc}
e_{1,1}+u+\rho_1 &\cdots & e_{1,r}&e_{1,r}\\
1&\ddots&&\vdots\\
\vdots&&e_{r,r}+u+\rho_r&e_{r,r}\\
0&\cdots&1&1
\end{array}
\right)\right)\\
\times T_{i,k} ( \rdet \Omega_{r+2,l}(u)).
\end{multline*}
Making the obvious row and column operations 
gives that
\begin{align*}
\rdet\left(
\begin{array}{ccccc}
1&e_{r+1,r+1}&\cdots&e_{r+1,l}\\
1&e_{r+1,r+1}+u +\rho_r &\cdots & e_{r+1,l}\\\vdots&&\ddots&\vdots\\
0&\cdots&1&e_{l,l}+u+\rho_l
\end{array}
\right)
&
=
\left(u+\rho_r\right) \rdet \Omega_{r+2,l}(u),\\
\rdet \left(
\begin{array}{ccccc}
e_{1,1}+u+\rho_1 &\cdots & e_{1,r}&e_{1,r}\\
1&\ddots&&\vdots\\
\vdots&&e_{r,r}+u+\rho_r&e_{r,r}\\
0&\cdots&1&1
\end{array}
\right)
&=\left(u+\rho_r\right) \rdet \Omega_{1,r-1}(u).
\end{align*}
Now substituting these into the above formula for
$\pr_\chi \left([T_{i,j}(e_{r+1,r}),T_{h,k}(\rdet \Omega(u))]\right)$
shows that it is zero.
\end{proof}

Denoting the matrix unit $e_{r,r} \otimes e_{i,j}$ instead by
$e_{i,j}^{[r]}$,
the subalgebra 
$\mathfrak h \cong \mathfrak{gl}_n^{\oplus l}$
of $\mathfrak{g}$ has basis
$\{e_{i,j}^{[r]}\}_{1 \leq r \leq l, 1 \leq i,j \leq n}$.
Let
\begin{equation}
\eta:U(\mathfrak h) \rightarrow U(\mathfrak h), \quad
e^{[r]}_{i,j} \mapsto e^{[r]}_{i,j} +
\delta_{i,j}(l-r)n.
\end{equation}
Let $\xi:U(\mathfrak p) \twoheadrightarrow U(\mathfrak h)$ be the
algebra homomorphism induced by the natural projection $\mathfrak p \twoheadrightarrow \mathfrak h$.
The composite $\mu := \eta \circ \xi$ is precisely the Miura
transform from (\ref{miura2}).

\begin{Lemma}\label{mim}
For $1 \leq i,j \leq n$ and $r > 0$, the Miura transform 
$\mu$ maps the element
$T_{i,j}^{(r)}$ of $U(\mathfrak{p})$ to the element
\begin{equation}\label{miuraimage}
\sum_{1 \leq s_1 < \cdots < s_r \leq l}
\sum_{\substack{1 \leq i_0,\cdots,i_r \leq n \\ i_0 = i, i_r = j}}
e_{i_0,i_1}^{[s_1]}
e_{i_1,i_2}^{[s_2]}
\cdots
e_{i_{r-1},i_r}^{[s_r]}
\end{equation}
of $U(\mathfrak{h})$.
\end{Lemma}

\begin{proof}
Applying the map $e_{r,s} \mapsto \delta_{r,s} (e_{r,r} + (l-r)n)$
to the matrix $\Omega(u)$ gives a diagonal matrix with
determinant $(u+e_{1,1}) (u+e_{2,2}) \cdots (u+e_{r,r})$.
The $u^{l-r}$-coefficient of this is the (non-commutative)
elementary symmetric function
$\sum_{1 \leq s_1 < \cdots < s_r \leq l} 
e_{s_1,s_1}
\cdots e_{s_r,s_r}$.
Now the lemma follows from (\ref{ijv3}).
\end{proof}

Finally 
let $Y_{n,l}$ denote the Yangian of level $l$ from \cite{Ch}.
This may be defined as the algebra on
generators $\{T_{i,j}^{(r)}\}_{1 \leq i,j \leq n, r >0}$
subject to the relations (\ref{mr}) 
and in addition
$T_{i,j}^{(r)} = 0$
for all $1 \leq i,j \leq n$ and $r > l$.
We will work now with the canonical filtration
on $Y_{n,l}$ and the Kazhdan filtration of 
$W(\pi)$, as defined in sections~\ref{scan}
and \ref{sslice} respectively.
The main theorem in our special case is as follows.

\begin{Theorem}
There is an isomorphism $Y_{n,l} \stackrel{\sim}{\rightarrow}
W(\pi)$ of filtered algebras such that the generators 
$\{T_{i,j}^{(r)}\}_{1 \leq i,j \leq n, r > 0}$
of $Y_{n,l}$ map to the elements of $W(\pi)$ with the same names.
Moreover, the Miura transform
$\mu:W(\pi) \rightarrow U(\mathfrak{h})$ is injective, and the map
$\theta$ from the diagram (\ref{thed})
is an isomorphism.
\end{Theorem}

\begin{proof}
By the PBW theorem for $Y_{n,l}$ from \cite[Theorem 3.1]{BK}
(which is due originally to Cherednik \cite{Ch0,Ch})
the set of all monomials in the elements
$\{T_{i,j}^{(r)}\}_{1 \leq i,j \leq n, 1 \leq r \leq l}$
taken in some fixed order and of total degree $\leq d$
form a basis for $\F_d Y_{n,l}$.
Moreover, {\em loc. cit.} shows that
there is an injective
algebra homomorphism
$\kappa_l:Y_{n,l} \hookrightarrow U(\mathfrak{h})$
mapping $T_{i,j}^{(r)} \in Y_{n,l}$ to precisely the
element (\ref{miuraimage}).
We also note from Lemma~\ref{centbase} that for each $d \geq 0$,
\begin{equation}\label{equaldims}
\dim \kappa_l(\F_d Y_{n,l}) =
\dim \F_d Y_{n,l}
= \dim \F_d S(\mathfrak{c}_{\mathfrak{g}}(e)),
\end{equation}
where $\F_d S(\mathfrak{c}_{\mathfrak{g}}(e))$
denotes the sum of all the graded pieces of 
$S(\mathfrak{c}_{\mathfrak{g}}(e))$ of degree $\leq d$
for the Kazhdan grading.

Let $X_d$ denote the subspace of $U(\mathfrak{p})$
spanned by all monomials in the elements
$\{T_{i,j}^{(r)}\}_{1 \leq i,j \leq n, 0 < r \leq l}$
taken in some fixed order and of total degree $\leq d$.
By Lemma~\ref{mim} and the previous paragraph, 
$\mu(X_d) = \kappa_l (\F_d Y_{n,l})$.
Hence, since we know by Lemma~\ref{inv} that $X_d \subseteq \F_d W(\pi)$, 
we get by 
(\ref{equaldims}) and
Corollary~\ref{phewy} that
$$
\dim \F_d S(\mathfrak{c}_{\mathfrak{g}}(e))
= \dim \mu(X_d) \leq \dim X_d \leq \dim \F_d W(\pi)
\leq \dim \F_d S(\mathfrak{c}_{\mathfrak{g}}(e)).
$$
Hence equality holds everywhere, which means
that $X_d = \F_d W(\pi)$, the Miura transform $\mu$ is injective
with the same image as $\kappa_l$, and the map $\theta$ is 
an isomorphism. Hence the composite $\mu^{-1} \circ \kappa_l$
gives the required isomorphism $Y_{n,l} \rightarrow W(\pi)$
of filtered algebras.
\end{proof}

\end{document}